\input{style/arxiv-general.cfg}
\documentclass[MSNbibl,citesort,number,seceqn,rotating,dvips]{arxbj}
\makeatletter
   \@ifpackageloaded{graphicx}{}{\usepackage{graphicx}}
\makeatother
\usepackage{upgreek,dcolumn}
\usepackage{graphicx}

\aid{0}
\volume{21}
\issue{2}
\pubyear{2015}
\firstpage{697}
\lastpage{739}
\doi{10.3150/13-BEJ584} 

\makeatletter

\newcommand{\rright}{\right}
\newcommand{\lleft}{\left}
\newremark{example}{Example}[section]
\newtheorem{theorem}{Theorem}[section]
\newtheorem{corollary}{Corollary}[section]
\newremark{remark}{Remark}[section]
\newtheorem{lemma}{Lemma}[section]
\newproclaim{assumption}{Assumption}[section]

\newcommand\un{\underline}
\newcommand{\N}{\mathbb{N}}
\newcommand{\R}{\mathbb{R}}
\newcommand{\Z}{\mathbb{Z}}
\newcommand{\C}{\mathbb{C}}
\newcommand{\tr}{\operatorname{tr}}
\renewcommand{\vec}{\operatorname{vec}}
\newcommand{\Var}{\operatorname{Var}}
\newcommand{\Cov}{\operatorname{Cov}}

\newcommand{\schwach}{\stackrel{\mathcal D}{\longrightarrow}}
\newcommand{\leb}{\lambda\! \! \lambda}

\newcommand{\eqref}[1]{(\ref{#1})}
\newcommand{\ARR}{\operatorname{AR}}
\newcommand{\AR}{\mathit{AR}}
\newcommand{\MA}{\mathit{MA}}
\newcommand{\MAA}{\operatorname{MA}}
\newcommand{\GARCH}{\mathit{GARCH}}
\newcommand{\GARCHA}{\operatorname{GARCH}}
\newcommand{\TAR}{\mathit{TAR}}
\newcommand{\TARA}{ \operatorname{TAR}}
\newcommand{\RCA}{\mathit{RCA}}
\newcommand{\RCAA}{\operatorname{RCA}}

\renewcommand{\emptyset}{\varnothing}
\newcommand{\ppi}{\uppi}
\newcommand{\pii}{\pi}

\def\sfrac#1#2{#1/#2}
\def\vfrac#1#2{(#1)/#2}

\newcolumntype{d}[1]{D{.}{.}{#1}}
\makeatother

\begin{document}
\begin{frontmatter}

\title{Testing equality of spectral densities using randomization techniques}
\runtitle{Testing equality of spectral densities}

\begin{aug}
\author[A]{\inits{C.}\fnms{Carsten} \snm{Jentsch}\corref{}\thanksref{A}\ead[label=e1]{cjentsch@mail.uni-mannheim.de}} \and
\author[B]{\inits{M.}\fnms{Markus} \snm{Pauly}\thanksref{B}\ead[label=e2]{markus.pauly@uni-duesseldorf.de}}
\address[A]{Department of Economics, University of Mannheim, L7, 3-5,
68131 Mannheim, Germany.\\ \printead{e1}}
\address[B]{Institute of Mathematics, University of D\"usseldorf,
Universit\"atsstrasse 1, 40225 D\"usseldorf, Germany. \printead{e2}}
\end{aug}

\received{\smonth{2} \syear{2012}}
\revised{\smonth{11} \syear{2013}}

%
\begin{abstract}
In this paper, we investigate the testing problem that the spectral
density matrices of several,
not necessarily independent, stationary processes are equal. Based on
an $L_2$-type test statistic,
we propose a new nonparametric approach, where the critical values of
the tests are calculated with
the help of randomization methods. We analyze asymptotic exactness and
consistency of these
randomization tests and show in simulation studies that the new
procedures posses very good size and power characteristics.
\end{abstract}

%

\end{frontmatter}

\section{Introduction and motivation}\label{sec1}

Suppose, one observes two stretches of time series data $X_1,\ldots
,X_n$ and $Y_1,\ldots,Y_n$ generated by some unknown stationary time
series processes. The question whether the models where these two
stretches come from are identical is of considerable interest and has
been investigated elaborately in statistical literature. For instance,
by imposing certain parametric assumptions on the data generation
process (DGP), it becomes possible to compare suitable estimated
parameters to judge how ``near'' the related models are. The assumption
of underlying autoregressive processes of some finite order, say,
allows to compare estimated coefficients to quantify the nearness of
both DGPs. But these models are usually not valid if the parametric
assumption is not fulfilled. Therefore, nonparametric methods are
desired that can be justified theoretically for a broader class of
stochastic processes.

If one is interested in second order properties, it seems obvious to
focus on autocovariances $\gamma(h)$, autocorrelations $\rho(h)$ or
spectral densities $f(\omega)$, where the latter choice has some
advantages over the others, see, for example, Paparoditis
\cite{Papa:2000}. Over
the years, several test statistics have been introduced that could be
used for testing the null hypothesis
%
\begin{eqnarray}
\label{H0special} \bigl\{ f_X(\omega)=f_Y(\omega) \mbox{
for all } \omega\in[-\ppi,\ppi] \bigr\}
\end{eqnarray}
versus the alternative $\{f_X\neq f_Y \mbox{ on a set of positive
$\leb$-measure}\}$. These may be perio\-dogram-based, see, for example,
Coates and Diggle \cite{Coates:1986}, P{\"o}tscher and Reschenhofer \cite{Poetscher:1988}, Diggle and Fisher
\cite{Diggle:1991}, Caiado \textit{et al.} \cite{Caiado:2006,Caiado:2009}
and Preu{ss} and Hildebrandt \cite{Preuss:2012}, or based on
kernel spectral density estimators, see, for example, Maharaj \cite{Maharaj:2002},
Eichler \cite{Eichler:2008} or Dette and
Paparoditis \cite{Dette:2009}. Jentsch and Pauly \cite{Jentsch/Pauly:2010} have investigated the
asymptotic properties of
test statistics closely related to ``periodogram-based distances''
proposed by Caiado \textit{et al.} \cite{Caiado:2009}.
They show that contrary to integrated periodograms these tests are
usually not consistent which is caused by the use of non-consistent
estimators for the spectral densities, the periodograms.
In contrast to that, $L_2$-type statistics based on the smoothed
periodogram are shown to be consistent, see Paparoditis \cite
{Papa:2000,Papa:2005},
Eichler \cite{Eichler:2008}, Dette and Paparoditis \cite{Dette:2009} and Jentsch \cite{Jentsch:2010}.
A possible test statistic for testing equality of spectral densities is
therefore a suitable inflated version of
%
\begin{eqnarray}
\label{L2TypForm} \int_{-\ppi}^\ppi\bigl(\widehat
f_X(\omega)-\widehat f_Y(\omega )\bigr)^2\,\mathrm{d}
\omega,
\end{eqnarray}
where $\widehat f_X$ and $\widehat f_Y$ are nonparametric kernel
estimators of the spectras.
Central limit theorems for such statistics can be found in the above
mentioned literature.
For instance, Paparoditis \cite{Papa:2000} has emphasized
that statistics of this
kind converge very slowly to a normal distribution, which results in a
poor performance of the corresponding asymptotic test, particularly for
small sample sizes.
Therefore, other approaches as bootstrap or resampling techniques in
general are needed to improve the small sample behavior of these
tests. Paparoditis \cite{Papa:2000,Papa:2005} has used a
(parametric) ARMA-bootstrap
for goodness-of-fit testing and Dette and Paparoditis \cite
{Dette:2009} have
proposed a different approach based on the asymptotic
Wishart-distribution of periodogram matrices.

In this paper, we consider the more general testing problem of
comparing the spectral density matrices of several, not necessarily
independent, stationary processes.
As motivated above, we will use an $L_2$-type statistic of the form
(\ref{L2TypForm}) which is adjusted for this multiple sample testing
problem. To overcome the slow convergence speed of our test statistic,
we propose a resampling method to compute critical values, the
randomization technique, that exploits the well-known asymptotic
independence of periodograms $I_n(\omega)$ and $I_n(\lambda)$ for
different frequencies $\omega\neq\lambda$ with $\omega,\lambda\in
[0,\ppi]$.
A nice feature of this approach is that only one tuning parameter, the
bandwidth of the involved kernel spectral densities, has to be selected.
This can be done automatically by selection procedures based on cross
validation as has been proposed by Beltr\~{a}o and
Bloomfield \cite{Beltrao:1987} and Robinson \cite
{Robinson:1991}. In comparison to other
mainstream resampling techniques as for example, the autoregressive
sieve or block bootstrap techniques, it has the advantage that no other
tuning parameter as order or block length has to be assessed in
addition to the bandwidth.
It is worth noting that these choices appear to be very crucial in
applications and the performance of the corresponding tests usually
reacts very sensitive on the choice of these quantities. Nevertheless,
these procedures are not applicable (at least) directly, because they
do not mimic the desired null distribution under the alternative automatically.
Moreover, the usage of randomization methods is very natural and has
also been proposed by Diggle and Fisher \cite{Diggle:1991} for
the special (and
included) testing problem of comparing the spectral densities of two
independent univariate stationary time series.
This idea has been adopted by Maharaj \cite
{Maharaj:2002} for comparing evolutionary
spectra of univariate nonstationary time series. In both papers, the
authors (only) account for the application of this procedure with the
help of simulations and the following nice heuristic.
Recall that a randomization test holds the prescribed level $\alpha$
exact for finite sample sizes if the data is invariant under the
corresponding randomization group, see, for example, Lehmann and Romano \cite{Lehmann:2005}. Since $I_{n,X}(\omega)$ and
$I_{n,Y}(\omega)$ are
asymptotically exchangeable in the case of independent time series, one
would suggest that a randomization test, say $\varphi_n^*$, for (\ref
{H0special}) based on periodograms or functions of periodograms is at
least asymptotically exact, that is, $E(\varphi_n^*)\rightarrow\alpha
$ holds as $n\rightarrow\infty$.
However, this is only a conjecture and so far the randomization
technique has not been analyzed in this context from a mathematical
point of view in the literature. This is without question meaningful
since it is offhand not clear in which situations it is applicable.
We will close this gap in the next sections by investigating more
general and in mathematical detail whether at all and, if true, under
which possibly needed additional assumptions randomization-based
testing procedures lead to success, that is, asymptotically exact and
consistent tests.
Amongst others it will turn out that this can be done effectively for
testing the special null hypothesis~(\ref{H0special}).\looseness=1

The paper is organized as follows. In the next Section~\ref
{subsecnotation}, we introduce our notation, state the general testing
problem and define the $L_2$-type test statistic. Section~\ref{sec2}
presents the main assumptions on the process and states asymptotic
results for the test statistic which lead to an asymptotically exact
and consistent benchmark test. The randomization procedure together
with a corresponding discussion about its consistency is described in
Section~\ref{sec3}. Based on our theoretical results we compare the
proposed tests in the special situation of testing for (\ref
{H0special}) in an extensive simulation study in Section~\ref{sec5},
where the effects of dependence and bandwidth choice on the size and
power behavior of both types of tests are investigated. All proofs are
deferred to Section~\ref{sec6}.

%
\subsection{Notation and formulation of the problem}\label
{subsecnotation}

We consider a $d$-dimensional zero mean stationary process $(\un
{X}_t,t\in\Z)$ with $d=pq$, $p,q\in\N$ and $q\geq2$. Under
suitable assumptions, the process $(\un{X}_t,t\in\Z)$ posses a
continuous $(d\times d)$ spectral density matrix $\mathbf{f}$ given by
\begin{eqnarray*}
\mathbf{f}(\omega):= \frac{1}{2\ppi} \sum_{h=-\infty}^\infty
\boldsymbol{\Gamma}(h) \exp(-\mathrm{i}h\omega),\qquad  \omega\in[-\ppi,\ppi],
\end{eqnarray*}
where $\boldsymbol{\Gamma}(h):=E(\un{X}_{t+h}\un{X}_t^T),h\in\Z$, are
the corresponding autocovariance matrices. To introduce our testing
problem of interest, we write
\begin{eqnarray*}
\mathbf{f}(\omega)&=&\lleft( %
\begin{array} {c@{\quad}c@{\quad}c@{\quad}c}f_{11}(
\omega) & f_{12}(\omega) & \cdots& f_{1d}(\omega)
\\
f_{21}(\omega) & f_{22}(\omega) & & f_{2d}(\omega
)
\\
\vdots& & \ddots& \vdots
\\
f_{d1}(\omega) & f_{d2}(\omega) & \cdots& f_{dd}(
\omega) \end{array} %
\rright)\\
&=&\lleft( %
\begin{array}
{c@{\quad}c@{\quad}c@{\quad}c}\mathbf{F}_{11}(\omega) & \mathbf{F}_{12}(\omega ) &
\cdots& \mathbf{F}_{1q}(\omega)
\\
\mathbf{F}_{21}(\omega) & \mathbf{F}_{22}(\omega) & &
\mathbf{F}_{2q}(\omega)
\\
\vdots& & \ddots& \vdots
\\
\mathbf{F}_{q1}(\omega) & \mathbf {F}_{q2}(\omega) & \cdots&
\mathbf{F}_{qq}(\omega) \end{array} %
\rright),\qquad  \omega
\in[-\ppi,\ppi],
\end{eqnarray*}
where $\mathbf{F}_{mn}(\omega)$ are $(p\times p)$ block matrices with
$\mathbf{F}_{mn}(\omega)=\overline{\mathbf{F}_{nm}(\omega)}^T$ for
all $m,n$. Here and throughout the paper, all matrix-valued quantities
are written as bold letters, all vector-valued quantities are
underlined, $\overline{\mathbf{A}}$ denotes the complex conjugation
and $\mathbf{A}^T$ the transposition of a complex matrix $\mathbf
{A}$.

In this general setup, we want to test whether all $(p\times p)$
spectral density matrices $\mathbf{F}_{kk}$ of the $q$ (sub-)processes
$(\un{X}_{t,k},t\in\Z)$, $k=1,\ldots,q$ with
\begin{eqnarray*}
\un{X}_{t,k}=(X_{t,(k-1)p+1},\ldots,X_{t,kp})^T
\end{eqnarray*}
are identical. Precisely, suppose we have observations $\un
{X}_1,\ldots,\un{X}_n$ at hand and we want to test the null hypothesis
%
\begin{eqnarray}
\label{H0} H_0\dvt  \bigl\{ \mathbf{F}_{11}(\omega)=
\mathbf{F}_{22}(\omega)=\cdots =\mathbf{F}_{qq}(\omega) \mbox{
for all } \omega\in[-\ppi,\ppi]\bigr\}
\end{eqnarray}
against the alternative
\begin{eqnarray*}
H_1\dvt  \bigl\{\exists k_1,k_2\in\{1,\ldots,q
\}, A\subset\mathcal{B} \mbox { with } \leb(A)>0 \mbox{ such that }
\mathbf{F}_{k_1k_1}(\omega )\neq\mathbf{F}_{k_2k_2}(\omega)\ \forall
\omega\in A\bigr\},
\end{eqnarray*}
where $\leb$ denotes the one-dimensional Lebesgue-measure on the Borel
$\sigma$-algebra $\mathcal{B}$. Observe that this general framework
includes \eqref{H0special}.

Consider now the periodogram matrix $\mathbf{I}(\omega):= \un
{J}(\omega)\overline{\un{J}(\omega)}^T$ based on $\un{X}_1,\ldots
,\un{X}_n$, where
%
\begin{eqnarray}
\label{DFT} \un{J}(\omega):= \frac{1}{\sqrt{2\ppi n}} \sum
_{t=1}^n \un{X}_t \mathrm{e}^{-\mathrm{i}t\omega},\qquad
\omega\in[-\ppi,\ppi]
\end{eqnarray}
is the corresponding $d$-variate discrete Fourier transform (DFT). For
a better comparison with the spectral density matrix $\mathbf{f}$, we
write as above
\begin{eqnarray*}
\mathbf{I}(\omega)&=&\lleft( %
\begin{array} {c@{\quad}c@{\quad}c@{\quad}c}I_{11}(
\omega) & I_{12}(\omega) & \cdots& I_{1d}(\omega)
\\
I_{21}(\omega) & I_{22}(\omega) & & I_{2d}(\omega
)
\\
\vdots& & \ddots& \vdots
\\
I_{d1}(\omega) & I_{d2}(\omega) & \cdots& I_{dd}(
\omega) \end{array} %
\rright)\\
&=&\lleft( %
\begin{array}
{c@{\quad}c@{\quad}c@{\quad}c}\mathbf{I}_{11}(\omega) & \mathbf{I}_{12}(\omega ) &
\cdots& \mathbf{I}_{1q}(\omega)
\\
\mathbf{I}_{21}(\omega) & \mathbf{I}_{22}(\omega) & &
\mathbf{I}_{2q}(\omega)
\\
\vdots& & \ddots& \vdots
\\
\mathbf{I}_{q1}(\omega) & \mathbf {I}_{q2}(\omega) & \cdots&
\mathbf{I}_{qq}(\omega) \end{array} %
\rright),\qquad  \omega
\in[-\ppi,\ppi],
\end{eqnarray*}
where $\mathbf{I}_{mn}(\omega)$ are $(p\times p)$ block matrices with
$\mathbf{I}_{mn}(\omega)=\overline{\mathbf{I}_{nm}(\omega)}^T$ for
all $m,n$. Moreover, we define its pooled block diagonal periodogram
matrix by
\begin{eqnarray*}
\widetilde{\mathbf{ I}}(\omega_k):=\frac{1}{q}\sum
_{j=1}^q \mathbf {I}_{jj}(
\omega_k).
\end{eqnarray*}
With that we can introduce the kernel estimator
\begin{eqnarray*}
\widehat{\mathbf{f}}(\omega) := \frac{1}{n} \sum
_{k=-\lfloor\vfrac
{n-1}{2}\rfloor}^{\lfloor\sfrac{n}{2}\rfloor} K_h(\omega-
\omega_k) \mathbf{I}(\omega_k),\qquad  \omega\in[-\ppi,\ppi]
\end{eqnarray*}
for the spectral density matrix $\mathbf{f}(\omega)$, where $\lfloor
x\rfloor$ is the integer part of $x\in\R$, $\omega_k:= 2\ppi\frac
{k}{n}, k=-\lfloor\frac{n-1}{2}\rfloor,\ldots,\lfloor\frac
{n}{2}\rfloor$ are the Fourier frequencies, $K$ is a nonnegative
symmetric kernel function satisfying $\int K(x)\,\mathrm{d}x=2\ppi$, $h$ is the
bandwidth and $K_h(\cdot):=\frac{1}{h}K(\frac{\cdot}{h})$. Note
that the periodogram is understood as a $2\ppi$-periodically extended
function on the real line throughout the paper. Moreover, we write as above
\begin{eqnarray*}
\widehat{\mathbf{f}}(\omega)&=&\lleft( %
\begin{array} {c@{\quad}c@{\quad}c@{\quad}c}
\widehat{f}_{11}(\omega) & \widehat {f}_{12}(\omega) & \cdots&
\widehat{f}_{1d}(\omega)
\\[2pt]
\widehat {f}_{21}(\omega) & \widehat{f}_{22}(\omega) & &
\widehat {f}_{2d}(\omega)
\\
\vdots& & \ddots& \vdots
\\
\widehat {f}_{d1}(\omega) & \widehat{f}_{d2}(\omega) &
\cdots& \widehat {f}_{dd}(\omega) \end{array} %
\rright) \\
&=&
\lleft( %
\begin{array} {c@{\quad}c@{\quad}c@{\quad}c} \widehat{\mathbf{F}}_{11}(
\omega) & \widehat {\mathbf{F}}_{12}(\omega) & \cdots& \widehat{\mathbf
{F}}_{1q}(\omega)
\\
\widehat{\mathbf{F}}_{21}(\omega) & \widehat {\mathbf{F}}_{22}(
\omega) & & \widehat{\mathbf{F}}_{2q}(\omega)
\\
\vdots& & \ddots& \vdots
\\
\widehat{\mathbf{F}}_{q1}(\omega) & \widehat{\mathbf{F}}_{q2}(
\omega) & \cdots& \widehat{\mathbf {F}}_{qq}(\omega) \end{array}
\rright),\qquad  \omega\in[-\ppi,\ppi],
\end{eqnarray*}
where $\widehat{\mathbf{F}}_{mn}(\omega)$ are uniformly consistent
estimators of the $(p\times p)$ block matrices $\mathbf{F}_{mn}(\omega
)$ for all $m,n$, that is, convergence $\sup_\omega\|\mathbf
{\widetilde F}_{mn}(\omega)-\mathbf{F}_{mn}(\omega)\|\rightarrow0$
holds almost surely.
Here and throughout the paper, $\|\cdot\|$ denotes the
Frobenius norm, that is, for a matrix $\mathbf{A}=(a_{ij})_{1\leq
i,j\leq p}\in\C^{p\times p}$ we set
\begin{eqnarray*}
\|\mathbf{A}\|^2 := \tr\bigl(\mathbf{A}\overline{
\mathbf{A}}^T\bigr) = \sum_{i,j=1}^p
|a_{ij}|^2,
\end{eqnarray*}
where $\tr(\cdot)$ stands for the trace of a matrix, i.e. $\tr
(\mathbf{A}):=\sum_{i=1}^p a_{ii}$.

For testing $H_0$, we now propose the $L_2$-type test statistic
%
\begin{eqnarray}
\label{teststat} T_n &:=& nh^{\sfrac{1}{2}}\int_{-\ppi}^\ppi
\sum_{r=1}^q \Biggl\|\frac{1}{n}\sum
_{k=-\lfloor\vfrac{n-1}{2}\rfloor}^{\lfloor
\sfrac{n}{2}\rfloor} K_h(\omega-
\omega_k) \bigl(\mathbf{I}_{rr}(\omega_k)-
\widetilde{ \mathbf{I}}(\omega_k)\bigr) \Biggr\|^2\,\mathrm{d}\omega\nonumber
\\[-8pt]\\[-8pt]
&=& nh^{\sfrac{1}{2}}\int_{-\ppi}^\ppi\sum
_{r=1}^q \bigl\|\widehat {\mathbf{F}}_{rr}(
\omega) - \widetilde{\mathbf{F}}(\omega) \bigr\| ^2 \,\mathrm{d}\omega,
\nonumber
\end{eqnarray}
where $\widetilde{\mathbf{F}}(\omega) := \frac{1}{q} \sum_{j=1}^q
\widehat{\mathbf{F}}_{jj}(\omega)$ is the pooled spectral density
estimator.

For a better understanding of the asymptotic results in the following
sections, it is interesting to note that
\begin{eqnarray*}
\mathbf{I}_{rr}(\omega_k)-\widetilde{\mathbf{I}}(
\omega_k)=-\frac
{1}{q}\sum_{j=1}^q
(1-q\delta_{jr}) \mathbf{I}_{jj}(\omega_k)
\end{eqnarray*}
holds, where $\delta_{jr}=1$ if $j=r$ and $\delta_{jr}=0$ otherwise.

\section{The unconditional test}\label{sec2}

For deriving an asymptotically exact test for (\ref{H0}), we will in
the sequel apply a result of Eichler \cite
{Eichler:2008} for a general class of
$L_2$-type statistics that includes the form (\ref{teststat}).
Therefore, we have to impose some assumptions.

%
\subsection{Assumptions}

First, we assume some strong mixing conditions for the DGP (compare
Assumption~3.1 in Eichler \cite{Eichler:2008}).

%
\begin{assumption}\label{ass1}
$(\un{X}_t,t\in\Z)$ is a zero mean $d$-variate (strictly)  stationary stochastic
process defined on a probability space $(\Omega,\mathcal{F},\mathbb
{P})$. Furthermore, for any $k>0$, the $k$th order cumulants of $(\un
{X}_t,t\in\Z)$ satisfy the mixing conditions
%
\begin{eqnarray}
\label{mixingcond} \sum_{u_1,\ldots,u_{k-1}\in\Z} \bigl(1+|u_j|^2
\bigr)\bigl|c_{a_1,\ldots
,a_k}(u_1,\ldots,u_{k-1})\bigr|< \infty
\end{eqnarray}
for all $j=1,\ldots,k-1$ and $a_1,\ldots,a_k=1,\ldots,d$, where
\begin{eqnarray*}
c_{a_1,\ldots,a_k}(u_1,\ldots,u_{k-1})=\operatorname{cum}(X_{u_1,a_1},
\ldots ,X_{u_{k-1},a_{k-1}},X_{0,a_k})
\end{eqnarray*}
is the $k$th order joint cumulant of $X_{u_1,a_1},\ldots
,X_{u_{k-1},a_{k-1}},X_{0,a_k}$ (cf. Brillinger
\cite{Brillinger:1981} for the definition).
\end{assumption}

Note that the above assumption requires the existence of the moments of
all orders. However, in the case that we presume that our process has a
linear structure, it can be weakened. See Remark~\ref{linear
processes} below.

Our second assumption on the kernel function $K$ and the bandwidth $h$
ensures the consistency of the kernel density estimators. It is similar
to and implies Assumption~3.3 of Eichler \cite
{Eichler:2008}.

\begin{assumption}\label{ass2}
\begin{enumerate}[(ii)]
\item[(i)] The kernel $K$ is a bounded, symmetric, nonnegative and
Lipschitz-continuous
function with compact support $[-\ppi,\ppi]$ and
\begin{eqnarray*}
\int_{-\ppi}^\ppi K(\omega)\,\mathrm{d}\omega=2\ppi.
\end{eqnarray*}
Furthermore, $K(\omega)$ has continuous Fourier transform $k(u)$ such that
%
\begin{eqnarray}
\label{Fouriertransform} \int k^2(u)\,\mathrm{d}u<\infty\quad \mbox{and}\quad  \int
k^4(u)\,\mathrm{d}u<\infty.
\end{eqnarray}
\item[(ii)] The bandwidth $h=h(n)$ is such that $h^{9/2}n\rightarrow
0$ and $h^2n\rightarrow\infty$ as $n\rightarrow\infty$.
\end{enumerate}
\end{assumption}

In comparison to the notation in Eichler \cite
{Eichler:2008}, remark that a factor
$2\ppi$ is incorporated in $K$. Moreover, for a better lucidity we
define the positive constants
%
\begin{eqnarray}\label{constants}
A_K:=\frac{1}{2\ppi}\int_{-\ppi}^\ppi
K^2(v)\,\mathrm{d}v \quad \mbox{and}\quad  B_K:=\frac{1}{\ppi^2}\int
_{-2\ppi}^{2\ppi} \biggl(\int_{-\ppi}^\ppi
K(v)K(v+z)\,\mathrm{d}v \biggr)^2\,\mathrm{d}z.
\end{eqnarray}
Furthermore, we like to point out that the optimal rate
of $h\approx n^{-1/5}$ (by means of an averaged mean integrated squared
error type criterion) for estimating the spectral density
is not covered by our assumptions.
However, assumptions of such kind are commonly imposed in the
literature in order to reduce the bias of the smoothed kernel estimator
leading to a central limit theorem for $T_n$, see Theorem~\ref
{Theorem 21} below.
Compare for instance Taniguchi and Kondo \cite
{Taniguchi:1993}, Taniguchi \textit{et al.}
\cite{Taniguchi:1996}, Taniguchi and Kakizawa \cite{Taniguchi:2000},
Eichler \cite{Eichler:2008} or
Dette and Paparoditis \cite{Dette:2009}, who used similar
(non-optimal)
assumptions on the rate of the bandwidth.

\subsection{Asymptotic results for $T_n$}

We are now ready to state a CLT for the test statistic $T_n$ given in
(\ref{teststat}).

%
\begin{theorem}[(Asymptotic null distribution of $T_n$)]\label
{theorem1}\label{Theorem 21}
Suppose that Assumptions \ref{ass1} and \ref{ass2} are fulfilled. If
$H_0$ is true, it holds
\begin{eqnarray*}
T_n-\frac{\mu_0}{\sqrt{h}}\stackrel{\mathcal{D}} {\longrightarrow }
\mathcal{N}\bigl(0,\tau_0^2\bigr)
\end{eqnarray*}
as $n\rightarrow\infty$, where
%
\begin{eqnarray}
\label{mu} \mu_0= A_K \int_{-\ppi}^\ppi
\Biggl(\frac{1}{q}\sum_{j_1,j_2=1}^q(-1+q
\delta_{j_1j_2})\bigl|\tr\bigl(\mathbf{F}_{j_1j_2}(\omega )
\bigr)\bigr|^2 \Biggr)\,\mathrm{d}\omega
\end{eqnarray}
and
%
\begin{eqnarray}\label{eq7}
\label{tau} \tau_0^2 &=& B_K \int
_{-\ppi}^\ppi \Biggl(\frac{1}{q^2}\sum
_{j_1,j_2,j_3,j_4=1}^q (-1+q\delta_{j_1j_2}) (-1+q
\delta_{j_3j_4})\nonumber\\[-8pt]\\[-8pt]
&&\hphantom{B_K \int
_{-\ppi}^\ppi \Biggl(\frac{1}{q^2}\sum
_{j_1,j_2,j_3,j_4=1}^q}{}\times\bigl|\tr\bigl(\mathbf {F}_{j_1j_3}(\omega)\overline{
\mathbf{F}_{j_2j_4}(\omega )}^T\bigr)\bigr|^2 \Biggr)\,\mathrm{d}
\omega, \nonumber
\end{eqnarray}
where $\mathbf{F}_{jj}(\omega)=\mathbf{F}_{11}(\omega)$ holds for
all $\omega$ and $j=1,\ldots,q$. Here and throughout the paper,
$\stackrel{\mathcal{D}}{\longrightarrow}$ denotes convergence in
distribution.
\end{theorem}

Remark that we have chosen the above representation of $\mu_0$ and
$\tau^2_0$ for notational convenience. It is of course possible to
split the sums above and to rewrite the expressions by using $\mathbf
{F}_{jj}=\mathbf{F}_{11}$.

It is interesting to note that we do not need the assumption of
positive definite spectral density block matrices since we do not work
with ratio-type test statistics, see, for example, Eichler \cite{Eichler:2008} for
$p=1$ and $q=2$ or Dette and Paparoditis \cite{Dette:2009}
for the case $p=1$.

We give one important example which is possibly of most relevant
interest, where we compare two real-valued spectral densities.

%
\begin{example}[(The case $p=1$ and $q=2$)]\label{exH0special}
For $p=1$ and $q=2$, the quantities $\mu_0$ and $\tau^2_0$ defined in
Theorem~\ref{theorem1} become
\begin{eqnarray*}
\mu_0= A_K \int_{-\ppi}^\ppi
f_{11}^2(\omega) \bigl(1-C_{12}(\omega )\bigr)\,\mathrm{d}
\omega
\end{eqnarray*}
and
\begin{eqnarray*}
\tau^2_0=B_K \int_{-\ppi}^\ppi
f_{11}^4(\omega) \bigl(1-C_{12}(\omega )
\bigr)^2\,\mathrm{d}\omega,
\end{eqnarray*}
where $C_{jk}(\omega)=\frac{|f_{jk}(\omega)|^2}{f_{jj}(\omega
)f_{kk}(\omega)}1(f_{jj}(\omega)f_{kk}(\omega)>0)$ is the squared
coherence between the two components $j$ and $k$ of $(\un{X}_t,t\in\Z
)$ (cf. Hannan \cite{Hannan:1970}, page 43).
\end{example}

Since $\mu_0$ and $\tau^2_0$ are in general unknown we have to
estimate them before we can apply the above results for testing (\ref{H0}).
Following Eichler \cite{Eichler:2008}, Remark~3.7, we can estimate
both quantities
by substituting $\mathbf{F}_{jk}$ in their expressions by the
corresponding consistent
estimators $\widehat{\mathbf{F}}_{jk}$.
Rewriting the right-hand side of the equations (\ref{mu}) and (\ref
{tau}) under $H_0$ by means of the equality $\mathbf{F}_{jj}= \frac
{1}{q} \sum_{i=1}^q\mathbf{F}_{ii}$ this leads to
%
\begin{eqnarray}\label{muwidehat}
\label{muHut} \widehat{\mu}:= A_K\int_{-\ppi}^\ppi
\biggl((q-1)\bigl|\tr\bigl(\widetilde{\mathbf{F}}(\omega)\bigr)\bigr|^2 -
\frac{1}{q} \mathop{\mathop{\sum}_{j_1,j_2=1}}_{j_1\neq j_2}^q
\bigl|\tr\bigl(\widehat{\mathbf{F}}_{j_1j_2}(\omega)\bigr)\bigr|^2
\biggr)\,\mathrm{d}\omega
\end{eqnarray}
as an adequate estimator for $\mu_0$. Similarly, we construct a
consistent estimator $\widehat{\tau}^2$ for $\tau^2_0$.

Note that it is in applications more convenient to use a discretized
version of $\widehat{\mu}$, where the last integral is approximated
by its Riemann sum.

Now we are ready to define the unconditional test
\begin{eqnarray*}
\varphi_n:=\mathbf{1}_{(u_{1-\alpha},\infty)}\bigl(\bigl(T_n-h^{-1/2}
\widehat {\mu}\bigr)/\widehat{\tau}\bigr),
\end{eqnarray*}
where $u_{1-\alpha}$ denotes the $(1-\alpha)$-quantile of the
standard normal distribution. Its properties are summarized in the
following theorem.

\begin{theorem}\label{theorem11}
Suppose that the Assumptions \ref{ass1} and \ref{ass2} hold. Then the
test $\varphi_n$ is not only of asymptotic level $\alpha$, that is,
$E(\varphi_n)\rightarrow\alpha$ holds if $H_0$ is true, but also
consistent for testing $H_0$ versus $H_1$, that is, $E(\varphi
_n)\rightarrow1$ under $H_1$ as $n\rightarrow\infty$.
\end{theorem}

\subsection{Power under local alternatives}

When studying the behavior of the test $\varphi_n$ under local
alternatives we have to consider observations $\underline{X}_1^n,\dots
,\underline{X}_n^n $ that come from sequences of $d$-dimensional
zero-mean processes $(\underline{X}_t^n, t\in\Z)$ with spectral
density matrices given by
%
\begin{eqnarray}
\label{localAlternatives} \mathbf{f}^n (\omega) = \mathbf{f} (\omega) +
\alpha_n \mathbf{g} (\omega).
\end{eqnarray}
Following Eichler \cite{Eichler:2008} we suppose
that $\mathbf{f}$ is nonnegative
definite, Hermitian and satisfies $H_0$ and that $\alpha_n\rightarrow
0$ as $n \rightarrow\infty$. Moreover, $\mathbf{f}$ and $\mathbf
{g}$ are assumed to be
twice continuously differentiable.
In addition, for every $n\in\N$ the sequence of processes
$(\underline{X}_t^n, t\in\Z)$ satisfies the strong mixing condition
(\ref{mixingcond}) uniformly in $n$,
where we replace
$c_{a_1,\ldots,a_k}(u_1,\ldots,u_{k-1})$ by $\operatorname{cum}(X_{u_1,a_1}^n,\ldots
,X_{u_{k-1},a_{k-1}}^n,X_{0,a_k}^n)$
(the $k$th order joint cumulant of $X_{u_1,a_1}^n,\ldots
,X_{u_{k-1},a_{k-1}}^n,\allowbreak X_{0,a_k}^n$).
Under these assumptions it can be deduced from Theorem $5.4$ in
Eichler \cite{Eichler:2008} that the test $\varphi
_n$ can detect local alternatives up to
an order of $\alpha_n = h^{-1/4}n^{-1/2}$. To be concrete, we have the
following result.

%
\begin{theorem}
\label{theorem: loclaAlternatives}
Under the conditions stated above the asymptotic power of the test
$\varphi_n$ under local alternatives (\ref{localAlternatives}) with
$\alpha_n = h^{-1/4}n^{-1/2}$ is given by
$
\lim_{n\rightarrow\infty} E(\varphi_n) = 1- \Phi(u_{1-\alpha} -
\nu/\tau_0)$,
where $\Phi$ denotes the c.d.f. of the standard normal distribution and
the detection shift $\nu$ is defined as
\[
\nu= \int_{-\ppi}^\ppi\vec\bigl(\overline{\mathbf{g}(
\omega)}\bigr)^T \boldsymbol{\Gamma} \vec\bigl(\mathbf{g}(\omega)
\bigr) \,\mathrm{d}\omega.
\]
Here the $(d^2\times d^2$) matrix $\boldsymbol{\Gamma} = (\Gamma
_{r,s})_{r,s}$ has the entries
\[
\Gamma_{i+(j-1)d, k + (\ell-1)d} = %
\cases{ \displaystyle \frac{q-1}{q} &\quad  \mbox{for}
$|i-k|=|j-l|=0$ \mbox{and} $(i,j)\in\Xi$,
\vspace*{4pt}\cr
-\displaystyle \frac{1}{q} & \quad \mbox{for}
$|i-k|=|j-l|\in p\N$ \mbox{and} $(i,j)\in \Xi$,\vspace*{2pt}
\cr
0 & \quad \mbox{otherwise},}
\]
with $\Xi= \bigcup_{c=1}^q \{(r,s)\dvt  1+ (c-1)p \leq r,s\leq cp \} $.
\end{theorem}

Note, that the same detection order has also been found by Paparoditis \cite{Papa:2000}.

%
\section{The randomization tests}\label{sec3}

Although the asymptotic performance of the unconditional test proposed
above seems to be satisfactory, it is well known that the convergence
speed of $L_2$-type statistics like (\ref{teststat}) is rather slow
(see H{\"a}rdle and Mammen \cite{Haerdle/Mammen:1993} as
well as Paparoditis (\cite{Papa:2000,Papa:2005}) for
more details). For small or moderate sample sizes, we therefore expect
that the above Gaussian approximation will not perform pretty well in
applications. Our simulations in Section~\ref{sec5} support this presumption.
To overcome this hitch, we propose in the following a suitable
randomization technique for approximating the finite sample
distribution of our test statistic. This method leads to so called
randomization tests that use more appropriate, data-dependent critical
values.
Let $\pi_k\dvtx  (\widetilde{\Omega},\widetilde{\mathcal{F}},\widetilde
{\mathbb{P}}) \rightarrow\mathcal{S}_q,k=0,\ldots,\lfloor{\frac
{n}{2}}\rfloor$, with
%
\begin{eqnarray}
\label{random permutation defi} \pii_k=\bigl(\pii_k(1),\ldots,
\pii_k(q)\bigr)^T
\end{eqnarray}
be a sequence of\vspace*{2pt} independent and uniformly distributed random
variables on the symmetric group $\mathcal{S}_q$ (the set of all
permutations of $(1,\dots,q)$) defined on some further probability
space $(\widetilde{\Omega},\widetilde{\mathcal{F}},\widetilde
{\mathbb{P}})$. In what follows, we assume that $(\pii_k)_k$ and $(\un
{X}_t)_t$ are independent (defined as random variables on the joint
probability space $(\Omega\times\widetilde{\Omega},\mathcal
{F}\otimes\widetilde{\mathcal{F}},\mathbb{P}\otimes\widetilde
{\mathbb{P}})$).
Conditional on the observations $\un{X}_1,\ldots,\un{X}_n$ the
corresponding randomization test statistic is given by
\begin{eqnarray*}
T_n^*=nh^{\sfrac{1}{2}}\int_{-\ppi}^\ppi
\sum_{r=1}^q \Biggl\|\frac
{1}{n}\sum
_{k=-\lfloor\vfrac{n-1}{2}\rfloor}^{\lfloor\sfrac
{n}{2}\rfloor} K_h(\omega-
\omega_k) \bigl(\mathbf{I}_{\pii_k(r),\pii
_k(r)}(\omega_k)-
\widetilde{ \mathbf{I}}(\omega_k)\bigr) \Biggr\|^2\,\mathrm{d}\omega,
\end{eqnarray*}
where we set $\pii_{-k}:=\pii_k$ and $\pii_{k+sn}:=\pii_k$ for $s\in\Z
$ to maintain the symmetry properties of a spectral density matrix.
Note that similar to the unconditional case, it holds
\begin{eqnarray*}
\mathbf{I}_{\pii_k(r),\pii_k(r)}(\omega_k)-\widetilde{\mathbf{ I}}(
\omega_k) = -\frac{1}{q}\sum_{j=1}^q
(1-q\delta_{j,\pii_k(r)}) \mathbf {I}_{jj}(\omega_k).
\end{eqnarray*}

\subsection{Asymptotic results for $T_n^*$}

In the following, we will analyze in which situations our randomization
method leads to asymptotically valid tests, that is, whether some
suitable centered $T_n^*$ converges to the same distribution as
$T_n-h^{-1/2}\mu_0$.
Therefore, we have to exploit the limiting behavior of the
randomization statistic $T_n^*$. It will turn out that we have to
center $T_n^*$ at
%
\begin{eqnarray}
\label{eq61} \widehat{\mu}^* &:=& A_K\int_{-\ppi}^\ppi
\Biggl(\frac{1}{q}\sum_{j_1,j_2=1}^q(-1+q
\delta_{j_1j_2}) \nonumber\\[-8pt]\\[-8pt]
&&\hphantom{A_K\int_{-\ppi}^\ppi
\Biggl(\frac{1}{q}\sum_{j_1,j_2=1}^q}{}\times\bigl\{ \bigl|\tr\bigl(\widehat{\mathbf{F}}_{j_1j_2}(
\omega)\bigr)\bigr|^2 +\tr\bigl(\widehat{\mathbf{F}}_{j_1j_1}(
\omega)\overline{\widehat {\mathbf{F}}_{j_2j_2}(\omega)}^T
\bigr) \bigr\} \Biggr)\,\mathrm{d}\omega\nonumber
\end{eqnarray}
to gain asymptotic normality.

\begin{theorem}[(Asymptotic distribution of $T_n^*$)]\label{theorem4}
Suppose that the mixing condition \eqref{mixingcond} holds for all
$k\leq32$. Under the Assumption~\ref{ass2}, we have (conditioned on
$\un{X}_1,\dots,\un{X}_n$) convergence in distribution
%
\begin{eqnarray}
\label{CLT Tn*} T_n^*-\frac{ \widehat{\mu}^*}{\sqrt{h}}\stackrel{\mathcal {D}} {
\longrightarrow}\mathcal{N}\bigl(0,\tau^{*2}\bigr)
\end{eqnarray}
in probability as $n\rightarrow\infty$, where $ \widehat{\mu}^*$ is
as in (\ref{eq61}) and $\tau^{*2}$ is defined by
%
\begin{eqnarray}\label{eq62}
\tau^{* 2} &:=& B_K \int_{-\ppi}^\ppi
\frac{1}{q^2}\sum_{j_1,j_2,j_3,j_4=1}^q
\biggl(-1+q\delta_{j_1j_3}\delta_{j_2j_4}+\frac{q}{q-1}(1-
\delta _{j_1j_3}) (1-\delta_{j_2j_4}) \biggr)
\nonumber
\\
& &\hphantom{B_K \int_{-\ppi}^\ppi
\frac{1}{q^2}\sum_{j_1,j_2,j_3,j_4=1}^q}{} \times \bigl\{\bigl|\tr\bigl(\mathbf{F}_{j_1j_3}(\omega)\overline {
\mathbf{F}_{j_2j_4}(\omega)}^T\bigr)\bigr|^2\\
&&\hphantom{B_K \int_{-\ppi}^\ppi
\frac{1}{q^2}\sum_{j_1,j_2,j_3,j_4=1}^q{} \times \bigl\{}{}+\tr\bigl(
\mathbf{F}_{j_1j_1}(\omega )\overline{\mathbf{F}_{j_2j_2}(
\omega)}^T\bigr)\tr\bigl(\overline{\mathbf {F}_{j_3j_3}(
\omega)}\mathbf{F}_{j_4j_4}(\omega)^T\bigr) \bigr\}\,\mathrm{d}\omega .
\nonumber
\end{eqnarray}
\end{theorem}

Note, that the conditions of the above theorem are weaker than the
assumptions of Theorem~\ref{Theorem 21}.
Moreover, remark that the centering term $ \widehat{\mu}^*$ defined
in \eqref{eq61} fulfills $\widehat{\mu}^* = \mu^* + \mathrm{o}_P(\sqrt {h})$, where
%
\begin{eqnarray}\label{mustern}
\mu^* &:= &A_K\int_{-\ppi}^\ppi \Biggl(
\frac{1}{q}\sum_{j_1,j_2=1}^q(-1+q
\delta_{j_1j_2})\nonumber \\[-8pt]\\[-8pt]
&&\hphantom{A_K\int_{-\ppi}^\ppi \Biggl(
\frac{1}{q}\sum_{j_1,j_2=1}^q}{}\times\bigl\{ \bigl|\tr\bigl(\mathbf{F}_{j_1j_2}(\omega)
\bigr)\bigr|^2 +\tr\bigl(\mathbf{F}_{j_1j_1}(\omega)\overline{
\mathbf {F}_{j_2j_2}(\omega)}^T\bigr) \bigr\} \Biggr)\,\mathrm{d}\omega.\nonumber
\end{eqnarray}
Hence, to compare the above results with Theorem~\ref{theorem1} we
have to analyze $\mu^*$ and $\tau^{* 2}$ under $H_0$ which leads to
$\mu_0^*$ and $\tau_0^{*2}$ in the following remark.

%
\begin{remark}\label{H0remark}
Under $H_0$ the constants $\mu^*$ and $\tau^{* 2}$ of Theorem~\ref
{theorem4}, reduce to
%
\begin{eqnarray}\label{eq6a}
\mu_0^* = A_K\int_{-\ppi}^\ppi
\Biggl(\frac{1}{q}\sum_{j_1,j_2=1}^q(-1+q
\delta_{j_1j_2}) \bigl\{\bigl|\tr\bigl(\mathbf {F}_{j_1j_2}(\omega)
\bigr)\bigr|^2 \bigr\} \Biggr)\,\mathrm{d}\omega= \mu_0
\end{eqnarray}
and
%
\begin{eqnarray}\label{eq6}
\tau_0^{* 2} &=& B_K \int
_{-\ppi}^\ppi \frac{1}{q^2}\sum
_{j_1,j_2,j_3,j_4=1}^q \biggl(-1+q\delta _{j_1j_3}
\delta_{j_2j_4}+\frac{q}{q-1}(1-\delta_{j_1j_3}) (1-\delta
_{j_2j_4}) \biggr)
\nonumber
\\[-8pt]\\[-8pt]
& &\hphantom{B_K \int
_{-\ppi}^\ppi \frac{1}{q^2}\sum
_{j_1,j_2,j_3,j_4=1}^q}{} \times\bigl|\tr\bigl(\mathbf{F}_{j_1j_3}(\omega)\overline{\mathbf
{F}_{j_2j_4}(\omega)}^T\bigr)\bigr|^2\,\mathrm{d}\omega,
\nonumber
\end{eqnarray}
where $\mathbf{F}_{jj}(\omega)=\mathbf{F}_{11}(\omega)$ holds for
all $\omega$ and $j=1,\dots,q$.
\end{remark}

It is interesting to note that the unconditional and conditional
centering parts $\mu_0$ in \eqref{mu} and $\mu_0^*$ in \eqref{eq6a}
coincide.
Moreover, $T_n-h^{-1/2}\mu_0$ as well as $T_n^*-h^{-1/2}\widehat{\mu
}^*$ posses a Gaussian limit distribution under the null.
However, in order to construct an asymptotically exact randomization
test based on $T_n^*-h^{-1/2}\widehat{\mu}^*$ we have to be sure that
the limit variances $\tau^2_0$ in \eqref{tau} and $\tau_0^{* 2}$ in
\eqref{eq6} under $H_0$ are equal as well. This will be discussed in
more detail in the next subsection. Here we just state an example
(corresponding to Example~\ref{exH0special}) where this requirement is
fulfilled.

%
\begin{example}[(The case $p=1$ and $q=2$)]
For $p=1$ and $q=2$ the quantities $\mu^*$ and $\tau^{*2}$ defined in
Theorem~\ref{theorem4} become
\begin{eqnarray*}
\mu^*= A_K\int_{-\ppi}^\ppi
\frac{1}{2} \bigl\{\bigl(f_{11}(\omega )-f_{22}(\omega)
\bigr)^2+f_{11}^2(\omega)+f_{22}^2(
\omega )-2\bigl|f_{12}(\omega)\bigr|^2 \bigr\}\,\mathrm{d}\omega
\end{eqnarray*}
and
\begin{eqnarray*}
\tau^{* 2}=B_K\int_{-\ppi}^\ppi
\frac{1}{4} \bigl\{\bigl(f_{11}(\omega )-f_{22}(\omega)
\bigr)^4+ \bigl(f_{11}^2(\omega)+f_{22}^2(
\omega )-2\bigl|f_{12}(\omega)\bigr|^2 \bigr)^2 \bigr\}\, \mathrm{d}
\omega.
\end{eqnarray*}
Under $H_0$, we have
\begin{eqnarray*}
\mu_0^*= A_K \int_{-\ppi}^\ppi
f_{11}^2(\omega) \bigl(1-C_{12}(\omega )\bigr)\,\mathrm{d}
\omega=\mu_0
\end{eqnarray*}
and
\begin{eqnarray*}
\tau_0^{* 2}=B_K \int_{-\ppi}^\ppi
f_{11}^4(\omega) \bigl(1-C_{12}(\omega)
\bigr)^2\,\mathrm{d}\omega=\tau_0^2.
\end{eqnarray*}
\end{example}

\subsection{The randomization test procedures}

Based on the conditional CLTs above we can now define different
randomization tests.
The first natural approach is to use the test
$
\varphi_{n,\mathrm{cent}}^*:= \mathbf{1}_{(c_{n,\mathrm{cent}}^*(\alpha),\infty
)}(T_n-h^{-1/2}\widehat{\mu})$,
where $c_{n,\mathrm{cent}}^*(\alpha)$ is the data-dependent $(1-\alpha
)$-quantile of the conditional distribution of $T_n^* -
h^{-1/2}\widehat{\mu}^*$ given the data. As typical for resampling
methods, note that we still use the same test statistic as for the
unconditional case and only apply the randomization statistic to
calculate the critical value. By Theorems \ref{theorem1} and \ref
{theorem4}, this test will be asymptotically exact, i.e. $E(\varphi
_{n,\mathrm{cent}}^*)\rightarrow\alpha$ holds under $H_0$, if the asymptotic
variances $\tau^2_0$ and $\tau_0^{*2}$ of the test statistic
$T_n-h^{-1/2}\widehat{\mu}$ and its randomization version
$T_n^*-h^{-1/2}\widehat{\mu}^*$ posses the same limit under the null.
Although $T_n^*-h^{-1/2}\widehat{\mu}^*$ does in general not mimic
the null distribution under the alternative, the following corollary
shows that $\varphi_{n,\mathrm{cent}}^*$ will also be asymptotically consistent
in these situations. Before we state these properties,
we like to introduce a computational less demanding version of the
above test. Therefore note that Remark~\ref{H0remark} implies that the
difference $h^{-1/2}(\widehat{\mu}^*-\widehat{\mu})$ converges to
zero in probability under $H_0$.
Hence, it may be convenient to use the test
%
\begin{equation}
\label{rand.test} \varphi_n^*:= \mathbf{1}_{(c_n^*(\alpha),\infty)}(T_n)
\end{equation}
without centering part, where $c_n^*(\alpha)$ is the data-dependent
$(1-\alpha)$-quantile of the conditional distribution of $T_n^*$.
For completeness, we shortly explain the numerical algorithm for the
implementation of $\varphi_n^*$.
The algorithm for $\varphi_{n,\mathrm{cent}}^*$ is analogue.
\begin{enumerate}[Step 5:]
\item[Step 1:] Compute the test statistic $T_n$ as given in (\ref
{teststat}) based on data $\un{X}_1,\dots,\un{X}_n$.
\item[Step 2:] Generate independent random permutations $\pii
_0,\ldots,\pii_{\lfloor{\sfrac{n}{2}}\rfloor}$ as in (\ref{random
permutation defi}) that are uniformly distributed on the symmetric
group $\mathcal{S}_q$.
\item[Step 3:] Calculate the randomization statistic $T_n^{*}$
(given $\un{X}_1,\dots,\un{X}_n$).
\item[Step 4:] Repeat the Steps $2$ and $3$ $B$-times, where $B$ is
large, which leads to $T_n^{* (1)},\dots,T_n^{* (B)}$.
\item[Step 5:] Reject the null hypothesis $H_0$ if $B^{-1} \sum_{b=1}^B \mathbf{1}\{T_n > T_n^{* (b)} \} >\alpha$.
\end{enumerate}

In the sequel, we analyze the asymptotic properties of both tests
$\varphi_n^*$ and $\varphi_{n,\mathrm{cent}}^*$ and compare it with the
unconditional benchmark test $\varphi_n$.

\begin{corollary}[(Exactness and consistency of $\varphi_n^*$ and
$\varphi_{n,\mathrm{cent}}^*$)]\label{corollary1}
Suppose the assumptions of Theorem~\ref{Theorem 21} are satisfied.
\begin{longlist}[(b)]
\item[(a)] If $H_0$ is true, the following assertions are equivalent:
\begin{enumerate}[(ii)]
\item[(i)] The randomization test $\varphi_n^*$ is asymptotically
exact and equivalent to $\varphi_n$, that is,
%
\begin{eqnarray}
\label{asequiv} E\bigl(\bigl|\varphi_n-\varphi_n^*\bigr|\bigr)
\rightarrow0 \qquad \mbox{as $n\rightarrow \infty$.}
\end{eqnarray}
\item[(ii)] It holds
%
\begin{eqnarray}\label{eq11}
0 &=& \mathop{\mathop{\sum}_{j_1,j_3=1}}_{j_1\neq j_3}^q
\bigl((q-1)^3-1\bigr)\bigl|\tr \bigl(\mathbf {F}_{j_1j_3}(\omega)
\overline{\mathbf{F}_{j_1j_3}(\omega)}^T\bigr)\bigr|^2
\nonumber
\\
& &{} -2\mathop{\mathop{\sum}_{j_1,j_3,j_4=1}}_{all \neq}^q
\bigl((q-1)^2+1\bigr)\bigl|\tr\bigl(\mathbf{F}_{j_1j_3}(\omega)
\overline{\mathbf {F}_{j_1j_4}(\omega)}^T\bigr)\bigr|^2
\\
& &{} +\mathop{\mathop{\sum}_{j_1,j_2,j_3,j_4=1}}_{j_1\neq j_3,j_2\neq
j_4,j_1\neq j_2,j_3\neq j_4}^q
(q-2)\bigl|\tr\bigl(\mathbf {F}_{j_1j_3}(\omega)\overline{\mathbf{F}_{j_2j_4}(
\omega)}^T\bigr)\bigr|^2.
\nonumber
\end{eqnarray}
\end{enumerate}
Moreover, if \eqref{eq11} holds, $\varphi_n^*$ will also be
asymptotically consistent, that is, $E(\varphi_n^*)\rightarrow1$
under $H_1$ as $n\rightarrow\infty$.

\item[(b)] The above statement \textup{(a)} also holds true for
$\varphi_{n,\mathrm{cent}}^*$ instead of $\varphi_n^*$.
\end{longlist}
\end{corollary}

Note that the asymptotic equivalence in (a) implies that both
tests posses the same power for contiguous alternatives.
Moreover, remark that the above stated consistency of the tests is
caused by the fact that $T_n-h^{-1/2}\widehat{\mu}$ converges to
$+\infty$ in probability for fixed alternatives, see the proof section
for more details, and the non-degenerated limit law of its
randomization version.
In the following, we will give some necessary and sufficient conditions
for (ii) above.

\begin{corollary}[(Necessary and sufficient conditions for exactness of
$\varphi_n^*$ and $\varphi_{n,\mathrm{cent}}^*$)]\label{corollary2}
Suppose that the assumptions of Theorem~\ref{Theorem 21} hold.
\begin{enumerate}[(iii)]
\item[(i)] In the case $q=2$, that is, we are testing for the equality
of two $(p\times p)$ spectral density matrices, condition \textup{(ii)} of
Corollary~\ref{corollary1} is fulfilled for all $p\in\N$.
\item[(ii)] For any $q,p\in\N$, $q\geq2$, a sufficient condition
for condition \textup{(ii)} of Corollary~\ref{corollary1} is
%
\begin{eqnarray}
\mathbf{F}_{ij}(\omega)=\mathbf{F}_{12}(\omega)
\label{eq10}
\end{eqnarray}
for all $\omega\in[-\ppi,\ppi]$ and all $i,j\in\{1,\ldots,q\}$ with
$i\neq j$.
This means that all $(p\times p)$ block matrices on all secondary
diagonals have to be equal and Hermitian everywhere.
\item[(iii)] For $p=1$ and $q\geq2$, a sufficient condition weaker
than \textup{(ii)} above is
\begin{eqnarray*}
\bigl|f_{ij}(\omega)\bigr|^2=\bigl|f_{12}(\omega)\bigr|^2
\end{eqnarray*}
for all $\omega\in[-\ppi,\ppi]$ and all $i,j\in\{1,\ldots,q\}$ with
$i\neq j$.
\item[(iv)] For $p=1$ and $q=3$, the condition
\begin{eqnarray*}
\bigl|f_{ij}(\omega)\bigr|^2=\bigl|f_{12}(\omega)\bigr|^2
\end{eqnarray*}
for all $\omega\in[-\ppi,\ppi]$ and all $i,j\in\{1,2,3\}$ with
$i\neq j$ is not only sufficient, but also necessary for condition
\textup{(ii)} of Corollary~\ref{corollary1}.
\end{enumerate}
\end{corollary}

Note that the randomization tests $\varphi_n^*$ and $\varphi
_{n,\mathrm{cent}}^*$ are asymptotically exact particularly in the case of
uncorrelated $p$-variate (sub-)time series $\un{X}_1,\ldots,\un
{X}_q$ due to (ii) above. Intuitively, this makes sense because
permuting the block diagonal matrices distorts the correlation
structure between these time series, if there is any. This explains
also why both randomization tests are asymptotically exact in the more
general case of equal covariance structure between the time series $\un
{X}_1,\ldots,\un{X}_q$ as shown in (ii) and (iii) of Corollary~\ref
{corollary2} above, which is often denoted by equicorrelatedness.

Observe that due to the sophisticated condition in Corollary~\ref
{corollary1}, a result for general $q\geq3$ equivalent to (iv) in
Corollary~\ref{corollary2} does not seem to be achievable.

Nevertheless, there is an obvious possibility for constructing another
randomization test that works even if the presumption \eqref{eq11} of
Corollary~\ref{corollary1} is not fulfilled. We just have to impose an
estimator for the conditional limit variance $\tau^{* 2}$ of $T_n^*$
as introduced in Theorem~\ref{theorem4}. An appropriate candidate is
given by
%
\begin{eqnarray}\label{tauStern}
\widehat{\tau}^{* 2} &:=& B_K \int_{-\ppi}^\ppi
\frac{1}{q^2}\sum_{j_1,j_2,j_3,j_4=1}^q
\biggl(-1+q\delta_{j_1j_3}\delta_{j_2j_4}+\frac{q}{q-1}(1-
\delta _{j_1j_3}) (1-\delta_{j_2j_4}) \biggr)
\nonumber
\\
& &\hphantom{B_K \int_{-\ppi}^\ppi
\frac{1}{q^2}\sum_{j_1,j_2,j_3,j_4=1}^q} {} \times \bigl\{\tr\bigl(\widehat{\mathbf {F}}_{j_1j_1}(
\omega)\overline{\widehat{\mathbf{F}}_{j_2j_2}(\omega)}^T\bigr)
\tr\bigl(\overline{\widehat{\mathbf{F}}_{j_3j_3}(\omega)}\widehat {
\mathbf{F}}_{j_4j_4}(\omega)^T\bigr)\qquad\\
&&\hphantom{B_K \int_{-\ppi}^\ppi
\frac{1}{q^2}\sum_{j_1,j_2,j_3,j_4=1}^q{} \times \bigl\{} {} +\bigl|\tr\bigl(\widehat{
\mathbf{F}}_{j_1j_3}(\omega)\overline{\widehat {\mathbf{F}}_{j_2j_4}(
\omega)}^T\bigr)\bigr|^2 \bigr\}\,\mathrm{d}\omega\nonumber
\end{eqnarray}
and the corresponding randomization test is
$
\varphi_{n,\mathrm{stud}}^*:= \mathbf{1}_{(c_{n,\mathrm{stud}}^*(\alpha),\infty
)}((T_n-h^{-1/2}\widehat{\mu})/\widehat{\tau})$,
where $c_{n,\mathrm{stud}}^*(\alpha)$ denotes the data-dependent $(1-\alpha
)$-quantile of the conditional distribution of $(T_n^*-h^{-1/2}\widehat
{\mu}^*)/\widehat{\tau}^*$ given the data.
Let us shortly state its properties.

%
\begin{corollary}\label{corollary3}
Suppose that the assumptions of Theorem~\ref{Theorem 21} hold. Then the
test $\varphi_{n,\mathrm{stud}}^*$ is not only of asymptotic level $\alpha$
under $H_0$ but also consistent for testing $H_0$ versus $H_1$, that
is, $E(\varphi_{n,\mathrm{stud}}^*)\rightarrow1$ under $H_1$ as $n\rightarrow
\infty$.
\end{corollary}

Finally, we like to note that some of the assumptions can be weakened
for linear processes. In particular, by combining Theorem~\ref
{theorem4} above with results of Dette and Paparoditis
\cite{Dette:2009} we gain
the following remark.

%
\begin{remark}[(Exactness and consistency of the tests for linear
processes)]\label{linear processes}
If the process posses a linear structure, that is, $\un{X}_t=\sum_{j=-\infty}^\infty\boldsymbol{\Psi}_j \un{e}_{t-j}$ for a
$d$-dimensional i.i.d. white noise $(\un{e}_t,t\in\Z)$ and a
sequence of $(d\times d)$-matrices $(\boldsymbol{\Psi}_j)_j=((\psi
_j(r,s))_{r,s})_j$, Condition \eqref{Fouriertransform} is not needed
and the mixing Assumption~\ref{ass1} in the above Corollaries \ref
{corollary1}--\ref{corollary3} can be substituted by the summable
condition $\sum_j |j|^{1/2} |\psi_j(r,s)|<\infty$ and the moment
condition $E(\|e_t\|^{32})<\infty$.
\end{remark}

For the two sample testing problem as described in the introductory
part our randomization test has a nice reading as a symmetry test. This
is part of the following remark.

%
\begin{remark}[(The case $q=2$: Interpretation as a conditional symmetry
test)]
Remark that we can rewrite (in distributional equality) the
summands of our randomization statistic for $q=2$ as
\begin{eqnarray*}
\mathbf{I}_{\pii_k(r),\pii_k(r)}(\omega_k)-\widetilde{ \mathbf{I}}(
\omega_k) \stackrel{\mathcal{D}} {=} \frac{1}{2}
e_k \bigl(\mathbf{I}_{11}(\omega_k)-
\mathbf{I}_{22}(\omega_k)\bigr),
\end{eqnarray*}
where $(e_k)_{k\geq0}$ are i.i.d. signs, that is, independent and on
$\{+1,-1\}$ uniformly distributed r.v.s, and we set $e_{-k}:=e_k$ as
well as $e_{k+sn}:=e_k $ for $s\in\Z$. Hence, the above defined
randomization tests can be interpreted as some general kind of
conditional symmetry tests in this situation.
\end{remark}

%
\begin{remark}[(Choice and influence of the bandwidth)]\label
{Sensitivity Remark}
As already mentioned in the \hyperref[sec1]{Introduction}, our randomization approach
has the nice advantage that the bandwidth is the only tuning parameter
that has to be assessed.
This feature becomes even better as we will see in our extensive
simulation study, see Section~\ref{sec5} below, where our approach
does not react very sensitive to variations of the bandwidth $h$.
However, it is of course desirable to have a detached principle for
selecting the bandwidth.
Therefore, we use a data driven cross validation
method for choosing $h$ for the kernel spectral
density estimation.
This method is due to Beltr\~{a}o and Bloomfield \cite{Beltrao:1987}
and Robinson \cite{Robinson:1991} and has also been
applied by Paparoditis \cite{Papa:2000,Papa:2005,Papa:2009}.
\end{remark}

Finally, we also like to discuss the local power of all randomization tests.

\begin{corollary}\label{corollary4}
Suppose that the assumptions of Theorem~\ref{theorem:
loclaAlternatives} hold.
\begin{enumerate}[(ii)]
\item[(i)] The test $\varphi_{n,\mathrm{stud}}^*$ has the same
local power as the asymptotic test $\varphi_n$, i.e. we have
$E(\varphi_{n,\mathrm{stud}}^*)\rightarrow1- \Phi(u_{1-\alpha} - \nu/\tau
_0)$ under local alternatives (\ref{localAlternatives}) with $\alpha
_n = h^{-1/4}n^{-1/2}$, where $\nu$ is given in Theorem~\ref{theorem:
loclaAlternatives}.
\item[(ii)] In addition, suppose that $\bf{ f}$ in (\ref
{localAlternatives}) satisfies (\ref{eq11}).
Then $\varphi_{n}^*$ and $\varphi_{n,\mathrm{cent}}^*$ also posses the same
power under local alternatives (\ref{localAlternatives}) with $\alpha
_n = h^{-1/4}n^{-1/2}$.
\end{enumerate}
\end{corollary}

\section{Simulation studies}\label{sec5}

In this section, we illustrate the performance of the randomization
procedure in comparison with the asymptotic (unconditional) benchmark
test $\varphi_n$ as described in
the previous sections.

For better lucidity, we thereby only analyze the finite sample behavior
of the computational least-demanding randomization test $\varphi_{n}^*$
as proposed in (\ref{rand.test}).
The other two randomization procedures from Section~\ref{sec3} behave
similar to or slightly worse than $\varphi_{n}^*$.
For more details, we refer the reader to our supplementary material
(cf. Jentsch and Pauly \cite{Jentsch/Pauly:2013}).

\subsection{The setup}
Suppose we observe bivariate time series data $(\un
{X}_t=(X_{t,1},X_{t,2})^T, t=1,\ldots,n)$ and we want to test the null
hypothesis $H_0$ of equality of both corresponding one-dimensional
spectral densities $f_1(\omega)$ and $f_2(\omega)$. In the setup of
Section~\ref{sec1}, this means $q=2$, $p=1$ and $f_j(\omega)=\mathbf
{F}_{jj}(\omega)$, $j=1,2$ and we test
\begin{eqnarray*}
H_0\dvt  \bigl\{f_1(\omega)=f_2(\omega)\mbox{ for all }\omega
\in[-\ppi,\ppi]\bigr\}
\end{eqnarray*}
against
\begin{eqnarray*}
H_1\dvt \bigl\{\exists A\subset\mathcal{B}\bigl([-\ppi,\ppi]\bigr) \mbox{
with } \leb(A)>0: f_1(\omega)\neq f_2(\omega) \mbox{ for all }
\omega\in A\bigr\}.
\end{eqnarray*}

In the following, we consider data from several well-established time
series models. In particular, our analysis includes Gaussian and
non-Gaussian linear time series as well as non-linear time series models.
The linear models under consideration cover moving average (MA) models
and autoregressive (AR) models with innovations following Gaussian,
logistic and double-exponential distributions, respectively.
GARCH models, threshold AR (TAR) and random coefficient autoregressive
(RCA) models are investigated to cover important classes of non-linear
time series.

Although GARCH processes are known to have power law tails (see Basrak
\textit{et al.} \cite{Basrak:2002}, Section~4) and, consequently, only moments up to some
finite order exists, we include GARCH models to investigate the general
performance of the randomization approach for processes that go beyond
our Assumption \ref{ass1}. For the same reasons, we consider also RCA models
in our simulation study.

The performance of the randomization test $\varphi_{n}^*$ in
comparison to the unconditional benchmark test $\varphi_n$ is
investigated under the null and under the alternative.

For all models under consideration, we have generated $T=400$ time
series. For evaluation of the test statistic, the bandwidth has been
chosen by cross validation as proposed in Remark~\ref{Sensitivity
Remark} and is denoted by $h_{\mathrm{CV}}$. Further, we use the
Bartlett--Priestley kernel (see Priestley \cite
{Priestley:1984}, page 448) for which the
constants in \eqref{constants} become $A_K=\frac{6}{5}$ and
$B_K=\frac{2672\ppi}{385}$ for this particular kernel function. For
each time series, the test $\varphi_n$ has been executed with critical
values from normal approximation as discussed in Section~\ref{sec2}
and the randomization test $\varphi_n^*$ as discussed in Section~\ref
{sec3}, where $B=300$ randomization replications have been used.

\subsection{Analysis of the size}\label{Analysis_size}

To investigate the behavior of the tests under the null, we consider
realizations from vector autoregressive models (VAR), vector
moving-average models (VMA) to cover linear time series and from GARCH,
TAR and RCA models to cover non-linear cases. We consider data from the
bivariate vector $\ARR(1)$ model
%
\begin{eqnarray}\label{AReq}
\ARR(1){:} \quad  \un{X}_t=\mathbf{A}\un{X}_{t-1}+
\un{e}_t,\qquad  t\in\Z,
\end{eqnarray}
where $\un{e}_t\sim(0,\boldsymbol{\Sigma}_1)$ is an independent and
identically distributed (i.i.d.) bivariate white noise with covariance
matrix $\boldsymbol{\Sigma}_1$, $\mathbf{A}$ is chosen from
\begin{eqnarray*}
\mathbf{A}_1=\lleft( %
\begin{array} {c@{\quad}c} 0.1 & 0
\\
0 & 0.1 \end{array} %
\rright), \qquad \mathbf{A}_2=\lleft(
\begin{array} {c@{\quad}c} 0.5 & 0
\\
0 & 0.5 \end{array} %
\rright),\qquad  \mathbf{A}_3=\lleft(
\begin{array} {c@{\quad}c} 0.9 & 0
\\
0 & 0.9 \end{array} %
\rright)
\end{eqnarray*}
and $\boldsymbol{\Sigma}_1=\mathbf{Id}$ is the unit matrix. The VAR
model corresponding to coefficient matrix $\mathbf{A}_i$ is denoted by
$\AR_i$.
All models have been investigated for i.i.d. Gaussian innovations $(\un
{e}_t,t\in\Z)$, whereas the most critical model $\AR_3$ has also been
analyzed for
logistic (with c.d.f. $F(x)= (1+\exp(-x))^{-1} $) and double-exponential
distributions (with p.d.f. $f(x) = \exp(-|x|)/2$) of the i.i.d.
innovations $(\un{e}_t,t\in\Z)$, respectively.
Observe that due to the diagonal shape of all involved matrices
$\boldsymbol{\Sigma}_1$ and $\mathbf{A}_i$, we are dealing with two
\emph{independent} univariate time series here. Furthermore, we
consider data from the bivariate vector $\MAA(1)$ model
%
\begin{eqnarray}\label{MAeq}
\MAA(1){:}\quad  \un{X}_t=\mathbf{B}\un{e}_{t-1}+
\un{e}_t,\qquad  t\in\Z,
\end{eqnarray}
where $\un{e}_t\sim(0,\boldsymbol{\Sigma}_2)$ is an i.i.d. bivariate
white noise with covariance matrix $\boldsymbol{\Sigma}_2$, $\mathbf{B}$
is chosen from
\begin{eqnarray*}
\mathbf{B}_1=\lleft( %
\begin{array} {c@{\quad}c} 0.1 & 0.5
\\
0.5 & 0.1 \end{array} %
\rright), \qquad \mathbf{B}_2=\lleft(
\begin{array} {c@{\quad}c} 0.5 & 0.5
\\
0.5 & 0.5 \end{array} %
\rright), \qquad \mathbf{B}_3=\lleft(
\begin{array} {c@{\quad}c} 0.9 & 0.5
\\
0.5 & 0.9 \end{array} %
\rright),
\end{eqnarray*}
and
%
\begin{eqnarray}\label{Sigma_2}
\boldsymbol{\Sigma}_2=\lleft( %
\begin{array} {c@{\quad}c} 1 & 0.5
\\
0.5 & 1 \end{array} %
\rright).
\end{eqnarray}
The VMA model corresponding to $\mathbf{B}_i$ is denoted by $\MA_i$.
Again all models have been investigated for i.i.d. Gaussian innovations
$(\un{e}_t,t\in\Z)$, whereas the model $\MA_3$ has also been
analyzed for
logistic and double-exponential distributions of the innovations, respectively,
as discussed above, but with covariance matrix $\boldsymbol{\Sigma}_2$.
In this setting, we are dealing with two \emph{dependent} time series
whose marginal spectral densities are equal due to the symmetric shape
of $\boldsymbol{\Sigma}_2$ and
$\mathbf{B}_i$. Note that the application of the randomization
technique to the case of two dependent time series is justified by
Corollary~\ref{corollary2}(i).
Also, we investigate three different non-linear time series models.
First, we consider bivariate data $(\un{X}_t=(X_{t,1},X_{t,2})^T,
t=1,\ldots,n)$ from two independent,
but identically distributed univariate $\GARCHA(1,1)$ processes $\{
X_{t,i},t\in\Z\}$, $i=1,2$, with
%
\begin{eqnarray}\label{GARCHeq}
\GARCHA(1,1){:}\quad  X_{t,i}=\sigma_{t,i}e_{t,i},\qquad
\sigma_{t,i}^2=\omega +aX_{t-1,i}^2+b
\sigma_{t-1,i}^2,\qquad  t\in\Z,
\end{eqnarray}
where $\omega=0.01$, $a=0.1$ and the coefficient $b$ is chosen from
\begin{eqnarray*}
b_1=0.2,\qquad  b_2=0.3,\qquad  b_3=0.4.
\end{eqnarray*}
The corresponding models are denoted by $\GARCH_i$, $i=1,2,3$. Further,
two (centered) independent, but identically distributed univariate
$\TARA(1)$ models $\{X_{t,i},t\in\Z\}$, $i=1,2$, that follow the
model equation
%
\begin{eqnarray}\label{TAReq}
\TARA(1){:}\quad  X_{t,i}= %
\cases{ a(1)X_{t-1,i}+e_{t,i},
&\quad  $X_{t-1,i}<0$,
\cr
a(2)X_{t-1,i}+e_{t,i}, &
\quad $X_{t-1,i}\geq0$, } %
 \qquad  t\in\Z
\end{eqnarray}
with coefficients $\un{a}=(a(1),a(2))'$ chosen from
\begin{eqnarray*}
\un{a}_1=(-0.2,0.1)^T,\qquad  \un{a}_2=(-0.3,0.2)^T,\qquad
\un{a}_3=(-0.4,0.3)^T
\end{eqnarray*}
are studied and denoted by $\TAR_i$, $i=1,2,3$. Also two independent,
but identically distributed univariate $\RCAA(1)$ models $\{
X_{t,i},t\in\Z
\}$, $i=1,2$, that follow the model equation
%
\begin{eqnarray}\label{RCAeq}
\RCAA(1){:}\quad  X_{t,i}=a_tX_{t-1}+e_{t,i}, \qquad t
\in\Z,
\end{eqnarray}
where $a_t\sim\mathcal{N}(0,\sigma^2)$ are i.i.d. centered normally
distributed random variables with standard deviation $\sigma$ chosen from
\begin{eqnarray*}
\sigma_1=0.1, \qquad \sigma_2=0.2,\qquad  \sigma_3=0.3
\end{eqnarray*}
are considered and denoted by $\RCA_i$, $i=1,2,3$. For all non-linear
models above, we have used (independent) standard normal i.i.d. white
noise processes $\{e_{t,i}\}$.
%
\begin{sidewaystable}
\tablewidth=\textwidth
\tabcolsep=0pt
\caption{Actual size of $\varphi_n$ and $\varphi_n^*$ for nominal size
$\alpha\in\{1\%,5\%,10\%\}$, sample size $n\in\{50,100,200\}$,
bandwidth $h=c\cdot h_{\mathrm{CV}}$ for $c\in\{0.5,1,1.5\}$, autoregressive
models $\AR_1$--$\AR_3$
and moving-average models $\MA_1$--$\MA_3$, all with Gaussian innovations}
\label{table_AR_H0}
\begin{tabular*}{\textwidth}{@{\extracolsep{\fill}}d{2.0}d{1.1}d{2.1}d{2.1}d{2.1}d{2.1}d{2.1}d{2.1}d{2.1}d{2.1}d{2.1}d{2.1}d{2.1}d{2.1}d{2.1}d{2.1}d{2.1}d{2.1}d{2.1}d{2.1}@{}}
\hline
& & \multicolumn{6}{l}{$\AR_1$} & \multicolumn{6}{l}{$\AR_2$} &
\multicolumn{6}{l}{$\AR_3$}\\[-5pt]
& & \multicolumn{6}{l}{\hrulefill} & \multicolumn{6}{l}{\hrulefill} &
\multicolumn{6}{l}{\hrulefill}\\
& \multicolumn{1}{l}{$n$:} & \multicolumn{2}{l}{50} & \multicolumn
{2}{l}{100} & \multicolumn{2}{l}{200} & \multicolumn{2}{l}{50}
& \multicolumn{2}{l}{100} & \multicolumn{2}{l}{200} &
\multicolumn{2}{l}{50} & \multicolumn{2}{l}{100} & \multicolumn
{2}{l}{200} \\[-5pt]
&  & \multicolumn{2}{l}{\hrulefill} & \multicolumn
{2}{l}{\hrulefill} & \multicolumn{2}{l}{\hrulefill} & \multicolumn{2}{l}{\hrulefill}
& \multicolumn{2}{l}{\hrulefill} & \multicolumn{2}{l}{\hrulefill} &
\multicolumn{2}{l}{\hrulefill} & \multicolumn{2}{l}{\hrulefill} & \multicolumn
{2}{l}{\hrulefill} \\
\multicolumn{1}{l}{$\alpha$} & \multicolumn{1}{l}{$c$} & \multicolumn{1}{l}{$\varphi_n$}
& \multicolumn{1}{l}{$\varphi_n^*$} & \multicolumn{1}{l}{$\varphi_n$} &
\multicolumn{1}{l}{$\varphi_n^*$} & \multicolumn{1}{l}{$\varphi_n$} & \multicolumn{1}{l}{$\varphi_n^*$}
& \multicolumn{1}{l}{$\varphi_n$} &
\multicolumn{1}{l}{$\varphi_n^*$} & \multicolumn{1}{l}{$\varphi_n$} & \multicolumn{1}{l}{$\varphi_n^*$} &
\multicolumn{1}{l}{$\varphi_n$} &
\multicolumn{1}{l}{$\varphi_n^*$} & \multicolumn{1}{l}{$\varphi_n$} & \multicolumn{1}{l}{$\varphi_n^*$} & \multicolumn{1}{l}{$\varphi_n$} &
\multicolumn{1}{l}{$\varphi_n^*$} & \multicolumn{1}{l}{$\varphi_n$} & \multicolumn{1}{l}{$\varphi_n^*$} \\\hline
1 & 0.5 & 5.5 & 1.5 & 8.5 & 2.0 & 6.5 & 1.0 & 10.0 & 1.5 & 7.5 & 1.3 &
7.3 & 1.8 & 17.8 & 2.8 & 13.0 & 1.3 & 11.3 & 1.0 \\
& 1 & 11.3 & 1.3 & 10.5 & 1.0 & 11.3 & 0.8 & 13.8 & 1.5 & 11.3 & 2.3 &
12.0 & 1.3 & 27.0 & 2.0 & 20.0 & 1.5 & 21.8 & 1.8 \\
& 1.5 & 9.8 & 1.5 & 17.0 & 1.8 & 10.5 & 1.8 & 16.8 & 0.5 & 13.3 & 1.0 &
11.5 & 1.5 & 36.3 & 2.0 & 26.5 & 1.3 & 24.0 & 1.0 \\[3pt]
5 & 0.5 & 10.5 & 5.0 & 14.5 & 5.8 & 12.3 & 4.8 & 14.3 & 6.0 & 14.8 &
4.5 & 12.3 & 6.8 & 27.3 & 7.5 & 20.8 & 5.8 & 17.5 & 4.8 \\
& 1 & 17.0 & 8.3 & 17.3 & 6.3 & 17.3 & 4.5 & 20.5 & 5.5 & 17.8 & 6.3 &
18.5 & 7.0 & 33.5 & 6.5 & 27.3 & 5.3 & 27.3 & 6.3 \\
& 1.5 & 18.0 & 5.3 & 21.8 & 8.3 & 17.5 & 5.3 & 23.5 & 4.8 & 18.5 & 4.8
& 17.5 & 4.8 & 42.5 & 7.3 & 33.8 & 4.3 & 30.8 & 5.8 \\[3pt]
10 & 0.5 & 15.0 & 8.3 & 20.3 & 11.8 & 18.3 & 10.8 & 20.0 & 11.0 & 19.0
& 10.0 & 16.3 & 11.0 & 33.0 & 12.8 & 28.8 & 11.5 & 22.5 & 10.5 \\
& 1 & 19.8 & 12.3 & 22.0 & 12.0 & 23.8 & 11.3 & 24.8 & 11.5 & 23.5 &
11.5 & 21.3 & 12.8 & 38.0 & 13.3 & 33.5 & 9.8 & 32.8 & 11.5 \\
& 1.5 & 22.3 & 8.0 & 27.3 & 13.3 & 21.8 & 11.8 & 27.3 & 11.3 & 24.3 &
9.0 & 23.8 & 9.3 & 48.0 & 13.8 & 41.0 & 8.5 & 34.3 & 10.0 \\[9pt]
& & \multicolumn{6}{l}{$\MA_1$} & \multicolumn{6}{l}{$\MA_2$} &
\multicolumn{6}{l}{$\MA_3$}\\[-5pt]
& & \multicolumn{6}{l}{\hrulefill} & \multicolumn{6}{l}{\hrulefill} &
\multicolumn{6}{l}{\hrulefill}\\
& \multicolumn{1}{l}{$n$:} & \multicolumn{2}{l}{50} & \multicolumn
{2}{l}{100} & \multicolumn{2}{l}{200} & \multicolumn{2}{l}{50}
& \multicolumn{2}{l}{100} & \multicolumn{2}{l}{200} &
\multicolumn{2}{l}{50} & \multicolumn{2}{l}{100} & \multicolumn
{2}{l}{200} \\[-5pt]
&  & \multicolumn{2}{l}{\hrulefill} & \multicolumn
{2}{l}{\hrulefill} & \multicolumn{2}{l}{\hrulefill} & \multicolumn{2}{l}{\hrulefill}
& \multicolumn{2}{l}{\hrulefill} & \multicolumn{2}{l}{\hrulefill} &
\multicolumn{2}{l}{\hrulefill} & \multicolumn{2}{l}{\hrulefill} & \multicolumn
{2}{l}{\hrulefill} \\
\multicolumn{1}{l}{$\alpha$} & \multicolumn{1}{l}{$c$} & \multicolumn{1}{l}{$\varphi_n$}
& \multicolumn{1}{l}{$\varphi_n^*$} & \multicolumn{1}{l}{$\varphi_n$} &
\multicolumn{1}{l}{$\varphi_n^*$} & \multicolumn{1}{l}{$\varphi_n$} & \multicolumn{1}{l}{$\varphi_n^*$}
& \multicolumn{1}{l}{$\varphi_n$} &
\multicolumn{1}{l}{$\varphi_n^*$} & \multicolumn{1}{l}{$\varphi_n$} & \multicolumn{1}{l}{$\varphi_n^*$} &
\multicolumn{1}{l}{$\varphi_n$} &
\multicolumn{1}{l}{$\varphi_n^*$} & \multicolumn{1}{l}{$\varphi_n$} & \multicolumn{1}{l}{$\varphi_n^*$} & \multicolumn{1}{l}{$\varphi_n$} &
\multicolumn{1}{l}{$\varphi_n^*$} & \multicolumn{1}{l}{$\varphi_n$} & \multicolumn{1}{l}{$\varphi_n^*$} \\\hline
1 & 0.5 & 4.0 & 2.0 & 4.0 & 1.5 & 4.8 & 1.8 & 6.0 & 1.8 & 5.8 & 0.5 &
4.0 & 1.3 & 6.5 & 1.3 & 5.3 & 1.8 & 5.3 & 1.3 \\
& 1 & 5.0 & 1.5 & 4.3 & 2.0 & 5.5 & 1.5 & 7.5 & 1.8 & 4.3 & 0.5 & 5.5 &
1.8 & 8.3 & 2.0 & 6.5 & 1.0 & 5.8 & 1.8 \\
& 1.5 & 3.3 & 0.5 & 4.0 & 1.5 & 3.3 & 0.5 & 9.5 & 1.3 & 6.8 & 1.3 & 7.0
& 1.0 & 11.5 & 1.0 & 8.5 & 2.0 & 11.0 & 2.5 \\[3pt]
5 & 0.5 & 7.8 & 5.8 & 8.5 & 5.0 & 8.5 & 6.8 & 11.0 & 5.5 & 10.0 & 5.3 &
11.3 & 5.3 & 11.0 & 5.8 & 10.5 & 6.3 & 8.8 & 5.5 \\
& 1 & 8.8 & 5.0 & 9.5 & 4.8 & 10.0 & 6.5 & 13.5 & 7.3 & 8.0 & 4.8 &
11.3 & 6.0 & 13.5 & 7.3 & 12.0 & 6.5 & 11.8 & 5.5 \\
& 1.5 & 8.3 & 4.8 & 7.5 & 4.3 & 9.5 & 6.3 & 14.3 & 6.5 & 11.0 & 5.3 &
11.5 & 5.3 & 17.8 & 5.0 & 13.5 & 5.5 & 17.0 & 7.5 \\[3pt]
10 & 0.5 & 13.0 & 8.3 & 13.5 & 10.3 & 12.3 & 11.5 & 13.8 & 9.5 & 16.0 &
10.3 & 15.8 & 10.3 & 16.0 & 9.8 & 16.5 & 10.8 & 12.8 & 11.3 \\
& 1 & 12.8 & 11.3 & 14.0 & 11.5 & 12.0 & 11.5 & 18.5 & 12.8 & 13.0 &
11.0 & 16.8 & 12.3 & 17.8 & 11.0 & 15.8 & 11.5 & 16.3 & 9.5 \\
& 1.5 & 11.5 & 10.5 & 12.8 & 10.3 & 13.3 & 13.0 & 18.0 & 13.5 & 14.5 &
11.0 & 15.3 & 10.8 & 23.0 & 8.8 & 18.0 & 9.3 & 20.8 & 12.8 \\
\hline
\end{tabular*}
\end{sidewaystable}
%
\begin{table}
\tablewidth=\textwidth
\tabcolsep=0pt
\caption{Actual size of $\varphi_n$ and $\varphi_n^*$ for nominal size
$\alpha\in\{1\%,5\%,10\%\}$, sample size $n\in\{50,100,200\}$,
bandwidth $h=c\cdot h_{\mathrm{CV}}$ for $c\in\{0.5,1,1.5\}$, autoregressive
model $\AR_3$ and moving-average model $\MA_3$, both with logistic and
double-exponential distribution of the innovations}
\label{table_MA_H0}
\begin{tabular*}{\textwidth}{@{\extracolsep{\fill}}ld{2.0}d{1.1}d{2.1}d{2.1}d{2.1}d{2.1}d{2.1}d{2.1}d{2.1}d{2.1}d{2.1}d{2.1}d{2.1}d{2.1}@{}}
\hline
& & & \multicolumn{6}{l}{$\AR_3$} & \multicolumn{6}{l}{$\MA
_3$}\\[-5pt]
& & & \multicolumn{6}{l}{\hrulefill} & \multicolumn{6}{l}{\hrulefill}\\
& & \multicolumn{1}{l}{$n$:} & \multicolumn{2}{l}{50} &
\multicolumn{2}{l}{100} & \multicolumn{2}{l}{200} & \multicolumn
{2}{l}{50} & \multicolumn{2}{l}{100} & \multicolumn{2}{l}{200} \\[-5pt]
& &  & \multicolumn{2}{l}{\hrulefill} &
\multicolumn{2}{l}{\hrulefill} & \multicolumn{2}{l}{\hrulefill} & \multicolumn
{2}{l}{\hrulefill} & \multicolumn{2}{l}{\hrulefill} & \multicolumn{2}{l}{\hrulefill} \\
\multicolumn{1}{@{}l}{$\{e_t\}$} & \multicolumn{1}{c}{$\alpha$} & \multicolumn{1}{c}{$c$} & \multicolumn{1}{c}{$\varphi_n$}
& \multicolumn{1}{c}{$\varphi_n^*$} & \multicolumn{1}{c}{$\varphi
_n$} & \multicolumn{1}{c}{$\varphi_n^*$} & \multicolumn{1}{c}{$\varphi_n$} & \multicolumn{1}{c}{$\varphi_n^*$} &
\multicolumn{1}{c}{$\varphi_n$} &
\multicolumn{1}{c}{$\varphi_n^*$} & \multicolumn{1}{c}{$\varphi_n$} & \multicolumn{1}{c}{$\varphi_n^*$} & \multicolumn{1}{c}{$\varphi_n$} &
\multicolumn{1}{c}{$\varphi_n^*$} \\\hline
Logistic & 1 & 0.5 & 16.3 & 3.3 & 15.8 & 2.5 & 14.3 & 2.0 & 6.8 & 2.8 &
6.3 & 1.5 & 8.0 & 2.0 \\
& & 1 & 33.8 & 4.0 & 23.3 & 0.8 & 21.0 & 1.8 & 9.8 & 1.8 & 9.0 & 2.0 &
10.0 & 2.8 \\
& & 1.5 & 39.3 & 3.3 & 36.3 & 2.0 & 29.8 & 2.5 & 16.5 & 4.8 & 14.5 &
1.8 & 10.5 & 2.5 \\[3pt]
& 5 & 0.5 & 26.8 & 6.5 & 23.8 & 7.3 & 21.0 & 7.0 & 12.8 & 9.3 & 12.0 &
6.8 & 13.5 & 8.3 \\
& & 1 & 40.8 & 8.0 & 30.5 & 7.3 & 27.5 & 4.8 & 17.0 & 8.5 & 16.5 & 8.3
& 17.0 & 9.3 \\
& & 1.5 & 46.3 & 7.5 & 45.3 & 7.5 & 34.8 & 8.8 & 24.0 & 11.0 & 20.0 &
7.8 & 17.5 & 6.5 \\[3pt]
& 10 & 0.5 & 33.3 & 13.5 & 27.8 & 11.8 & 24.5 & 12.8 & 17.5 & 13.3 &
16.5 & 13.5 & 17.0 & 15.3 \\
& & 1 & 45.5 & 12.3 & 35.0 & 11.5 & 33.0 & 10.3 & 19.8 & 13.5 & 21.5 &
14.0 & 21.8 & 16.3 \\
& & 1.5 & 51.5 & 12.5 & 49.0 & 15.3 & 39.5 & 13.5 & 27.3 & 17.5 & 25.0
& 15.5 & 21.8 & 11.5 \\ [6pt]
Double & 1 & 0.5 & 18.0 & 2.3 & 11.3 & 0.5 & 11.0 & 0.8 & 11.5 & 3.3 &
10.5 & 3.8 & 10.3 & 4.0 \\
\quad -exp. & & 1 & 27.8 & 2.8 & 27.8 & 2.0 & 20.0 & 2.0 & 16.0 & 4.0 & 17.0
& 4.8 & 14.3 & 4.0 \\
& & 1.5 & 34.0 & 2.5 & 31.8 & 3.8 & 30.8 & 1.3 & 17.0 & 1.5 & 18.3 &
2.8 & 17.5 & 4.0 \\[3pt]
& 5 & 0.5 & 29.0 & 6.5 & 20.0 & 5.5 & 19.0 & 4.8 & 17.8 & 10.3 & 17.3 &
10.8 & 16.8 & 8.5 \\
& & 1 & 35.5 & 7.3 & 36.5 & 7.3 & 25.5 & 4.3 & 20.8 & 12.0 & 23.5 &
13.5 & 20.3 & 14.8 \\
& & 1.5 & 42.3 & 9.0 & 38.8 & 7.5 & 37.3 & 4.8 & 23.8 & 8.0 & 25.0 &
9.5 & 23.0 & 12.8 \\  [3pt]
& 10 & 0.5 & 36.0 & 11.0 & 26.0 & 11.8 & 22.5 & 8.8 & 23.0 & 16.5 &
22.3 & 15.5 & 19.8 & 15.8 \\
& & 1 & 43.8 & 13.5 & 40.5 & 13.3 & 29.8 & 10.3 & 26.0 & 17.8 & 28.3 &
23.0 & 23.5 & 20.0 \\
& & 1.5 & 46.5 & 14.3 & 42.5 & 14.5 & 40.5 & 12.3 & 29.0 & 15.3 & 29.5
& 16.5 & 27.5 & 20.5 \\
\hline
\end{tabular*}
\end{table}
%
\begin{sidewaystable}
\vspace*{6pt}
\tablewidth=\textwidth
\tabcolsep=0pt
\caption{Actual size of $\varphi_n$ and $\varphi_n^*$ for nominal
size $\alpha\in\{1\%,5\%,10\%\}$, sample size $n\in\{50,100,200\}$,
bandwidth $h=c\cdot h_{\mathrm{CV}}$ for $c\in\{0.5,1,1.5\}$ and non-linear
models $\GARCH_i$, $\TAR_i$ and $\RCA_i$, respectively}
\label{table_nonlinear_H0}
{\fontsize{8.2pt}{10.2pt}\selectfont{\begin{tabular*}{\textwidth}{@{\extracolsep{\fill}}ld{2.0}d{1.1}d{2.1}d{2.1}d{2.1}d{2.1}d{2.1}d{2.1}d{2.1}d{2.1}d{2.1}d{2.1}d{2.1}d{2.1}d{2.1}d{2.1}d{2.1}d{2.1}d{2.1}d{2.1}@{}}
\hline
& & \multicolumn{1}{l}{$i$:} & \multicolumn{6}{l}{$1$} & \multicolumn{6}{l}{$2$} &
\multicolumn{6}{l}{$3$}\\[-5pt]
& &  & \multicolumn{6}{l}{\hrulefill} & \multicolumn{6}{l}{\hrulefill} &
\multicolumn{6}{l}{\hrulefill}\\
& & \multicolumn{1}{l}{$n$:} & \multicolumn{2}{l}{50} &
\multicolumn{2}{l}{100} & \multicolumn{2}{l}{200} & \multicolumn
{2}{l}{50} & \multicolumn{2}{l}{100} & \multicolumn{2}{l}{200}
& \multicolumn{2}{l}{50} & \multicolumn{2}{l}{100} & \multicolumn
{2}{l}{200} \\[-5pt]
& &  & \multicolumn{2}{l}{\hrulefill} &
\multicolumn{2}{l}{\hrulefill} & \multicolumn{2}{l}{\hrulefill} & \multicolumn
{2}{l}{\hrulefill} & \multicolumn{2}{l}{\hrulefill} & \multicolumn{2}{l}{\hrulefill}
& \multicolumn{2}{l}{\hrulefill} & \multicolumn{2}{l}{\hrulefill} & \multicolumn
{2}{l}{\hrulefill} \\
\multicolumn{1}{l}{Model} & \multicolumn{1}{l}{$\alpha$} & \multicolumn{1}{l}{$c$} & \multicolumn{1}{l}{$\varphi_n$} &
\multicolumn{1}{l}{$\varphi_n^*$} & \multicolumn{1}{l}{$\varphi_n$}
& \multicolumn{1}{l}{$\varphi_n^*$} & \multicolumn{1}{l}{$\varphi_n$} & \multicolumn{1}{l}{$\varphi_n^*$} &
\multicolumn{1}{l}{$\varphi_n$} &
\multicolumn{1}{l}{$\varphi_n^*$} & \multicolumn{1}{l}{$\varphi_n$} & \multicolumn{1}{l}{$\varphi_n^*$} & \multicolumn{1}{l}{$\varphi_n$} &
\multicolumn{1}{l}{$\varphi_n^*$} & \multicolumn{1}{l}{$\varphi_n$} & \multicolumn{1}{l}{$\varphi_n^*$} & \multicolumn{1}{l}{$\varphi_n$} &
\multicolumn{1}{l}{$\varphi_n^*$} & \multicolumn{1}{l}{$\varphi_n$} & \multicolumn{1}{l}{$\varphi_n^*$} \\\hline
$\GARCH_i$ & 1 & 0.5 & 13.0 & 3.5 & 11.8 & 3.8 & 14.0 & 3.3 & 11.0 & 3.5
& 15.0 & 2.3 & 14.5 & 2.5 & 12.8 & 2.5 & 14.8 & 4.5 & 16.8 & 6.5 \\
& & 1 & 15.0 & 2.8 & 14.5 & 2.8 & 14.3 & 3.5 & 12.5 & 2.5 & 21.8 & 4.5
& 17.3 & 3.5 & 16.3 & 4.0 & 16.0 & 3.8 & 17.8 & 4.3 \\
& & 1.5 & 17.3 & 3.0 & 16.0 & 3.5 & 21.5 & 3.0 & 22.3 & 3.8 & 18.8 &
4.0 & 22.0 & 3.8 & 22.3 & 3.8 & 21.5 & 4.0 & 20.5 & 3.5 \\[3pt]
& 5 & 0.5 & 20.5 & 11.0 & 17.5 & 9.8 & 22.8 & 10.8 & 19.5 & 12.0 & 24.0
& 11.8 & 22.3 & 9.8 & 20.0 & 11.5 & 23.0 & 13.0 & 23.5 & 14.0 \\
& & 1 & 21.3 & 12.0 & 20.3 & 10.3 & 21.5 & 9.0 & 18.8 & 8.3 & 28.0 &
14.3 & 23.3 & 12.5 & 22.8 & 9.5 & 22.5 & 11.0 & 25.8 & 11.3 \\
& & 1.5 & 25.8 & 7.3 & 21.3 & 8.8 & 26.5 & 11.3 & 29.8 & 11.0 & 24.8 &
12.0 & 26.8 & 13.0 & 25.8 & 12.0 & 28.0 & 14.0 & 27.0 & 9.8 \\[3pt]
& 10 & 0.5 & 24.0 & 16.8 & 22.5 & 15.3 & 27.3 & 19.0 & 26.5 & 18.0 &
28.3 & 17.8 & 28.8 & 19.5 & 27.5 & 20.8 & 28.0 & 19.8 & 29.3 & 21.5 \\
& & 1 & 28.3 & 18.5 & 27.0 & 16.5 & 28.8 & 17.3 & 25.8 & 14.3 & 32.5 &
22.5 & 27.5 & 18.3 & 29.3 & 17.5 & 29.0 & 18.5 & 31.3 & 19.0 \\
& & 1.5 & 29.8 & 13.5 & 25.5 & 13.8 & 32.3 & 17.5 & 33.8 & 17.8 & 29.8
& 16.8 & 30.5 & 19.5 & 30.3 & 18.8 & 33.5 & 20.5 & 31.8 & 17.5 \\[6pt]
$\TAR_i$ & 1 & 0.5 & 7.0 & 1.3 & 8.3 & 1.8 & 6.0 & 1.3 & 8.5 & 1.5 & 8.3
& 0.8 & 8.3 & 1.0 & 10.5 & 1.8 & 10.0 & 1.8 & 11.0 & 3.8 \\
& & 1 & 11.5 & 1.5 & 10.3 & 1.0 & 11.0 & 1.5 & 11.8 & 1.0 & 11.0 & 1.8
& 11.8 & 1.3 & 13.3 & 1.3 & 12.0 & 2.0 & 13.0 & 2.3 \\
& & 1.5 & 15.0 & 1.0 & 12.0 & 1.5 & 12.0 & 0.5 & 16.0 & 2.8 & 13.5 &
2.3 & 12.3 & 2.3 & 11.8 & 0.8 & 13.8 & 2.0 & 13.0 & 2.5 \\[3pt]
& 5 & 0.5 & 13.0 & 5.3 & 14.0 & 5.3 & 11.5 & 5.0 & 15.0 & 6.0 & 12.3 &
6.0 & 14.0 & 5.3 & 15.8 & 6.3 & 15.8 & 8.5 & 19.3 & 10.8 \\
& & 1 & 16.3 & 5.3 & 15.8 & 6.3 & 15.8 & 6.3 & 17.0 & 4.8 & 17.5 & 7.0
& 17.3 & 7.3 & 21.3 & 7.0 & 17.5 & 6.8 & 21.0 & 9.5 \\
& & 1.5 & 20.8 & 6.3 & 17.5 & 4.5 & 18.0 & 5.3 & 23.5 & 8.5 & 21.5 &
6.8 & 17.5 & 5.3 & 17.8 & 5.0 & 19.5 & 5.8 & 22.5 & 6.8 \\[3pt]
& 10 & 0.5 & 16.5 & 8.8 & 19.3 & 13.3 & 15.8 & 10.3 & 19.5 & 13.3 &
19.0 & 11.5 & 18.0 & 12.0 & 22.3 & 10.8 & 21.3 & 14.8 & 23.8 & 15.3 \\
& & 1 & 19.8 & 11.5 & 19.8 & 12.8 & 20.5 & 11.5 & 21.3 & 10.3 & 21.5 &
12.5 & 21.8 & 14.0 & 27.8 & 14.5 & 22.5 & 12.8 & 26.8 & 16.3 \\
& & 1.5 & 23.8 & 10.8 & 20.5 & 8.8 & 22.8 & 9.5 & 29.0 & 14.0 & 26.3 &
12.0 & 24.0 & 11.8 & 22.0 & 9.3 & 24.0 & 10.3 & 27.5 & 11.0 \\[6pt]
$\RCA_i$ & 1 & 0.5 & 6.5 & 2.5 & 7.0 & 1.3 & 8.3 & 0.8 & 9.0 & 1.8 &
11.3 & 2.3 & 11.8 & 2.5 & 10.3 & 3.0 & 12.3 & 4.0 & 11.8 & 2.8 \\
& & 1 & 12.3 & 1.0 & 11.0 & 2.0 & 12.0 & 2.5 & 10.8 & 1.8 & 13.5 & 2.5
& 13.0 & 0.8 & 15.0 & 3.8 & 12.8 & 0.8 & 12.3 & 2.0 \\
& & 1.5 & 12.5 & 2.5 & 13.3 & 1.8 & 14.3 & 1.5 & 15.5 & 4.0 & 12.0 &
1.5 & 15.8 & 2.0 & 15.8 & 0.8 & 16.5 & 1.8 & 17.3 & 2.3 \\[3pt]
& 5 & 0.5 & 13.0 & 6.5 & 13.8 & 7.5 & 14.8 & 5.5 & 14.3 & 6.5 & 18.8 &
9.3 & 18.8 & 8.3 & 17.8 & 9.0 & 18.5 & 11.3 & 19.8 & 10.8 \\
& & 1 & 18.0 & 7.0 & 17.3 & 5.0 & 19.0 & 5.8 & 16.3 & 7.3 & 21.0 & 8.5
& 20.0 & 7.3 & 20.3 & 8.5 & 19.0 & 6.8 & 18.8 & 8.3 \\
& & 1.5 & 18.8 & 6.8 & 17.8 & 4.8 & 21.5 & 6.3 & 23.8 & 8.0 & 18.3 &
4.0 & 22.0 & 7.8 & 20.8 & 6.8 & 25.5 & 8.3 & 24.5 & 7.3 \\[3pt]
& 10 & 0.5 & 20.5 & 12.3 & 21.0 & 13.3 & 17.5 & 11.0 & 20.8 & 14.0 &
24.0 & 15.8 & 23.3 & 15.0 & 22.8 & 14.3 & 25.5 & 18.8 & 24.5 & 16.3 \\
& & 1 & 22.3 & 11.8 & 21.5 & 10.3 & 23.0 & 12.8 & 21.8 & 12.3 & 26.0 &
15.0 & 25.0 & 15.0 & 25.8 & 13.8 & 23.8 & 14.3 & 24.3 & 14.8 \\
& & 1.5 & 23.5 & 12.3 & 20.8 & 8.3 & 27.0 & 13.3 & 28.3 & 14.5 & 21.0 &
9.8 & 25.5 & 13.8 & 26.0 & 12.8 & 31.8 & 14.8 & 29.5 & 15.5 \\
\hline
\end{tabular*}}}
\end{sidewaystable}
%
\begin{table}
\tablewidth=\textwidth
\tabcolsep=0pt
\caption{Actual size of $\varphi_n$ for nominal size $\alpha\in\{1\%
,5\%,10\%\}$, sample size $n\in\{50,100,200,500,1000,\allowbreak 2000\}$,
bandwidth $h_{\mathrm{CV}}$ and models $\AR_3$, $\MA_3$, $\GARCH_3$, $\TAR
_3$ and
$\RCA_3$, respectively}
\label{table_sample_sizes}
\begin{tabular*}{\textwidth}{@{\extracolsep{\fill}}ld{2.0}d{2.1}d{2.1}d{2.1}d{2.1}d{2.1}d{2.1}@{}}
\hline
&&\multicolumn{6}{l}{$n$} \\[-5pt]
&&\multicolumn{6}{l}{\hrulefill} \\
\multicolumn{1}{l}{Model} & \multicolumn{1}{l}{$\alpha$} &  \multicolumn{1}{l}{50} & \multicolumn{1}{l}{100}
& \multicolumn{1}{l}{200} & \multicolumn{1}{l}{500} & \multicolumn{1}{l}{1000} & \multicolumn{1}{l}{2000} \\\hline
$\AR_3$ & 1 &  27.0 & 20.0 & 21.8 & 19.3 & 18.5 & 18.3\\
& 5 &  33.5 & 27.3 & 27.3 & 25.8 & 23.5 & 24.0 \\
& 10 &  38.0 & 33.5 & 32.8 & 30.0 & 28.5 & 27.3 \\[3pt]
$\MA_3$ & 1 &  8.3 & 6.5 & 5.8 & 6.3 & 5.8 & 4.8 \\
& 5 &  13.5 & 12.0 & 11.8 & 11.0 & 10.0 & 9.3 \\
& 10 &  17.8 & 15.8 & 16.3 & 17.0 & 14.0 & 13.3 \\[3pt]
$\GARCH_3$ & 1 &  16.3 & 16.0 & 17.8 & 16.8 & 13.3 & 13.5\\
& 5 &  22.8 & 22.5 & 25.8 & 26.3 & 19.5 & 20.8 \\
& 10 &  29.3 & 29.0 & 31.3 & 30.0 & 27.3 & 26.0 \\[3pt]
$\TAR_3$ & 1 &  13.3 & 12.0 & 13.0 & 10.0 & 8.5 & 7.8 \\
& 5 &  21.3 & 17.5 & 21.0 & 15.5 & 15.3 & 14.8 \\
& 10 &  27.8 & 22.5 & 26.8 & 20.3 & 21.3 & 17.8 \\[3pt]
$\RCA_3$ & 1 & 15.0 & 12.8 & 12.3 & 13.5 & 11.8 & 11.5\\
& 5 &  20.3 & 19.0 & 18.8 & 19.8 & 18.8 & 17.8 \\
& 10 &  25.8 & 23.8 & 24.3 & 25.3 & 24.3 & 22.0 \\
\hline
\end{tabular*}
\end{table}

For nominal sizes $\alpha\in\{1\%,5\%,10\%\}$ and sample sizes $n\in
\{50,100,200\}$, the corresponding results for all combinations are
displayed in Tables~\ref{table_AR_H0}--\ref{table_nonlinear_H0}. To
check how sensitive the tests react on the bandwidth choice, we report
the simulation results for bandwidths $c\cdot h_{\mathrm{CV}}$ and $c\in\{
0.5,1,1.5\}$ to cover under-smoothing and over-smoothing with respect
to the bandwidth $h_{\mathrm{CV}}$ chosen via cross validation.

To illustrate the slow convergence of the actual size of the
unconditional test $\varphi_n$ for the different models under the null
and to emphasize the need of resampling techniques to resolve this
issue, we report its performance also for larger sample sizes in
Table~\ref{table_sample_sizes}.

\subsection{Analysis of the power}

As can be seen in Tables~\ref{table_AR_H0}--\ref{table_nonlinear_H0},
the asymptotic test tends to overreject the null hypothesis
systematically in most situations and
should not be applied, at least for the simulated sample sizes.
In particular, for the autoregressive models where the actually
achieved size is far too large compared to nominal size (see
Tables~\ref{table_AR_H0}--\ref{table_MA_H0}),
the unconditional benchmark test $\varphi_n$ cannot be judged.
Moreover, even for larger sample sizes up to $n=2000$ it does not keep
the prescribed level satisfactorily,
see Table~\ref{table_sample_sizes}. Thereforee, we present here
only  the
power behavior of the randomization test $\varphi_n^*$.
However, we present a small power simulation study of $\varphi_n$ as
well as all other randomization tests in the supplementary material
(cf. Jentsch and Pauly \cite{Jentsch/Pauly:2013}).
It can be seen that there is actually no big difference in the power
behavior (measured as achieved power) between the asymptotic test
$\varphi_n$ and all randomization tests.

To illustrate the behavior of $\varphi_n^*$ under the alternative,
that is for inequality of both spectral densities,
we consider several models belonging to the same model classes that
have already been considered above under the null.
First, we consider realizations from the autoregressive model in \eqref
{AReq}. Here, $\mathbf{A}$ is chosen from
\begin{eqnarray*}
\mathbf{A}_4=\lleft( %
\begin{array} {c@{\quad}c} 0.9 & 0
\\
0 & 0.8 \end{array} %
\rright), \qquad \mathbf{A}_5=\lleft(
\begin{array} {c@{\quad}c} 0.9 & 0
\\
0 & 0.7 \end{array} %
\rright), \qquad \mathbf{A}_6=\lleft(
\begin{array} {c@{\quad}c} 0.9 & 0
\\
0 & 0.6 \end{array} %
\rright)
\end{eqnarray*}
and $\boldsymbol{\Sigma}=\mathbf{Id}$. The corresponding models are now
denoted by $\AR_i$, $i=4,5,6$ and due to the diagonal shape, we are
dealing with two independent time series. In the second case, we
generate realizations from the moving average model in \eqref{MAeq}
with $\mathbf{B}$ chosen from
\begin{eqnarray*}
\mathbf{B}_4=\lleft( %
\begin{array} {c@{\quad}c} 0.5 & 0.5
\\
0.5 & 0.7 \end{array} %
\rright), \qquad \mathbf{B}_5=\lleft(
\begin{array} {c@{\quad}c} 0.5 & 0.5
\\
0.5 & 0.8 \end{array} %
\rright), \qquad \mathbf{B}_6=\lleft(
\begin{array} {c@{\quad}c} 0.5 & 0.5
\\
0.5 & 0.9 \end{array} %
\rright)
\end{eqnarray*}
and $\boldsymbol{\Sigma}=\boldsymbol{\Sigma}_2$ defined in \eqref{Sigma_2}
denoted by $\MA_i$, $i=4,5,6$. Due to non-diagonal shape, we are dealing
with two dependent time series in this case. To investigate the power
performance of the tests, we consider the same non-linear time series
models as above. First, we use independent $\GARCHA(1,1)$ models, where the
first process $\{X_{t,1},t\in\Z\}$ follows $\GARCHA(1,1)$ equation in
\eqref{GARCHeq} with $\omega=0.01$, $a=0.1$, $b=0.2$ and the second
process $\{X_{t,2},t\in\Z\}$ is generated by the same model, where
$\omega=0.01$, $a=0.1$, but $b$ is chosen from
\begin{eqnarray*}
b_4=0.3,\qquad  b_5=0.4,\qquad  b_6=0.5
\end{eqnarray*}
referred to as model $\GARCH_i$, $i=4,5,6$. Two (centered) independent
threshold AR models under the alternative are generated by equation
\eqref{TAReq}, where $\{X_{t,1},t\in\Z\}$ corresponds to
coefficients $(a(1),a(2))=(-0.2,0.1)$ and those for $\{X_{t,2},t\in\Z
\}$ are chosen from
\begin{eqnarray*}
\un{a}_4=(-0.3,0.2)^T,\qquad  \un{a}_5=(-0.4,0.3)^T,\qquad
\un{a}_6=(-0.5,0.4)^T.
\end{eqnarray*}

Remark in this context, that for the asymptotic test $\varphi_n$ we
would observe a larger absolute power since this test is very liberal
in the considered situations.
In the supplementary material (cf. Jentsch and Pauly
\cite{Jentsch/Pauly:2013}), we
therefore present its power performance in comparison to all three
randomization tests
by a comparison of actually achieved size obtained from the simulations
in Section~\ref{Analysis_size} with rejection rates in Tables~\ref
{table_AR_H1}--\ref{table_MA_H1} instead of using the nominal size.
There it can be seen that the power behavior of all three tests
(measured as power in comparison to actual size) appear to be quite
similar for all models under investigation.
%
\begin{table}
\tablewidth=\textwidth
\tabcolsep=0pt
\caption{Power of $\varphi_n^*$ for nominal size $\alpha\in\{1\%,5\%
,10\%\}$, sample size $n\in\{50,100,200\}$,
bandwidth $h=c\cdot h_{\mathrm{CV}}$ for $c\in\{0.5,1,1.5\}$ and autoregressive
models $\AR_4$, $\AR_5$ and $\AR_6$ with Gaussian, logistic and
double-exponential distribution of the innovations}
\label{table_AR_H1}
\begin{tabular*}{\textwidth}{@{\extracolsep{\fill}}ld{2.0}d{1.1}d{2.1}d{2.1}d{2.1}d{2.1}d{2.1}d{2.1}d{2.1}d{2.1}d{2.1}@{}}
\hline
& & & \multicolumn{3}{l}{$\AR_4$} & \multicolumn{3}{l}{$\AR
_5$} &
\multicolumn{3}{l}{$\AR_6$}\\[-5pt]
& & & \multicolumn{3}{l}{\hrulefill} & \multicolumn{3}{l}{\hrulefill} &
\multicolumn{3}{l}{\hrulefill}\\
& & \multicolumn{1}{l}{$n$:} & \multicolumn{1}{l}{50} &
\multicolumn{1}{l}{100} & \multicolumn{1}{l}{200} & \multicolumn
{1}{l}{50} & \multicolumn{1}{l}{100} & \multicolumn{1}{l}{200}
& \multicolumn{1}{l}{50} & \multicolumn{1}{l}{100} & \multicolumn
{1}{l}{200} \\
\multicolumn{1}{l}{$\{e_t\}$} & \multicolumn{1}{l}{$\alpha$} & \multicolumn{1}{l}{$c$} & & & & & & & & & \\\hline
Gauss & 1 & 0.5 & 1.5 & 3.3 & 10.8 & 5.3 & 9.0 & 35.0 & 4.8 & 19.3 &
64.3 \\
& & 1 & 3.0 & 5.3 & 10.0 & 4.5 & 11.8 & 37.0 & 6.8 & 19.0 & 67.0 \\
& & 1.5 & 3.8 & 2.5 & 11.5 & 3.8 & 10.0 & 36.8 & 5.3 & 20.5 & 67.5 \\[3pt]
& 5 & 0.5 & 7.5 & 14.3 & 32.3 & 16.5 & 29.0 & 64.8 & 19.0 & 49.5 & 89.0
\\
& & 1 & 8.5 & 15.3 & 25.8 & 14.5 & 33.8 & 70.3 & 23.3 & 51.5 & 91.0 \\
& & 1.5 & 12.5 & 12.3 & 33.0 & 12.8 & 29.5 & 71.0 & 20.3 & 48.3 & 89.3
\\[3pt]
& 10 & 0.5 & 16.5 & 24.0 & 46.3 & 28.5 & 46.8 & 80.8 & 31.5 & 72.5 &
95.0 \\
& & 1 & 17.0 & 26.0 & 39.0 & 26.0 & 47.5 & 83.8 & 39.5 & 71.8 & 96.3 \\
& & 1.5 & 18.8 & 25.0 & 46.5 & 23.3 & 43.3 & 84.5 & 36.0 & 64.3 & 97.3
\\[6pt]
Logistic & 1 & 0.5 & 3.3 & 5.0 & 10.8 & 6.3 & 10.0 & 36.8 & 8.0 & 24.0
& 64.5 \\
& & 1 & 3.8 & 4.8 & 11.0 & 6.5 & 12.5 & 31.5 & 7.3 & 16.3 & 70.0 \\
& & 1.5 & 2.8 & 5.8 & 10.5 & 6.0 & 12.3 & 34.3 & 10.3 & 22.3 & 66.5 \\[3pt]
& 5 & 0.5 & 10.0 & 17.0 & 29.3 & 16.5 & 31.5 & 65.5 & 22.8 & 55.5 &
88.3 \\
& & 1 & 11.3 & 16.5 & 26.5 & 20.3 & 33.3 & 64.5 & 21.8 & 53.0 & 90.0 \\
& & 1.5 & 9.0 & 18.3 & 26.8 & 16.3 & 30.0 & 64.5 & 24.5 & 52.5 & 91.0 \\[3pt]
& 10 & 0.5 & 18.0 & 26.5 & 45.3 & 25.8 & 48.5 & 77.3 & 37.8 & 69.8 &
94.8 \\
& & 1 & 18.5 & 26.0 & 44.0 & 29.5 & 50.5 & 81.0 & 35.5 & 71.5 & 96.0 \\
& & 1.5 & 18.3 & 26.5 & 39.0 & 26.5 & 44.0 & 81.5 & 39.5 & 70.5 & 98.3
\\[6pt]
Double & 1 & 0.5 & 4.0 & 5.0 & 9.8 & 6.3 & 10.3 & 39.8 & 10.3 & 24.0 &
59.3 \\
\quad -exp. & & 1 & 6.0 & 6.0 & 8.5 & 5.8 & 13.8 & 38.3 & 8.8 & 23.5 & 65.5 \\
& & 1.5 & 4.5 & 6.8 & 12.5 & 8.0 & 15.0 & 35.3 & 9.3 & 23.5 & 69.0 \\[3pt]
& 5 & 0.5 & 13.3 & 16.5 & 27.0 & 18.5 & 31.5 & 67.3 & 27.3 & 54.8 &
86.8 \\
& & 1 & 14.5 & 16.8 & 24.0 & 17.3 & 34.3 & 66.8 & 23.8 & 55.3 & 91.0 \\
& & 1.5 & 12.0 & 17.0 & 32.8 & 19.0 & 38.3 & 66.5 & 24.0 & 50.3 & 89.0
\\[3pt]
& 10 & 0.5 & 22.3 & 25.5 & 43.0 & 28.8 & 48.3 & 82.0 & 41.0 & 75.0 &
93.8 \\
& & 1 & 23.3 & 28.5 & 40.5 & 31.0 & 49.3 & 82.5 & 39.0 & 71.8 & 95.8 \\
& & 1.5 & 21.0 & 23.3 & 43.8 & 31.0 & 51.8 & 79.5 & 39.0 & 66.5 & 97.0
\\
\hline
\end{tabular*}
\end{table}
%
\begin{table}
\tablewidth=\textwidth
\tabcolsep=0pt
\caption{Power of $\varphi_n^*$ for nominal size $\alpha\in\{1\%,5\%
,10\%\}$, sample size $n\in\{50,100,200\}$,
bandwidth $h=c\cdot h_{\mathrm{CV}}$ for $c\in\{0.5,1,1.5\}$ and moving-average
models $\MA_4$, $\MA_5$ and $\MA_6$ with Gaussian, logistic and
double-exponential distribution of the innovations}
\label{table_MA_H1}
\begin{tabular*}{\textwidth}{@{\extracolsep{\fill}}ld{2.0}d{1.1}d{2.1}d{2.1}d{2.1}d{2.1}d{2.1}d{2.1}d{2.1}d{2.1}d{2.1}@{}}
\hline
& & & \multicolumn{3}{l}{$\MA_4$} & \multicolumn{3}{l}{$\MA
_5$} &
\multicolumn{3}{l}{$\MA_6$}\\
[-5pt]
& & & \multicolumn{3}{l}{\hrulefill} & \multicolumn{3}{l}{\hrulefill} &
\multicolumn{3}{l}{\hrulefill}\\
& & \multicolumn{1}{l}{$n$:} & \multicolumn{1}{l}{50} &
\multicolumn{1}{l}{100} & \multicolumn{1}{l}{200} & \multicolumn
{1}{l}{50} & \multicolumn{1}{l}{100} & \multicolumn{1}{l}{200}
& \multicolumn{1}{l}{50} & \multicolumn{1}{l}{100} & \multicolumn
{1}{l}{200} \\
\multicolumn{1}{l}{$\{e_t\}$} & \multicolumn{1}{l}{$\alpha$} & \multicolumn{1}{l}{$c$} & & & & & & & & & \\\hline
Gauss & 1 & 0.5 & 3.3 & 7.0 & 14.3 & 4.5 & 12.3 & 34.8 & 12.8 & 29.8 &
65.5 \\
& & 1 & 4.8 & 11.0 & 17.8 & 7.8 & 19.5 & 38.0 & 12.3 & 38.3 & 69.0 \\
& & 1.5 & 3.3 & 9.5 & 20.8 & 9.5 & 19.0 & 44.3 & 14.0 & 34.8 & 76.5 \\[3pt]
& 5 & 0.5 & 10.0 & 17.0 & 31.0 & 13.8 & 29.3 & 54.0 & 31.5 & 50.0 &
84.8 \\
& & 1 & 12.8 & 21.5 & 34.3 & 21.8 & 40.3 & 62.8 & 29.8 & 62.0 & 87.0 \\
& & 1.5 & 12.5 & 20.8 & 39.0 & 23.3 & 39.0 & 66.3 & 29.3 & 62.0 & 90.0
\\[3pt]
& 10 & 0.5 & 19.8 & 26.8 & 41.3 & 22.8 & 41.0 & 64.5 & 41.5 & 63.8 &
91.5 \\
& & 1 & 20.0 & 34.0 & 47.3 & 31.8 & 53.0 & 74.5 & 43.0 & 72.0 & 94.3 \\
& & 1.5 & 23.3 & 30.5 & 50.0 & 35.3 & 52.8 & 78.8 & 43.3 & 74.5 & 95.5
\\[6pt]
Logistic & 1 & 0.5 & 2.5 & 3.8 & 8.3 & 3.0 & 10.8 & 16.0 & 9.3 & 14.8 &
30.0 \\
& & 1 & 2.5 & 7.0 & 12.3 & 5.3 & 11.8 & 14.5 & 10.8 & 20.0 & 36.8 \\
& & 1.5 & 2.5 & 5.0 & 9.0 & 5.0 & 11.0 & 21.0 & 10.3 & 14.8 & 40.0 \\[3pt]
& 5 & 0.5 & 13.5 & 13.3 & 19.0 & 13.5 & 23.5 & 29.5 & 22.5 & 29.0 &
49.0 \\
& & 1 & 11.8 & 13.3 & 22.8 & 14.8 & 26.0 & 34.3 & 24.3 & 39.3 & 59.5 \\
& & 1.5 & 10.8 & 12.3 & 22.0 & 15.3 & 24.8 & 39.3 & 21.0 & 33.5 & 61.8
\\[3pt]
& 10 & 0.5 & 19.0 & 21.3 & 28.3 & 21.3 & 34.3 & 41.8 & 31.3 & 39.3 &
59.8 \\
& & 1 & 19.8 & 22.5 & 33.0 & 24.0 & 36.8 & 47.3 & 32.3 & 51.0 & 71.3 \\
& & 1.5 & 18.5 & 19.8 & 31.0 & 21.5 & 37.0 & 51.5 & 34.5 & 45.5 & 72.0
\\[6pt]
Double & 1 & 0.5 & 3.5 & 4.8 & 9.0 & 6.3 & 11.8 & 18.5 & 11.3 & 19.8 &
35.0 \\
\quad -exp. & & 1 & 6.5 & 10.5 & 11.3 & 8.0 & 13.3 & 23.0 & 10.3 & 20.5 &
37.3 \\
& & 1.5 & 5.3 & 9.8 & 13.5 & 7.0 & 12.0 & 22.3 & 12.3 & 23.5 & 43.5 \\[3pt]
& 5 & 0.5 & 12.5 & 14.8 & 20.5 & 19.5 & 24.5 & 32.8 & 26.8 & 34.8 &
51.8 \\
& & 1 & 13.3 & 19.8 & 25.3 & 19.5 & 26.0 & 39.3 & 23.8 & 39.3 & 55.5 \\
& & 1.5 & 18.0 & 21.5 & 24.8 & 18.3 & 25.8 & 40.8 & 25.0 & 40.5 & 64.5
\\[3pt]
& 10 & 0.5 & 21.5 & 23.8 & 28.3 & 29.8 & 36.3 & 42.3 & 35.8 & 45.0 &
60.5 \\
& & 1 & 23.0 & 27.3 & 36.3 & 29.3 & 35.8 & 48.5 & 33.0 & 50.3 & 64.5 \\
& & 1.5 & 26.5 & 34.0 & 36.5 & 26.8 & 36.5 & 49.3 & 34.0 & 50.3 & 73.8
\\
\hline
\end{tabular*}
\end{table}

\subsection{Discussion}

From Tables~\ref{table_AR_H0}--\ref{table_sample_sizes}, it can be
seen that the asymptotic test $\varphi_{n}$
has difficulties in keeping the prescribed level and tends to over
rejects the null systematically for all small ($n=50$) and moderate
($n=100$) sample sizes.
Its performance is not even desirable for larger sample sizes ($n\geq200$).
This is the case for all linear and non-linear as well as all three
innovation distributions under consideration.
In Table~\ref{table_AR_H0}, especially for the most critical
autoregressive models $\AR_3$, where the $\AR$ coefficient is near to
unity and the corresponding
spectral densities have a non-flat shape, the null approximation of
$\varphi_n$ is extremely poor and the performance is
unacceptable (see Table~\ref{table_AR_H0}, right panel).
Nevertheless, this generally poor performance is not surprising since
the slow convergence speed of $L_2$-type statistics is already known,
see for instance Paparoditis \cite{Papa:2000}, and
Table~\ref{table_sample_sizes} in this paper.

In comparison to that, the randomization test $\varphi_n^*$ performs
better than $\varphi_n$.
For all models under the null and innovation distributions under
investigation, the randomization test holds the prescribed level more
satisfactorily than $\varphi_n$.

Especially for Gaussian innovations and linear time series, the usage
of $\varphi_n^*$ can be recommended for all dependent and independent
settings and sample sizes.
In these cases, a close inspection of Tables~\ref{table_AR_H0}--\ref
{table_nonlinear_H0} shows also that
the bandwidth choice seems to have only a slight effect on the behavior
of the randomization test, where this choice appears to be more crucial
for $\varphi_n$.
To demonstrate this, compare for instance in Tables~\ref{table_AR_H0}
and \ref{table_MA_H0} the performance of model $\AR_3$ for $n=200$ and
$\alpha=10\%$, where the range of the actual size is
from $22.5\%$ to $40.5\%$ for $\varphi_n$ and from $8.8\%$ to $12.3\%$
for $\varphi_n^*$ over the bandwidths $0.5h_{\mathrm{CV}}, h_{\mathrm{CV}}$ and $1.5
h_{\mathrm{CV}}$.

For linear time series with non-Gaussian innovations the performance of
$\varphi_n$ is still worse than that of $\varphi_n^*$.
However, here the randomization test is slightly more effected by the
choice of the bandwidth and its accuracy under $H_0$ often needs larger
sample sizes ($n\geq100$ or even $n\geq200$).

For non-linear time series, see Table~\ref{table_nonlinear_H0}, a
similar observation can be made. The performance of $\varphi_n$ is
poor for all sample sizes from $50$ to $200$ and $\varphi_n^*$ again
keeps the prescribed level much better.
In particular, for the $\TAR$ models, the control of the nominal size is
quite accurate.
For the RCA and GARCH models under consideration, which are not covered
by our Assumption~\ref{ass1}, the performance of $\varphi_n^*$ is still quite
good and, in particular, improved in comparison to $\varphi_n$.
However, the finite sample performance seems to be more effected by the
bandwidth choice for RCA models and $\varphi_n^*$ still tends to
overreject the null for all GARCH models.

Tables~\ref{table_AR_H1}--\ref{table_nonlinear_H1} show the
power behavior of the randomization test $\varphi_n^*$, where
we compare its power to its nominal size. When studying the panels with
increasing sample sizes (from left to right) the consistency results of
Theorem~\ref{theorem11} and
Corollary~\ref{corollary1} under the alternative can be confirmed by
the simulations.
In particular, for the non-linear RCA and GARCH models a typical
consistency behavior can also be observed.
Similar to the situation under the null, for Gaussian innovations and
linear time series the bandwidth choice does not seem to be crucial for
$\varphi_n^*$,
but for other innovations and non-linear time series it has a
considerable effect on its power behavior.
To this end, we think that the typically applied cross-validation
selector leads to quite adequate finite sample performances.

%
\begin{table}
\tablewidth=\textwidth
\tabcolsep=0pt
\caption{Power of $\varphi_n^*$ for nominal size $\alpha\in\{1\%,5\%
,10\%\}$,
sample size $n\in\{50,100,200\}$, bandwidth $h=c\cdot h_{\mathrm{CV}}$ for
$c\in\{0.5,1,1.5\}$ and non-linear models $\GARCH_i$, $\TAR_i$ and
$\RCA_i$ for $i=4,5,6$, respectively}
\label{table_nonlinear_H1}%
\begin{tabular*}{\textwidth}{@{\extracolsep{\fill}}ld{2.0}d{1.1}d{2.1}d{2.1}d{2.1}d{2.1}d{2.1}d{2.1}d{2.1}d{2.1}d{2.1}@{}}
\hline
& & \multicolumn{1}{l}{$i$:} & \multicolumn{3}{l}{$4$} & \multicolumn{3}{l}{$5$} &
\multicolumn{3}{l}{$6$}\\[-5pt]
& & & \multicolumn{3}{l}{\hrulefill} & \multicolumn{3}{l}{\hrulefill} &
\multicolumn{3}{l}{\hrulefill}\\
& & \multicolumn{1}{l}{$n$:} & \multicolumn{1}{l}{50} &
\multicolumn{1}{l}{100} & \multicolumn{1}{l}{200} & \multicolumn
{1}{l}{50} & \multicolumn{1}{l}{100} & \multicolumn{1}{l}{200}
& \multicolumn{1}{l}{50} & \multicolumn{1}{l}{100} & \multicolumn
{1}{l}{200} \\
\multicolumn{1}{l}{Model} & \multicolumn{1}{l}{$\alpha$} & \multicolumn{1}{l}{$c$} & & & & & & & & & \\\hline
$\GARCH_i$ & 1 & 0.5 & 5.3 & 6.5 & 8.5 & 12.3 & 22.8 & 42.0 & 21.0 &
53.5 & 86.3 \\
& & 1 & 4.3 & 7.5 & 11.8 & 12.3 & 23.8 & 46.5 & 26.5 & 56.0 & 83.0 \\
& & 1.5 & 4.0 & 5.3 & 11.3 & 12.8 & 24.0 & 46.8 & 25.5 & 55.3 & 83.5 \\[3pt]
& 5 & 0.5 & 12.8 & 16.3 & 20.5 & 26.3 & 35.8 & 61.3 & 41.5 & 71.3 &
95.0 \\
& & 1 & 10.3 & 17.8 & 25.3 & 24.5 & 38.5 & 66.8 & 47.3 & 72.3 & 92.8 \\
& & 1.5 & 12.8 & 13.8 & 21.0 & 28.3 & 40.5 & 66.3 & 49.5 & 73.5 & 94.3
\\[3pt]
& 10 & 0.5 & 19.8 & 24.5 & 27.3 & 37.0 & 48.8 & 70.3 & 56.5 & 82.5 &
96.3 \\
& & 1 & 16.8 & 25.3 & 34.8 & 32.8 & 50.0 & 77.3 & 58.5 & 81.0 & 95.8 \\
& & 1.5 & 21.3 & 23.0 & 33.0 & 37.3 & 51.5 & 74.8 & 59.0 & 81.3 & 96.8
\\[6pt]
$\TAR_i$& 1 & 0.5 & 2.5 & 2.3 & 3.5 & 2.0 & 1.3 & 2.8 & 2.0 & 3.8 & 6.8
\\
& & 1 & 1.5 & 1.3 & 1.8 & 1.3 & 0.8 & 3.8 & 2.0 & 3.5 & 5.8 \\
& & 1.5 & 1.5 & 1.5 & 2.8 & 2.0 & 1.3 & 1.0 & 1.8 & 1.5 & 3.3 \\[3pt]
& 5 & 0.5 & 6.3 & 8.3 & 8.0 & 8.0 & 6.0 & 12.8 & 6.5 & 11.3 & 21.3 \\
& & 1 & 6.5 & 5.5 & 4.5 & 6.3 & 4.8 & 9.8 & 6.3 & 10.5 & 18.5 \\
& & 1.5 & 5.8 & 5.5 & 7.8 & 6.3 & 5.5 & 5.8 & 5.8 & 6.5 & 12.5 \\[3pt]
& 10 & 0.5 & 13.3 & 12.5 & 13.3 & 12.8 & 10.3 & 22.3 & 14.0 & 22.5 &
30.8 \\
& & 1 & 12.8 & 10.8 & 10.5 & 11.5 & 9.8 & 15.8 & 11.8 & 17.5 & 27.5 \\
& & 1.5 & 11.5 & 9.3 & 12.5 & 11.8 & 11.0 & 9.0 & 13.8 & 14.0 & 19.3 \\
[6pt]
$\RCA_i$& 1 & 0.5 & 1.8 & 1.3 & 2.3 & 2.5 & 3.0 & 2.8 & 2.3 & 6.0 & 10.5
\\
& & 1 & 1.3 & 2.5 & 2.3 & 1.5 & 3.8 & 2.5 & 2.8 & 5.8 & 12.0 \\
& & 1.5 & 2.3 & 0.8 & 1.3 & 1.5 & 1.5 & 2.5 & 3.8 & 5.8 & 12.0 \\[3pt]
& 5 & 0.5 & 6.0 & 7.3 & 7.0 & 6.8 & 11.3 & 10.8 & 8.5 & 15.0 & 20.3 \\
& & 1 & 5.8 & 7.0 & 8.3 & 6.5 & 11.3 & 9.3 & 13.0 & 15.3 & 27.3 \\
& & 1.5 & 6.8 & 7.5 & 5.5 & 6.3 & 8.8 & 9.8 & 10.3 & 16.0 & 24.5 \\[3pt]
& 10 & 0.5 & 14.8 & 12.5 & 12.5 & 12.3 & 14.8 & 16.5 & 17.3 & 23.3 &
30.3 \\
& & 1 & 10.0 & 14.3 & 12.3 & 12.5 & 16.5 & 15.0 & 23.3 & 23.8 & 39.3 \\
& & 1.5 & 14.0 & 13.0 & 13.5 & 11.0 & 14.8 & 17.0 & 19.8 & 22.8 & 33.0
\\
\hline
\end{tabular*}
\end{table}

Our simulation experience may be summarized as follows:
\begin{itemize}[$\bullet$]
\item[$\bullet$] The randomization technique makes sure that $\varphi_n^*$ keeps
the prescribed level well for very small sample sizes
(especially for most linear time series and $\TAR$ models under consideration)
as shown in Tables~\ref{table_AR_H0}--\ref{table_nonlinear_H0}.
Its performance is \emph{in general} significantly better than those
of the asymptotic benchmark test $\varphi_n$.
\item[$\bullet$] Even for larger sample sizes up to $n=2000$ the asymptotic test
$\varphi_n$ cannot be recommended due to its very slow convergence as
emphasized in Table~\ref{table_sample_sizes}.
\item[$\bullet$] In comparison to other resampling methods applied in time series
analysis, the randomization technique used here has the big advantage
that their performance does not depend on the choice of any tuning
parameter in addition to the bandwidth.
The choice of the bandwidth can be done automatically by standard
methods and does not seem to be as crucial as for $\varphi_n$.
\item[$\bullet$] The performance of $\varphi_{n}^*$ becomes even more excellent
if we compare its behavior for the very small sample size of $n=50$
with the poor performance of the unconditional test.
\item[$\bullet$] In comparison to that the only other known and mathematically
analyzed resampling test for $H_0$, the bootstrap test of Dette and
Paparoditis \cite{Dette:2009},
needs sample sizes of $n\geq512$ for gaining comparable results (see
Table $1$ in their paper).
Remark that they have also modeled a $\operatorname{VAR}(1)$ model
under the null with
Gaussian innovations.
When comparing their results with ours (see Table~\ref{table_AR_H0} in
this paper) note that $\rho=0$ in their paper corresponds to our model
$\AR_3$,
but with the (slightly less critical) choice
\begin{eqnarray*}
\mathbf{A}=\lleft( %
\begin{array} {c@{\quad}c} 0.8 & 0
\\
0 & 0.8 \end{array} %
\rright).
\end{eqnarray*}
\item Furthermore, the power performance of $\varphi_{n}^*$ improves
with increasing sample size (as usual under consistency).
The simulations in the supplementary material even show that the power
behavior (measured as power in comparison to actual size) is similar to
that of the benchmark test $\varphi_n$ and all other randomization tests.
\item Finally, to sum up, the randomization procedure helps to hold the
prescribed level under the null more satisfactorily and does not
forfeit power under the alternative in comparison to the unconditional case.
\end{itemize}
%

\section{Final remarks and outlook}
In this paper, we have introduced novel randomization-type tests for
comparing spectral density matrices.
Their theoretical properties have been analyzed in detail and we have
also studied their finite sample performance in extensive simulation studies.
The asymptotic behavior under the null as well as for fixed and local
alternatives have been developed for non-linear time series under a
joint cumulant condition \eqref{ass1}, whereas for linear processes
finite $32$nd moments are required.
Although these conditions seem to be rather strong, we like to mention
that only a few mathematical results on this topic can be found in the
current literature, where most of them
are only developed for the one-dimensional case.
Note, that such, or stronger linearity or even Gaussianity conditions
are typical assumptions in recent time series publications,
see, for example, Eichler \cite{Eichler:2008},
Dette and Paparoditis \cite{Dette:2009}, Dette and Hildebrandt \cite{Dette:2012}, Jentsch \cite
{Jentsch:2010}, Jentsch and Pauly \cite{Jentsch/Pauly:2010},
Preu{ss} and Hildebrandt \cite{Preuss:2012} or Preuss \textit{et al.} \cite{Preuss:2013}.
Moreover, one of the main contributions of the current paper is the
theoretical justification of the randomization approach, which has yet
not been achieved in the time series context.

In future research, we aim to get rid of the last tuning parameter
(i.e., the bandwidth) by studying a different test statistic that is
based on integrated periodograms rather than kernel spectral density estimators.
In doing so, we also plan to relax the conditions on the process in the
non-linear case and to investigate to what extend our theory can still
be established.
As pointed out by one referee, a promising approach may be given by
substituting our joint cumulant condition by
physical-dependence-type-conditions introduced in Wu
\cite{Wu:2005} and further
studied in for example, Liu and Wu \cite{Liu:2010} or
Xiao and Wu \cite{Xiao:2013}.
However, to our knowledge the theory is currently not available in full
strength for our setting, that is, for multivariate $k$-sample problems
with dependencies and triangular arrays. The latter is required to
derive the local power behavior of our test statistic.

\section{Proofs}\label{sec6}

In the following and for a better lucidity of the proofs, we will use
the abbreviate notation $I_{w,m_1,m_2;k}=I_{pw-m_1,pw-m_2}(\omega_k)$
and $f_{w_1,w_2,m_1,m_2;k}=f_{pw_1-m_1,pw_2-m_2}(\omega_k)$.

\subsection{Proofs of Section \texorpdfstring{\protect\ref{sec2}}{2}}

\begin{pf*}{Proof of Theorem \protect\ref{theorem1}}
Under the Assumptions \ref{ass1} and \ref{ass2}, which accord with
Assumptions 3.1 and 3.3 (with taper function $h(t)\equiv1$) in
Eichler \cite{Eichler:2008}, Theorem~3.5 
in his paper is applicable to the test statistic $T_n$, where the also
required Assumption~3.2
is obviously satisfied due to the quite simple shape of $T_n$. Also, by
direct computation, we get
\begin{eqnarray*}
E(T_n) &=& \frac{h^{\sfrac{1}{2}}}{n}\int_{-\ppi}^\ppi
\sum_{m_1,m_2=0}^{p-1} \sum
_{k_1,k_2} K_h(\omega-\omega_{k_1})K_h(
\omega -\omega_{k_2})
\\
& &\hphantom{\frac{h^{\sfrac{1}{2}}}{n}\int_{-\ppi}^\ppi
\sum_{m_1,m_2=0}^{p-1} \sum
_{k_1,k_2}}{} \times\frac{1}{q^2}\sum_{j_1,j_2=1}^q
\Biggl(\sum_{r=1}^q(1-q
\delta_{j_1r}) (1-q\delta_{j_2r}) \Biggr)\\
&&\hphantom{\frac{h^{\sfrac{1}{2}}}{n}\int_{-\ppi}^\ppi
\sum_{m_1,m_2=0}^{p-1} \sum
_{k_1,k_2}{} \times\frac{1}{q^2}\sum_{j_1,j_2=1}^q}{}\times E(I_{j_1,m_1,m_2;k_1}
\overline{I_{j_2,m_1,m_2;k_2}})\,\mathrm{d}\omega
\\
&=& \frac{1}{nh^{3/2}}\int_{-\ppi}^\ppi\sum
_{m_1,m_2=0}^{p-1} \sum_{k}
K^2 \biggl(\frac{\omega-\omega_{k}}{h} \biggr)
\\
& & \hphantom{\frac{1}{nh^{3/2}}\int_{-\ppi}^\ppi\sum
_{m_1,m_2=0}^{p-1} \sum_{k}}{}\times\frac{1}{q}\sum_{j_1,j_2=1}^q
(-1+q\delta _{j_1j_2})f_{j_1,j_2,m_1,m_1;k}\overline{f_{j_1,j_2,m_2,m_2;k}}\,\mathrm{d}
\omega+\mathrm{o}(1),
\end{eqnarray*}
which is asymptotically equivalent to
\begin{eqnarray*}
\frac{1}{h^{\sfrac{1}{2}}}\mu_0=\frac{1}{h^{\sfrac{1}{2}}}\frac
{1}{2\ppi}\int
_{-\ppi}^\ppi K^2(v)\,\mathrm{d}v\int
_{-\ppi}^\ppi \Biggl(\frac
{1}{q}\sum
_{j_1,j_2=1}^q(-1+q\delta_{j_1j_2})\bigl|\tr\bigl(
\mathbf {F}_{j_1j_2}(\omega)\bigr)\bigr|^2 \Biggr)\,\mathrm{d}\omega.
\end{eqnarray*}
For the variance, we have
\begin{eqnarray*}
& & \hspace*{-4pt}\Var(T_n)
\\
&&\hspace*{-4pt}\quad = \frac{h}{n^2}\int_{-\ppi}^\ppi\int
_{-\ppi}^\ppi\sum_{m_1,m_2,m_3,m_4=0}^{p-1}
\sum_{k_1,k_2,k_3,k_4} K_h(\omega-\omega
_{k_1})K_h(\omega-\omega_{k_2})K_h(
\lambda-\omega_{k_3})K_h(\lambda -\omega_{k_4})
\\
& &\hspace*{-4pt}\hphantom{\quad = \frac{h}{n^2}\int_{-\ppi}^\ppi\int
_{-\ppi}^\ppi\sum_{m_1,m_2,m_3,m_4=0}^{p-1}
\sum_{k_1,k_2,k_3,k_4}}{} \times\frac{1}{q^4} \Biggl(\sum_{j_1,j_2,j_3,j_4=1}^q
\Biggl(\sum_{r_1=1}^q (1-q
\delta_{j_1r_1}) (1-q\delta_{j_2r_1}) \Biggr)\\
&&\hspace*{-4pt}\hphantom{\quad = \frac{h}{n^2}\int_{-\ppi}^\ppi\int
_{-\ppi}^\ppi\sum_{m_1,m_2,m_3,m_4=0}^{p-1}
\sum_{k_1,k_2,k_3,k_4}{} \times\frac{1}{q^4} \Biggl(\sum_{j_1,j_2,j_3,j_4=1}^q}{}\times \Biggl(\sum
_{r_2=1}^q(1-q\delta_{j_3r_2}) (1-q
\delta_{j_4r_2}) \Biggr) \Biggr)
\\
& &\hspace*{-4pt}\hphantom{\quad = \frac{h}{n^2}\int_{-\ppi}^\ppi\int
_{-\ppi}^\ppi\sum_{m_1,m_2,m_3,m_4=0}^{p-1}
\sum_{k_1,k_2,k_3,k_4}}{} \times \bigl\{E(I_{j_1,m_1,m_2;k_1}\overline {I_{j_2,m_1,m_2;k_2}}I_{j_3,m_3,m_4;k_3}
\overline {I_{j_4,m_3,m_4;k_4}})
\\
& &\hspace*{-4pt}\hphantom{\quad = \frac{h}{n^2}\int_{-\ppi}^\ppi\int
_{-\ppi}^\ppi\sum_{m_1,m_2,m_3,m_4=0}^{p-1}
\sum_{k_1,k_2,k_3,k_4}{} \times \bigl\{}{} -E(I_{j_1,m_1,m_2;k_1}\overline {I_{j_2,m_1,m_2;k_2}})\\
&&\hspace*{-4pt}\hphantom{\quad = \frac{h}{n^2}\int_{-\ppi}^\ppi\int
_{-\ppi}^\ppi\sum_{m_1,m_2,m_3,m_4=0}^{p-1}
\sum_{k_1,k_2,k_3,k_4}{} \times \bigl\{{}-}{}\times E(I_{j_3,m_3,m_4;k_3}
\overline {I_{j_4,m_3,m_4;k_4}}) \bigr\}\,\mathrm{d}\omega \,\mathrm{d}\lambda.
\end{eqnarray*}
In evaluation of the difference of the expectations above, the only
asymptotically non-vanishing cases are $\omega_{k_1}=\omega_{k_3}\neq
\omega_{k_2}=\omega_{k_4}$, $\omega_{k_1}=\omega_{k_4}\neq\omega
_{k_2}=\omega_{k_3}$, $\omega_{k_1}=-\omega_{k_3}\neq\omega
_{k_2}=-\omega_{k_4}$ and $\omega_{k_1}=-\omega_{k_4}\neq\omega
_{k_2}=-\omega_{k_3}$, where all of them make the same contribution,
which enables us to consider the first case only equipped with a factor
$4$. In this situation, we have
\begin{eqnarray*}
\bigl\{E(\cdots)-E(\cdots)E(\cdots)\bigr\}=f_{j_1,j_3,m_1,m_4;k_1}\overline
{f_{j_1,j_3,m_2,m_3;k_1}}\overline {f_{j_2,j_4,m_1,m_4;k_2}}f_{j_2,j_4,m_2,m_3;k_2}
\delta_{k_1k_3}\delta _{k_2k_4},
\end{eqnarray*}
which results in
\begin{eqnarray*}
& & \Var(T_n)
\\
&&\quad = \frac{4h}{n^2} \sum_{k_1,k_2} \biggl(\int
_{-\ppi}^\ppi K_h(\omega -
\omega_{k_1})K_h(\omega-\omega_{k_2})\,\mathrm{d}\omega
\biggr)^2\frac
{1}{q^2}\\
&&\hphantom{\quad = \frac{4h}{n^2} \sum_{k_1,k_2}}{}\times\sum_{j_1,j_2,j_3,j_4=1}^q
(-1+q\delta_{j_1j_2}) (-1+q\delta _{j_3j_4})
\\
& &\hphantom{\quad = \frac{4h}{n^2} \sum_{k_1,k_2}{}\times\sum_{j_1,j_2,j_3,j_4=1}^q}{}
\times \Biggl(\sum_{m_1,m_4=0}^{p-1}f_{j_1,j_3,m_1,m_4;k_1}
\overline {f_{j_2,j_4,m_1,m_4;k_2}} \Biggr) \\
&&\hphantom{\quad = \frac{4h}{n^2} \sum_{k_1,k_2}{}\times\sum_{j_1,j_2,j_3,j_4=1}^q}{}\times \Biggl(\sum_{m_2,m_3=0}^{p-1}
\overline {f_{j_1,j_3,m_2,m_3;k_1}}f_{j_2,j_4,m_2,m_3;k_2} \Biggr)+\mathrm{o}(1).
\end{eqnarray*}
The term on the right-hand side is asymptotically equivalent to
\begin{eqnarray*}
\tau_0^2=B_K \int_{-\ppi}^\ppi
\Biggl(\frac{1}{q^2}\sum_{j_1,j_2,j_3,j_4=1}^q
(-1+q\delta_{j_1j_2}) (-1+q\delta_{j_3j_4})\bigl|\tr \bigl(
\mathbf{F}_{j_1j_3}(\omega)\overline{\mathbf{F}_{j_2j_4}(\omega
)}^T\bigr)\bigr|^2 \Biggr)\,\mathrm{d}\omega,
\end{eqnarray*}
which concludes this proof.
\end{pf*}

\begin{pf*}{Proof of Theorem \protect\ref{theorem11}}
Note first that the $\sqrt{nh}$-consistency of the kernel density
estimators, see, for example, Equation (5) in Eichler \cite{Eichler:2008},
implies $\widehat{\mu}-\mu_0 = \mathrm{o}_P(h^{1/2})$. Thus Slutsky's lemma
together with Theorem~\ref{theorem1} show that under $H_0$
\[
\frac{T_n - h^{-1/2}\widehat{\mu}}{\tau_0} \schwach N(0,1)
\]
which corresponds to the limit result for $Q_T$ on page 976 in
Eichler \cite{Eichler:2008}.
Together with the consistency of $\widehat{\tau}^2$ this proves the
asymptotic exactness of $\varphi_n$ under the null hypothesis. In
contrast to that, similar to Eichler's Theorem~5.1, it can be verified
by direct computations and Assumption~\ref{ass2} that $h^{1/2}T_n$
converges in probability to $+\infty$ under the alternative $H_1$.
Thus, the same hold true for $T_n - h^{-1/2}\widehat{\mu}$, which
implicates  consistency of the test.
\end{pf*}

\begin{pf*}{Proof of Theorem \protect\ref{theorem: loclaAlternatives}}
Remark that the imposed uniform mixing conditions still imply the
$\sqrt{nh}$-consistency of the kernel density estimators for the processes.
Hence, the estimators for mean and variance can be replaced by their
limit values and we can apply Theorem 5.4 in Eichler \cite{Eichler:2008} with the
special function
\[
\Psi(\mathbf{f})= \Biggl(\vec \Biggl(\mathbf{F}_{rr} -
\frac
{1}{q}\sum_{j=1}^q
\mathbf{F}_{jj} \Biggr) \Biggr)_{1\leq r\leq q}.
\]
Noting that
\[
\frac{\partial\Psi(\mathbf{Z})}{\partial Z_{i,j}}\biggl |_{\mathbf
{Z} = \mathbf{f}} = \vec \bigl(\bigl(\mathbf{1}\bigl
\{(r,s)=(i,j)\bigr\}\bigr)_{(r,s)\in\Xi} \bigr) -\frac
{1}{q} \vec \bigl(
\bigl(\mathbf{1}\bigl\{ (r,s) \in\Xi_{(i,j)}\bigr\}\bigr)_{(r,s)\in
\Xi
}
\bigr)
\]
for
\[
\Xi_{(i,j)} = %
\cases{ \bigl\{(r,s)\in\Xi\dvt  \exists c\in
\Z_{\neq0} \mbox{ with } (r,s) = (cp+i,cp+j) \bigr\}, &\quad  \mbox{for} $(i,j)
\in\Xi$,\vspace*{2pt}
\cr
\emptyset, &\quad  \mbox{elsewise}} %
\]
and with $\Xi= \bigcup_{c=1}^q \{(r,s)\dvt  1+ (c-1)p \leq r,s\leq cp \}$
that theorem implies, that our standardized test statistic converges in
distribution to a normally distributed random variable with mean
$\nu$ and variance $\tau_0^2$. Hence, the results follows.
\end{pf*}

\subsection{Proofs of Section \texorpdfstring{\protect\ref{sec3}}{3}}

Throughout this section, $E^*(\cdot):=E(\cdot|\un{X}_1,\ldots,\un
{X}_n), \Var^*(\cdot):=\Var(\cdot|\un{X}_1,\ldots,\un{X}_n)$
and
$\Cov^*(\cdot, \cdot):=\Cov(\cdot,\cdot|\un{X}_1,\ldots,\un{X}_n)$
denote mean, variance and covariance conditioned on the data,
respectively.

%
\begin{lemma}[(Brillinger \cite{Brillinger:1981},
Theorem~4.3.2)]\label{lemma2}
Suppose that the mixing condition (\ref{mixingcond}) is satisfied for
all $k\leq s\in\N$. Then for any $k\in\{1,\ldots,s\}$ and $\omega
_1,\ldots,\omega_k\in[0,2\ppi]$, we have
%
\begin{eqnarray}\label{brillinger}
&&\operatorname{cum}\bigl(J_{a_1}(\omega_{1}),\ldots,J_{a_k}(
\omega_k)\bigr)\nonumber\\[-8pt]\\[-8pt]
&&\quad  = \frac{(2\ppi)^{(k/2)-1}}{n^{k/2}}f_k^{(a_1,\ldots,a_k)}(
\omega _1,\ldots,\omega_{k-1})\sum
_{t=1}^{n} \mathrm{e}^{\mathrm{i}t\sum_{j=1}^k\omega_{j}} + \mathrm{O} \biggl(
\frac{1}{n^{k/2}} \biggr), \nonumber
\end{eqnarray}
where $J_{a}(\omega)$ denotes the $a$th component of the $d$-variate
discrete Fourier transform $\un{J}(\omega)$, see (\ref{DFT}),
$f_k^{(a_1,\ldots,a_k)}$ is defined as
\begin{eqnarray*}
&&f_k^{(a_1,\ldots,a_k)}(\omega_1,\ldots,
\omega_{k-1})\\
&&\quad  = \frac{1}{(2\ppi)^{k-1}} \sum_{t_{1},\ldots,t_{k-1}=-\infty}^\infty
c_{a_1,\ldots,a_k}(t_1,\ldots,t_{k-1}) \exp
\bigl(\mathrm{i}(t_{1}\omega_{1} + \cdots+ t_{k-1}
\omega_{k-1})\bigr)
\end{eqnarray*}
the $k$th order cumulant spectral density and $f_k$ is the $(d\times d
\times\cdots\times d)$ ($n$ times) $k$th order
cumulant spectra of the $d$-variate time series $(\un{X}_t,t\in\Z)$.
\end{lemma}

%
\begin{lemma}\label{lemma1}
It holds
\begin{eqnarray*}
&\mbox{\textup{(i)}}&\quad  \sum_{r=1}^q(1-q
\delta_{j_1,\pii_{k_1}(r)}) (1-q\delta_{j_2,\pii
_{k_2}(r)})=-q+q^2\sum
_{r=1}^q\delta_{j_1,\pii_{k_1}(r)}\delta
_{j_2,\pii_{k_2}(r)},
\\
&\mbox{\textup{(ii)}}&\quad  \sum_{r=1}^qE^* \bigl((1-q
\delta_{j_1,\pii
_{k_1}(r)}) (1-q\delta_{j_2,\pii_{k_2}(r)}) \bigr)=\bigl(-q+q^2
\delta _{j_1j_2}\bigr)\delta_{|k_1|,|k_2|},
\\
&\mbox{\textup{(iii)}}&\quad  \sum_{r_1,r_2=1}^qE^* \bigl((1-q
\delta_{j_1,\pii
_{k_1}(r_1)}) (1-q\delta_{j_2,\pii_{k_2}(r_1)}) (1-q\delta_{j_3,\pii
_{k_3}(r_2)})
(1-q\delta_{j_4,\pii_{k_4}(r_2)}) \bigr)
\\
& &\qquad  = q^2 (-1+q\delta_{j_1j_2} ) (-1+q\delta _{j_3j_4}
)\delta_{|k_1|,|k_2|}\delta_{|k_3|,|k_4|}(1-\delta _{|k_1|,|k_3|})
\\
& &\qquad \quad {} + q^2 \biggl(-1+q\delta_{j_1j_3}\delta_{j_2j_4}+
\frac
{q}{q-1}(1-\delta_{j_1j_3}) (1-\delta_{j_2j_4}) \biggr)\\
&&\hphantom{\qquad \quad {} +}{}\times
\delta _{|k_1|,|k_3|}\delta_{|k_2|,|k_4|}(1-\delta_{|k_1|,|k_2|})
\\
& & \qquad \quad {}+ q^2 \biggl(-1+q\delta_{j_1j_4}\delta_{j_2j_3}+
\frac
{q}{q-1}(1-\delta_{j_1j_4}) (1-\delta_{j_2j_3}) \biggr)\\
&&\hphantom{\qquad \quad {} +}{}\times
\delta _{|k_1|,|k_4|}\delta_{|k_2|,|k_3|}(1-\delta_{|k_1|,|k_2|})
\\
& &\qquad \quad {} + C\delta_{|k_1|,|k_2|,|k_3|,|k_4|},
\end{eqnarray*}
for some constant $C<\infty$, where $\delta
_{|k_1|,|k_2|,|k_3|,|k_4|}=1$ if $|k_1|=|k_2|=|k_3|=|k_4|$ and $0$ otherwise.
\end{lemma}

\begin{pf}
The first assertion (i) follows from $\sum_{r=1}^q\delta_{j,\pii
_{k}(r)}=1$. Due to $\widetilde{\mathbb{P}}(\pii_{k}(r)=j)=\frac
{1}{q}$ for all $k,r$ and $j$, we get
%
\begin{eqnarray}\label{eq1}
E^*(\delta_{j_1,\pii_{k_1}(r_1)}\delta_{j_2,\pii_{k_2}(r_2)}) &=& %
\cases{ E^*(
\delta_{j_1,\pii_{k}(r)}\delta_{j_2,\pii_{k}(r)}), &\quad\vspace*{2pt}  $|k_1|=|k_2|,r_1=r_2$,
\cr
E^*(\delta_{j_1,\pii_{k}(r_1)}\delta _{j_2,\pii_{k}(r_2)}), &\quad \vspace*{2pt} $|k_1|=|k_2|,r_1
\neq r_2$,
\cr
E^*(\delta _{j_1,\pii_{k_1}(r)})E^*(\delta_{j_2,\pii_{k_2}(r)}),
& \quad $|k_1|\neq|k_2|$ } %
\\
&=& %
\cases{\displaystyle \frac{1}{q}\delta_{j_1j_2}, &\quad
$|k_1|=|k_2|,r_1=r_2$,\vspace*{4pt}
\cr
\displaystyle \frac{1}{q(q-1)}(1-\delta_{j_1j_2}), &\quad  $|k_1|=|k_2|,r_1
\neq r_2$,\vspace*{4pt}
\cr
\displaystyle \frac{1}{q^2}, &\quad $|k_1|
\neq|k_2|$ } \nonumber%
\end{eqnarray}
which yields (ii). Now consider (iii). First, set $Z_{j,k,r}:=1-q\delta
_{j,\pii_{k}(r)}$ and note that $E^*(Z_{j,k,r})=0$ for all $j,k$ and
$r$. Therefore, we have to consider the following non-vanishing cases
(a) $|k_1|=|k_2|\neq|k_3|=|k_4|$, (b) $|k_1|=|k_3|\neq|k_2|=|k_4|$,
(c) $|k_1|=|k_4|\neq|k_2|=|k_3|$ and (d) $|k_1|=|k_2|=|k_3|=|k_4|$.
The first case (a) gives
\begin{eqnarray*}
q^2 (-1+q\delta_{j_1j_2} ) (-1+q\delta_{j_3j_4} )
\delta_{|k_1|,|k_2|}\delta_{|k_3|,|k_4|}(1-\delta_{|k_1|,|k_3|})
\end{eqnarray*}
which follows immediately from (ii). Case (b) becomes
\begin{eqnarray*}
\sum_{r_1,r_2=1}^qE^* (Z_{j_1,k_1,r_1}Z_{j_2,k_2,r_1}Z_{j_3,k_3,r_2}Z_{j_4,k_4,r_2}
)=\!\sum_{r_1,r_2=1}^qE^* (Z_{j_1,k_1,r_1}Z_{j_3,k_1,r_2}
)E^* (Z_{j_2,k_2,r_1}Z_{j_4,k_2,r_2} )
\end{eqnarray*}
and together with \eqref{eq1}, we get
\begin{eqnarray*}
E^* (Z_{j_1,k_1,r_1}Z_{j_3,k_1,r_2} )=-1+q\delta _{r_1r_2}
\delta_{j_1j_3}+\frac{q}{q-1}(1-\delta_{r_1r_2}) (1-
\delta_{j_1j_3})
\end{eqnarray*}
and an analogue result holds for $E^*
(Z_{j_2,k_2,r_1}Z_{j_4,k_2,r_2} )$. By multiplication and summing
up, the contribution of case (b) becomes
\begin{eqnarray*}
q^2 \biggl(-1+q\delta_{j_1j_3}\delta_{j_2j_4}+
\frac{q}{q-1}(1-\delta _{j_1j_3}) (1-\delta_{j_2j_4}) \biggr)
\end{eqnarray*}
and case (c) contributes analogously.
\end{pf}

With these results, we can analyze the conditional expectation and
variance of $T_n^*$.

\begin{theorem}[(Conditional mean and variance of $T_n^*$)]\label
{theorem2}
Under the assumptions of Theorem~\ref{theorem4}, it holds
\begin{eqnarray*}
E^*\bigl(T_n^*\bigr) = h^{-1/2}\mu^*_n
+\mathrm{o}_P(1),
\end{eqnarray*}
where
%
\begin{eqnarray}
\mu^*_n &:=& \frac{1}{n}\int_{-\ppi}^\ppi
\sum_{m_1,m_2=0}^{p-1} \sum
_{k} K_h^2(\omega-
\omega_{k})\frac{1}{q}
\sum_{j_1,j_2=1}^q
(-1+q\delta _{j_1j_2})
 I_{pj_1-m_1,pj_1-m_2}(\omega_{k})\nonumber
 \\[-8pt]\\[-8pt]
 &&\hphantom{\frac{1}{n}\int_{-\ppi}^\ppi
\sum_{m_1,m_2=0}^{p-1} \sum
_{k} K_h^2(\omega-
\omega_{k})\frac{1}{q}
\sum_{j_1,j_2=1}^q}{}\times\overline {I_{pj_2-m_1,pj_2-m_2}(
\omega_{k})}\,\mathrm{d}\omega.\nonumber
\end{eqnarray}
Moreover, we have
\begin{eqnarray*}
\Var^*\bigl(T_n^*\bigr) &=& \frac{4h}{n^2}\mathop{\mathop{\sum
}_{k_1,k_2}}_{|k_1|\neq|k_2|} \biggl(\int_{-\ppi}^\ppi
K_h(\omega -\omega _{k_1})K_h(\omega-
\omega_{k_2})\,\mathrm{d}\omega \biggr)^2
\\
& & \hphantom{\frac{4h}{n^2}\mathop{\mathop{\sum
}_{k_1,k_2}}_{|k_1|\neq|k_2|}}{}\times\sum
_{m_1,m_2,m_3,m_4=0}^{p-1}\frac{1}{q^2}\sum_{j_1,j_2,j_3,j_4=1}^q
\biggl(-1+q\delta _{j_1j_3}\delta_{j_2j_4}\\
&&\hphantom{\frac{4h}{n^2}\mathop{\mathop{\sum
}_{k_1,k_2}}_{|k_1|\neq|k_2|}{}\times\sum
_{m_1,m_2,m_3,m_4=0}^{p-1}\frac{1}{q^2}\sum_{j_1,j_2,j_3,j_4=1}^q
\biggl(}{}+\frac{q}{q-1}(1-
\delta_{j_1j_3}) (1-\delta _{j_2j_4}) \biggr)
\\
& &\hphantom{\frac{4h}{n^2}\mathop{\mathop{\sum
}_{k_1,k_2}}_{|k_1|\neq|k_2|}{}\times\sum
_{m_1,m_2,m_3,m_4=0}^{p-1}\frac{1}{q^2}\sum_{j_1,j_2,j_3,j_4=1}^q
}{} \times I_{j_1,m_1,m_2;k_1}\overline{I_{j_2,m_1,m_2;k_2}}\\
&&\hphantom{\frac{4h}{n^2}\mathop{\mathop{\sum
}_{k_1,k_2}}_{|k_1|\neq|k_2|}{}\times\sum
_{m_1,m_2,m_3,m_4=0}^{p-1}\frac{1}{q^2}\sum_{j_1,j_2,j_3,j_4=1}^q
}{} \times I_{j_3,m_3,m_4;k_1}\overline{I_{j_4,m_3,m_4;k_2}}+\mathrm{o}_P(1)
\\
&=:& \tau_n^{*2}+\mathrm{o}_P(1).
\end{eqnarray*}
\end{theorem}

\begin{pf}
By using Lemma~\ref{lemma1}(ii), we get
\begin{eqnarray*}
E^*\bigl(T_n^*\bigr) &=& \frac{h^{\sfrac{1}{2}}}{n}\int_{-\uppi}^\uppi
\sum_{m_1,m_2=0}^{p-1} \sum
_{k_1,k_2} K_h(\omega-\omega_{k_1})K_h(
\omega -\omega_{k_2})\frac{1}{q^2}
\\
& & \hphantom{\frac{h^{\sfrac{1}{2}}}{n}\int_{-\pii}^\pii
\sum_{m_1,m_2=0}^{p-1} \sum
_{k_1,k_2}}{}\times\sum
_{j_1,j_2=1}^q E^* \Biggl(\sum_{r=1}^q(1-q
\delta_{j_1,\pii
_{k_1}(r)}) (1-q\delta_{j_2,\pii_{k_2}(r)}) \Biggr)\\
&&\hphantom{\frac{h^{\sfrac{1}{2}}}{n}\int_{-\pii}^\pii
\sum_{m_1,m_2=0}^{p-1} \sum
_{k_1,k_2}{}\times\sum
_{j_1,j_2=1}^q}{}\times I_{j_1,m_1,m_2;k_1}
\overline{I_{j_2,m_1,m_2;k_2}}\,\mathrm{d}\omega
\\
&=& \frac{h^{\sfrac{1}{2}}}{n}\int_{-\uppi}^\uppi\sum
_{m_1,m_2=0}^{p-1} \sum_{k_1,k_2}
K_h(\omega-\omega_{k_1})K_h(\omega -
\omega_{k_2})\frac{1}{q}\\
&&\hphantom{\frac{h^{\sfrac{1}{2}}}{n}\int_{-\pii}^\pii\sum
_{m_1,m_2=0}^{p-1} \sum_{k_1,k_2}}{}\times\sum_{j_1,j_2=1}^q
(-1+q\delta _{j_1j_2} )\delta_{|k_1|,|k_2|}
\\
& &\hphantom{\frac{h^{\sfrac{1}{2}}}{n}\int_{-\pii}^\pii\sum
_{m_1,m_2=0}^{p-1} \sum_{k_1,k_2}{}\times\sum_{j_1,j_2=1}^q}{}
\times I_{j_1,m_1,m_2;k_1}\overline{I_{j_2,m_1,m_2;k_2}}\,\mathrm{d}\omega.
\end{eqnarray*}
Furthermore, the case $k_1=-k_2$ is asymptotically negligible, which
yields the first assertion of Theorem~\ref{theorem2}. Considering the
conditional variance of $T_n^*$, we get
\begin{eqnarray*}
& & \Var^*\bigl(T_n^{*2}\bigr)
\\
&&\quad = \frac{h}{n^2}\int_{-\ppi}^\ppi\int
_{-\ppi}^\ppi\sum_{m_1,m_2,m_3,m_4=0}^{p-1}
\sum_{k_1,k_2,k_3,k_4} K_h(\omega-\omega
_{k_1})K_h(\omega-\omega_{k_2}) K_h(
\lambda-\omega_{k_3})K_h(\lambda -\omega_{k_4})
\\
& &\hphantom{\quad = \frac{h}{n^2}\int_{-\ppi}^\ppi\int
_{-\ppi}^\ppi\sum_{m_1,m_2,m_3,m_4=0}^{p-1}
\sum_{k_1,k_2,k_3,k_4}}{} \times\frac{1}{q^4}\\
&&\hphantom{\quad = \frac{h}{n^2}\int_{-\ppi}^\ppi\int
_{-\ppi}^\ppi\sum_{m_1,m_2,m_3,m_4=0}^{p-1}
\sum_{k_1,k_2,k_3,k_4}}{} \times\sum_{j_1,j_2,j_3,j_4=1}^q
\sum_{r_1,r_2=1}^q \operatorname{Cov}^* (Z_{j_1,k_1,r_1}Z_{j_2,k_2,r_1},\\
&&\hphantom{\quad = \frac{h}{n^2}\int_{-\ppi}^\ppi\int
_{-\ppi}^\ppi\sum_{m_1,m_2,m_3,m_4=0}^{p-1}
\sum_{k_1,k_2,k_3,k_4 }{}\times\sum_{j_1,j_2,j_3,j_4=1}^q
\sum_{r_1,r_2=1}^q\operatorname{Cov}^* (} Z_{j_3,k_3,r_2}Z_{j_4,k_4,r_2}
)
\\
& &\hphantom{\quad = \frac{h}{n^2}\int_{-\ppi}^\ppi\int
_{-\ppi}^\ppi\sum_{m_1,m_2,m_3,m_4=0}^{p-1}
\sum_{k_1,k_2,k_3,k_4}{}\times\sum_{j_1,j_2,j_3,j_4=1}^q
\sum_{r_1,r_2=1}^q}{} \times I_{j_1,m_1,m_2;k_1}\overline{I_{j_2,m_1,m_2;k_2}} \\
&&\hphantom{\quad = \frac{h}{n^2}\int_{-\ppi}^\ppi\int
_{-\ppi}^\ppi\sum_{m_1,m_2,m_3,m_4=0}^{p-1}
\sum_{k_1,k_2,k_3,k_4}{}\times\sum_{j_1,j_2,j_3,j_4=1}^q
\sum_{r_1,r_2=1}^q}{}\times I_{j_3,m_3,m_4;k_3}
\overline{I_{j_4,m_3,m_4;k_4}}\,\mathrm{d}\omega \,\mathrm{d}\lambda.
\end{eqnarray*}
By using Lemma~\ref{lemma1}(ii) and (iii), only the cases
$|k_1|=|k_3|\neq|k_2|=|k_4|$ and $|k_1|=|k_4|\neq|k_2|=|k_3|$ play a
role asymptotically. More precisely, only the four cases $k_1=k_3\neq
k_2=k_4$ and $k_1=-k_3 \neq k_2=-k_4$ and $k_1=k_4 \neq k_2=k_3$ and
$k_1=-k_4 \neq k_2=-k_3$ do not vanish, but all of them make the same
contribution asymptotically. This gives a factor $4$ and we get up to
an $\mathrm{o}_P(1)$-term
\begin{eqnarray*}
& & \Var^*\bigl(T_n^{*2}\bigr)
\\
&&\quad = \frac{4h}{n^2}\int_{-\ppi}^\ppi\int
_{-\ppi}^\ppi\sum_{m_1,m_2,m_3,m_4=0}^{p-1}
\mathop{\mathop{\sum }_{k_1,k_2}}_{|k_1|\neq|k_2|} K_h(
\omega-\omega_{k_1})K_h(\omega-\omega_{k_2})K_h(
\lambda-\omega _{k_1})K_h(\lambda-\omega_{k_2})
\\
& & \hphantom{\quad = \frac{4h}{n^2}\int_{-\ppi}^\ppi\int
_{-\ppi}^\ppi\sum_{m_1,m_2,m_3,m_4=0}^{p-1}
\mathop{\mathop{\sum }_{k_1,k_2}}_{|k_1|\neq|k_2|}}{}\times\frac{1}{q^2}\sum_{j_1,j_2,j_3,j_4=1}^q
\biggl(-1+q\delta _{j_1j_3}\delta_{j_2j_4}\\
&& \hphantom{\quad = \frac{4h}{n^2}\int_{-\ppi}^\ppi\int
_{-\ppi}^\ppi\sum_{m_1,m_2,m_3,m_4=0}^{p-1}
\mathop{\mathop{\sum }_{k_1,k_2}}_{|k_1|\neq|k_2|}{}\times\frac{1}{q^2}\sum_{j_1,j_2,j_3,j_4=1}^q
\biggl(}{}+\frac{q}{q-1}(1-
\delta_{j_1j_3}) (1-\delta _{j_2j_4}) \biggr)
\\
& &\hphantom{\quad = \frac{4h}{n^2}\int_{-\ppi}^\ppi\int
_{-\ppi}^\ppi\sum_{m_1,m_2,m_3,m_4=0}^{p-1}
\mathop{\mathop{\sum }_{k_1,k_2}}_{|k_1|\neq|k_2|}{}\times\frac{1}{q^2}\sum_{j_1,j_2,j_3,j_4=1}^q
}{} \times I_{j_1,m_1,m_2;k_1}\overline{I_{j_2,m_1,m_2;k_2}}\\
& &\hphantom{\quad = \frac{4h}{n^2}\int_{-\ppi}^\ppi\int
_{-\ppi}^\ppi\sum_{m_1,m_2,m_3,m_4=0}^{p-1}
\mathop{\mathop{\sum }_{k_1,k_2}}_{|k_1|\neq|k_2|}{}\times\frac{1}{q^2}\sum_{j_1,j_2,j_3,j_4=1}^q
}{} \times I_{j_3,m_3,m_4;k_1}
\overline{I_{j_4,m_3,m_4;k_2}}\,\mathrm{d}\omega \,\mathrm{d}\lambda,
\end{eqnarray*}
which concludes this proof.
\end{pf}

\begin{theorem}[(On $E(\mu_n^*)$, $\Var(\mu_n^*)$, $E(\tau_n^{*2})$
and $\Var(\tau_n^{*2})$)]\label{theorem3}
Under the assumptions of Theorem~\ref{theorem4}, it holds
\begin{eqnarray}\label{tau_var}
E\bigl(\mu_n^*\bigr) &=& \mu^*+\mathrm{o}\bigl(h^{1/2}\bigr)\quad  \mbox{and}\quad  \Var\bigl(\mu _n^*\bigr)=\mathrm{o}(h),\nonumber
\\[-8pt]\\[-8pt]
E\bigl(\tau_n^{*2}\bigr) &= &\tau^{* 2}+\mathrm{o}(1) \quad \mbox{and} \quad \Var\bigl(\tau_n^{*2}\bigr)=\mathrm{o}(1) \nonumber
\end{eqnarray}
as $n\rightarrow\infty$.
\end{theorem}

\begin{pf}
Because arguments are completely analogue, we prove only the more
complicated part that deals with $\tau_n^{*2}$. First, by introducing
the notation $J_{j,m;k}=\frac{1}{\sqrt{2\ppi n}}\sum_{t=1}^n
X_{t,pj-m}\mathrm{e}^{-\mathrm{i}t\omega_k}$, we get
%
\begin{eqnarray}\label{eq2}
\!\!\!\!\!\!\!\!\!\!\!\!\!\!\!\!\!\!E\bigl(\tau_n^{*2}\bigr) &=& \frac{4h}{n^2}\mathop{
\mathop{\sum }_{k_1,k_2}}_{|k_1|\neq|k_2|} \biggl(\int
_{-\ppi}^\ppi K_h(\omega -\omega
_{k_1})K_h(\omega-\omega_{k_2})\,\mathrm{d}\omega
\biggr)^2 \nonumber\\
&&\hphantom{\frac{4h}{n^2}\mathop{
\mathop{\sum }_{k_1,k_2}}_{|k_1|\neq|k_2|}}{}\times\sum_{m_1,m_2,m_3,m_4=0}^{p-1}
\frac{1}{q^2}\sum_{j_1,j_2,j_3,j_4=1}^q
\biggl(-1+q\delta _{j_1j_3}\delta_{j_2j_4}\nonumber\\
&&\hphantom{\frac{4h}{n^2}\mathop{
\mathop{\sum }_{k_1,k_2}}_{|k_1|\neq|k_2|}{}\times\sum_{m_1,m_2,m_3,m_4=0}^{p-1}
\frac{1}{q^2}\sum_{j_1,j_2,j_3,j_4=1}^q
\biggl(}{}+\frac{q}{q-1}(1-
\delta_{j_1j_3}) (1-\delta _{j_2j_4}) \biggr)\nonumber
\\[-8pt]\\[-8pt]
& & \hphantom{\frac{4h}{n^2}\mathop{
\mathop{\sum }_{k_1,k_2}}_{|k_1|\neq|k_2|}{}\times\sum_{m_1,m_2,m_3,m_4=0}^{p-1}
\frac{1}{q^2}\sum_{j_1,j_2,j_3,j_4=1}^q
}{} \times E(J_{j_1,m_1;k_1}\overline{J_{j_1,m_2;k_1}} J_{j_3,m_3;k_1}\nonumber\\
&&\hphantom{\frac{4h}{n^2}\mathop{
\mathop{\sum }_{k_1,k_2}}_{|k_1|\neq|k_2|}{}\times
\sum_{m_1,m_2,m_3,m_4=0}^{p-1}
\frac{1}{q^2}\sum_{j_1,j_2,j_3,j_4=1}^q
{} \times E}{}\times
\overline{J_{j_3,m_4;k_1}}\overline {J_{j_2,m_1;k_2}}\nonumber\\
&&\hphantom{\frac{4h}{n^2}\mathop{
\mathop{\sum }_{k_1,k_2}}_{|k_1|\neq|k_2|}{}\times\sum_{m_1,m_2,m_3,m_4=0}^{p-1}
\frac{1}{q^2}\sum_{j_1,j_2,j_3,j_4=1}^q
{} \times E }{}\times J_{j_2,m_2;k_2}
\overline {J_{j_4,m_3;k_2}}J_{j_4,m_4;k_2}).
\nonumber
\end{eqnarray}
The last expectation above can be expressed in terms of cumulants and
we have
%
\begin{eqnarray}
E(J_1\cdots J_8)=\sum_{\sigma}
\prod_{B\in\sigma} \operatorname{cum}(J_i\dvt i\in B),
\label{eq3}
\end{eqnarray}
where we identify the DFTs with $J_1,\ldots,J_8$ for notational
convenience, $\operatorname{cum}(J_i\dvt i\in B)$ is the joint cumulant of $(J_i\dvt i\in B)$,
$\sigma$ runs through the list of all partitions of $\{1,\ldots,8\}$
and $B$ runs through the list of all blocks of the partition $\sigma$.
Note that the largest contribution in \eqref{eq3} is made by products
of four cumulants of second order. Moreover, all of those combinations
where cumulants $\operatorname{cum}(J_{\cdot,\cdot;k_i},J_{\cdot,\cdot;k_j})$ or
$\operatorname{cum}(J_{\cdot,\cdot;k_i},\overline{J_{\cdot,\cdot;k_j}})$ for
$i\neq j$ occur, contain a factor $\delta_{k_1k_2}$ or $\delta
_{k_1-k_2}$, respectively, but these cases are excluded in the
summation in \eqref{eq2}. The cases with $\operatorname{cum}(J_{\cdot,\cdot
;k_i},J_{\cdot,\cdot;k_i})$ or $\operatorname{cum}(\overline{J_{\cdot,\cdot
;k_i}},\overline{J_{\cdot,\cdot;k_i}})$ contain a factor $\delta
_{k_i0}$, which causes the sum over $k_i$ to collapse and, therefore,
they are of lower order and asymptotically negligible. The only two
remaining products consisting of four cumulants of second order are by
(\ref{brillinger})
\begin{eqnarray*}
& & \operatorname{cum}(J_{j_1,m_1;k_1},\overline{J_{j_1,m_2;k_1}})\operatorname{cum}( J_{j_3,m_3;k_1},
\overline{J_{j_3,m_4;k_1}})
\\
& &\qquad {} \times \operatorname{cum}(\overline {J_{j_2,m_1;k_2}},J_{j_2,m_2;k_2})\operatorname{cum}(\overline
{J_{j_4,m_3;k_2}},J_{j_4,m_4;k_2})
\\
&&\quad = f_{j_1,j_1,m_1,m_2;k_1}f_{j_3,j_3,m_3,m_4;k_1}\overline {f_{j_2,j_2,m_1,m_2;k_2}}
\overline{f_{j_4,j_4,m_3,m_4;k_2}}+\mathrm{O} \biggl(\frac{1}{n} \biggr)
\end{eqnarray*}
and
\begin{eqnarray*}
& & \operatorname{cum}(J_{j_1,m_1;k_1},\overline{J_{j_3,m_4;k_1}})\operatorname{cum}(\overline
{J_{j_1,m_2;k_1}},J_{j_3,m_3;k_1})
\\
& &\qquad {} \times \operatorname{cum}(\overline{J_{j_2,m_1;k_2}},J_{j_4,m_4;k_2})
\operatorname{cum}(J_{j_2,m_2;k_2},\overline{J_{j_4,m_3;k_2}})
\\
&&\quad = f_{j_1,j_3,m_1,m_4;k_1}\overline{f_{j_1,j_3,m_2,m_3;k_1}}\overline {f_{j_2,j_4,m_1,m_4;k_2}}f_{j_2,j_4,m_2,m_3;k_2}+\mathrm{O}
\biggl(\frac
{1}{n} \biggr).
\end{eqnarray*}
Now, by taking the sums over $m_1,m_2,m_3,m_4$ of both expressions, we get
%
\begin{eqnarray}\label{eq12}
\tr\bigl(\mathbf{F}_{j_1j_1}(\omega_{k_1})\overline{\mathbf
{F}_{j_2j_2}(\omega_{k_2})}^T\bigr) \tr\bigl(
\overline{\mathbf{F}_{j_3j_3}(\omega_{k_1})}\mathbf
{F}_{j_4j_4}(\omega_{k_2})^T\bigr)+\mathrm{O} \biggl(
\frac{1}{n} \biggr)
\end{eqnarray}
and
\begin{eqnarray*}
\tr\bigl(\mathbf{F}_{j_1j_3}(\omega_{k_1})\overline{\mathbf
{F}_{j_2j_4}(\omega_{k_2})}^T\bigr) \tr\bigl(
\overline{\mathbf{F}_{j_1j_3}(\omega_{k_1})} \mathbf
{F}_{j_2j_4}(\omega_{k_2})^T\bigr)+\mathrm{O} \biggl(
\frac{1}{n} \biggr),
\end{eqnarray*}
respectively. Asymptotically equivalent to \eqref{eq2}, this results in
\begin{eqnarray*}
& & B_K\int_{-\ppi}^\ppi
\frac{1}{q^2}\sum_{j_1,j_2,j_3,j_4=1}^q \biggl(-1+q
\delta _{j_1j_3}\delta_{j_2j_4}+\frac{q}{q-1}(1-
\delta_{j_1j_3}) (1-\delta _{j_2j_4}) \biggr)
\\
& &\hphantom{B_K\int_{-\ppi}^\ppi
\frac{1}{q^2}\sum_{j_1,j_2,j_3,j_4=1}^q}{} \times \bigl\{\tr\bigl(\mathbf{F}_{j_1j_1}(\omega_{k_1})
\overline {\mathbf{F}_{j_2j_2}(\omega_{k_2})}^T\bigr)
\tr\bigl(\mathbf{F}_{j_3j_3}(\omega_{k_1})\overline{\mathbf
{F}_{j_4j_4}(\omega_{k_2})}^T\bigr)\\
&&\hphantom{B_K\int_{-\ppi}^\ppi
\frac{1}{q^2}\sum_{j_1,j_2,j_3,j_4=1}^q{} \times \bigl\{}{} +\bigl|\tr\bigl(
\mathbf{F}_{j_1j_3}(\omega)\overline{\mathbf {F}_{j_2j_4}(\omega)}
\bigr)\bigr|^2 \bigr\}\,\mathrm{d}\omega,
\end{eqnarray*}
which shows the first assertion of \eqref{tau_var}. For its second
part, we have
%
\begin{eqnarray}\label{eq5}
\Var\bigl(\tau_n^{*2}\bigr) &=& \frac{16h^2}{n^4}\mathop{
\mathop{\sum }_{k_1,k_2}}_{|k_1|\neq|k_2|} \mathop{\mathop{\sum
}_{k_3,k_4}}_{|k_3|\neq|k_4|} \biggl(\int_{-\ppi
}^\ppi
K_h(\omega-\omega _{k_1})K_h(\omega-
\omega_{k_2})\,\mathrm{d}\omega \biggr)^2 \nonumber
\\
& &\hphantom{\frac{16h^2}{n^4}\mathop{
\mathop{\sum }_{k_1,k_2}}_{|k_1|\neq|k_2|} \mathop{\mathop{\sum
}_{k_3,k_4}}_{|k_3|\neq|k_4|}} {}\times \biggl(\int_{-\ppi}^\ppi K_h(
\lambda-\omega _{k_3})K_h(\lambda-\omega_{k_4})\,\mathrm{d}
\omega \biggr)^2
\nonumber
\\
& &\hphantom{\frac{16h^2}{n^4}\mathop{
\mathop{\sum }_{k_1,k_2}}_{|k_1|\neq|k_2|} \mathop{\mathop{\sum
}_{k_3,k_4}}_{|k_3|\neq|k_4|}} {} \times\sum_{m_1,\ldots
,m_8=0}^{p-1}\frac{1}{q^4}\sum_{j_1,\ldots,j_8=1}^q
\biggl(-1+q\delta _{j_1j_3}\delta_{j_2j_4}\nonumber\\
& &\hphantom{\frac{16h^2}{n^4}\mathop{
\mathop{\sum }_{k_1,k_2}}_{|k_1|\neq|k_2|} \mathop{\mathop{\sum
}_{k_3,k_4}}_{|k_3|\neq|k_4|}{} \times\sum_{m_1,\ldots
,m_8=0}^{p-1}\frac{1}{q^4}\sum_{j_1,\ldots,j_8=1}^q\biggl(} {}+\frac{q}{q-1}(1-
\delta_{j_1j_3}) (1-\delta _{j_2j_4}) \biggr)
\nonumber
\\
& & \hphantom{\frac{16h^2}{n^4}\mathop{
\mathop{\sum }_{k_1,k_2}}_{|k_1|\neq|k_2|} \mathop{\mathop{\sum
}_{k_3,k_4}}_{|k_3|\neq|k_4|}{} \times\sum_{m_1,\ldots
,m_8=0}^{p-1}\frac{1}{q^4}\sum_{j_1,\ldots,j_8=1}^q} {}\times \biggl(-1+q\delta_{j_5j_7}\delta_{j_6j_8}\nonumber\\
&&\hphantom{\frac{16h^2}{n^4}\mathop{
\mathop{\sum }_{k_1,k_2}}_{|k_1|\neq|k_2|} \mathop{\mathop{\sum
}_{k_3,k_4}}_{|k_3|\neq|k_4|}{} \times\sum_{m_1,\ldots
,m_8=0}^{p-1}\frac{1}{q^4}\sum_{j_1,\ldots,j_8=1}^q{}\times\biggl(} {}+
\frac
{q}{q-1}(1-\delta_{j_5j_7}) (1-\delta_{j_6j_8}) \biggr)
\nonumber
\\
& & \hphantom{\frac{16h^2}{n^4}\mathop{
\mathop{\sum }_{k_1,k_2}}_{|k_1|\neq|k_2|} \mathop{\mathop{\sum
}_{k_3,k_4}}_{|k_3|\neq|k_4|}{} \times\sum_{m_1,\ldots
,m_8=0}^{p-1}\frac{1}{q^4}\sum_{j_1,\ldots,j_8=1}^q} {}\times \bigl\{E(J_{j_1,m_1;k_1}\overline{J_{j_1,m_2;k_1}}
J_{j_3,m_3;k_1}\nonumber\\
& & \hphantom{\frac{16h^2}{n^4}\mathop{
\mathop{\sum }_{k_1,k_2}}_{|k_1|\neq|k_2|} \mathop{\mathop{\sum
}_{k_3,k_4}}_{|k_3|\neq|k_4|}{} \times\sum_{m_1,\ldots
,m_8=0}^{p-1}\frac{1}{q^4}\sum_{j_1,\ldots,j_8=1}^q {}\times \bigl\{E}{}\times\overline{J_{j_3,m_4;k_1}}\overline {J_{j_2,m_1;k_2}}J_{j_2,m_2;k_2}\nonumber\\
& & \hphantom{\frac{16h^2}{n^4}\mathop{
\mathop{\sum }_{k_1,k_2}}_{|k_1|\neq|k_2|} \mathop{\mathop{\sum
}_{k_3,k_4}}_{|k_3|\neq|k_4|}{} \times\sum_{m_1,\ldots
,m_8=0}^{p-1}\frac{1}{q^4}\sum_{j_1,\ldots,j_8=1}^q{}\times \bigl\{E}{}\times \overline {J_{j_4,m_3;k_2}}J_{j_4,m_4;k_2} J_{j_5,m_5;k_3}
\nonumber
\\
& & \hphantom{\frac{16h^2}{n^4}\mathop{
\mathop{\sum }_{k_1,k_2}}_{|k_1|\neq|k_2|} \mathop{\mathop{\sum
}_{k_3,k_4}}_{|k_3|\neq|k_4|}{} \times\sum_{m_1,\ldots
,m_8=0}^{p-1}\frac{1}{q^4}\sum_{j_1,\ldots,j_8=1}^q{}\times \bigl\{E}{}\times \overline{J_{j_5,m_6;k_3}} J_{j_7,m_7;k_3}\\
& & \hphantom{\frac{16h^2}{n^4}\mathop{
\mathop{\sum }_{k_1,k_2}}_{|k_1|\neq|k_2|} \mathop{\mathop{\sum
}_{k_3,k_4}}_{|k_3|\neq|k_4|}{} \times\sum_{m_1,\ldots
,m_8=0}^{p-1}\frac{1}{q^4}\sum_{j_1,\ldots,j_8=1}^q{}\times \bigl\{E}{}\times \overline{J_{j_7,m_8;k_3}}\overline {J_{j_6,m_5;k_4}}J_{j_6,m_6;k_4}\nonumber\\
& & \hphantom{\frac{16h^2}{n^4}\mathop{
\mathop{\sum }_{k_1,k_2}}_{|k_1|\neq|k_2|} \mathop{\mathop{\sum
}_{k_3,k_4}}_{|k_3|\neq|k_4|}{} \times\sum_{m_1,\ldots
,m_8=0}^{p-1}\frac{1}{q^4}\sum_{j_1,\ldots,j_8=1}^q{}\times \bigl\{E}{}\times \overline {J_{j_8,m_7;k_4}}J_{j_8,m_8;k_4})
\nonumber
\\
& & \hphantom{\frac{16h^2}{n^4}\mathop{
\mathop{\sum }_{k_1,k_2}}_{|k_1|\neq|k_2|} \mathop{\mathop{\sum
}_{k_3,k_4}}_{|k_3|\neq|k_4|}{} \times\sum_{m_1,\ldots
,m_8=0}^{p-1}\frac{1}{q^4}\sum_{j_1,\ldots,j_8=1}^q{}\times \bigl\{}{}  -E(J_{j_1,m_1;k_1}\overline{J_{j_1,m_2;k_1}} J_{j_3,m_3;k_1}\nonumber\\
& & \hphantom{\frac{16h^2}{n^4}\mathop{
\mathop{\sum }_{k_1,k_2}}_{|k_1|\neq|k_2|} \mathop{\mathop{\sum
}_{k_3,k_4}}_{|k_3|\neq|k_4|}{} \times\sum_{m_1,\ldots
,m_8=0}^{p-1}\frac{1}{q^4}\sum_{j_1,\ldots,j_8=1}^q{}\times \bigl\{{}-E }{}\times
\overline{J_{j_3,m_4;k_1}}\overline {J_{j_2,m_1;k_2}}\nonumber\\
& & \hphantom{\frac{16h^2}{n^4}\mathop{
\mathop{\sum }_{k_1,k_2}}_{|k_1|\neq|k_2|} \mathop{\mathop{\sum
}_{k_3,k_4}}_{|k_3|\neq|k_4|}{} \times\sum_{m_1,\ldots
,m_8=0}^{p-1}\frac{1}{q^4}\sum_{j_1,\ldots,j_8=1}^q{}\times \bigl\{{}-E }{}\times J_{j_2,m_2;k_2}\overline {J_{j_4,m_3;k_2}}\nonumber\\
& & \hphantom{\frac{16h^2}{n^4}\mathop{
\mathop{\sum }_{k_1,k_2}}_{|k_1|\neq|k_2|} \mathop{\mathop{\sum
}_{k_3,k_4}}_{|k_3|\neq|k_4|}{} \times\sum_{m_1,\ldots
,m_8=0}^{p-1}\frac{1}{q^4}\sum_{j_1,\ldots,j_8=1}^q{}\times \bigl\{{}-E }{}\times J_{j_4,m_4;k_2})
\nonumber
\\
& & \hphantom{\frac{16h^2}{n^4}\mathop{
\mathop{\sum }_{k_1,k_2}}_{|k_1|\neq|k_2|} \mathop{\mathop{\sum
}_{k_3,k_4}}_{|k_3|\neq|k_4|}{} \times\sum_{m_1,\ldots
,m_8=0}^{p-1}\frac{1}{q^4}\sum_{j_1,\ldots,j_8=1}^q{}\times \bigl\{{}-}{} \times E(J_{j_5,m_5;k_3}\overline{J_{j_5,m_6;k_3}}\nonumber\\
& & \hphantom{\frac{16h^2}{n^4}\mathop{
\mathop{\sum }_{k_1,k_2}}_{|k_1|\neq|k_2|} \mathop{\mathop{\sum
}_{k_3,k_4}}_{|k_3|\neq|k_4|}{} \times\sum_{m_1,\ldots
,m_8=0}^{p-1}\frac{1}{q^4}\sum_{j_1,\ldots,j_8=1}^q{}\times \bigl\{-{}\times E }{}\times
J_{j_7,m_7;k_3}\overline{J_{j_7,m_8;k_3}}\nonumber\\
& & \hphantom{\frac{16h^2}{n^4}\mathop{
\mathop{\sum }_{k_1,k_2}}_{|k_1|\neq|k_2|} \mathop{\mathop{\sum
}_{k_3,k_4}}_{|k_3|\neq|k_4|}{} \times\sum_{m_1,\ldots
,m_8=0}^{p-1}\frac{1}{q^4}\sum_{j_1,\ldots,j_8=1}^q{}\times \bigl\{-{}\times E }{}\times
\overline {J_{j_6,m_5;k_4}}J_{j_6,m_6;k_4}\nonumber\\
& & \hphantom{\frac{16h^2}{n^4}\mathop{
\mathop{\sum }_{k_1,k_2}}_{|k_1|\neq|k_2|} \mathop{\mathop{\sum
}_{k_3,k_4}}_{|k_3|\neq|k_4|}{} \times\sum_{m_1,\ldots
,m_8=0}^{p-1}\frac{1}{q^4}\sum_{j_1,\ldots,j_8=1}^q{}\times \bigl\{-{}\times E }{}\times
\overline {J_{j_8,m_7;k_4}}J_{j_8,m_8;k_4}) \bigr\}
\nonumber
\end{eqnarray}
and similar to the computations for $E(\tau_n^{*2})$, we are able to
express the $16$th and both $8$th moments with cumulants and the
difference above
becomes
%
\begin{eqnarray}\label{blub2}
\bigl\{E(J_1\cdots J_{16})-E(J_1\cdots
J_8)E(J_9\cdots J_{16}) \bigr\} =\sum
_{\widetilde\sigma}\prod_{B\in\widetilde{\sigma}}
\operatorname{cum}(J_i\dvt i\in B),
\end{eqnarray}
where we identify the DFTs $J_{l,m;k}$ from above with $J_1,\ldots
,J_{16}$ for notational convenience and now $\widetilde\sigma$ runs
through the list of all partitions of $\{1,\ldots,16\}$ that can not
be written as cumulants depending on subsets of $\{J_1,\ldots,J_8\}$
multiplied with cumulants depending on subsets of $\{J_9,\ldots
,J_{16}\}$ and $B$ runs through the list of all blocks of the partition
$\widetilde\sigma$. Now, we can use Lemma~\ref{lemma2} and insert
\eqref{brillinger} in \eqref{blub2}. Again the largest contribution
comes from products of eight cumulants of second order and due to
$\sum_{t=1}^{n}\mathrm{e}^{\mathrm{i}t\sum_{j=1}^{k}\omega_{j}}=n$ if $\sum_{j=1}^{k}\omega_{j}=0$ and $\sum_{t=1}^{n}\mathrm{e}^{\mathrm{i}t\sum_{j=1}^{k}\omega
_{j}}=0$ otherwise, all these combinations that occur in \eqref{blub2}
contain a factor $\delta_{|k_i|,|k_j|}$ for $i\neq j$. This is
either excluded in the summation in \eqref{eq5} or results in at least
one collapsing sum, which eventually implies $\Var(\tau_n^{*2})=\mathrm{o}(1)$.
\end{pf}

\begin{pf*}{Proof of Remark \protect\ref{H0remark}}
We only treat the variance. Under $H_0$, the second summand in \eqref
{eq62} does not depend on $j_1,\dots,j_4$, so that this case vanishes
thanks to
\begin{eqnarray*}
\sum_{j_1,j_2,j_3,j_4=1}^q \biggl(-1+q
\delta_{j_1j_3}\delta _{j_2j_4}+\frac{q}{q-1}(1-
\delta_{j_1j_3}) (1-\delta_{j_2j_4}) \biggr)=0,
\end{eqnarray*}
which yields the desired assertion. \end{pf*}

\begin{pf*}{Proof of Theorem \protect\ref{theorem4}}
Note first that by Theorem~\ref{theorem3} and the $\sqrt {nh}$-consistency of the spectral density estimators, see the proof of
Theorem~\ref{theorem11}, we have
$\mu_n^* - \widehat{\mu}^* = \mathrm{o}_P(h^{1/2})$. Hence, it is sufficient
to prove a central limit theorem for $T_n^*-h^{-1/2}\mu_n^*$.

By using Lemma~\ref{lemma1}(i), we have
\begin{eqnarray*}
T_n^* &=& \frac{h^{\sfrac{1}{2}}}{n} \sum_{k_1,k_2}
\int_{-\ppi}^\ppi K_h(\omega-
\omega_{k_1})K_h(\omega-\omega_{k_2})\,\mathrm{d}\omega
\\
& &\hphantom{\frac{h^{\sfrac{1}{2}}}{n} \sum_{k_1,k_2}}{} \times\sum
_{m_1,m_2=0}^{p-1}\sum
_{j_1,j_2=1}^q \Biggl(-\frac{1}{q}+\sum_{r=1}^q
\delta_{j_1,\pii
_{k_1}(r)}\delta_{j_2,\pii_{k_2}(r)} \Biggr) I_{j_1,m_1,m_2;k_1}
\overline{I_{j_2,m_1,m_2;k_2}}
\\
&=:& \sum_{k_1,k_2} w_{n,k_1,k_2}^*(
\pii_{k_1},\pii_{k_2})
\end{eqnarray*}
with an obvious notation for $w_{n,k_1,k_2}^*(\pii_{k_1},\pii_{k_2})=:
w_{k_1,k_2}^*(\pii_{k_1},\pii_{k_2})$. Because of\vspace*{-1pt}
%
\begin{eqnarray}
\label{clean condition} E^*\Biggl(\sum_{r=1}^q
\delta_{\pii_{k_1}(r),j_1}\delta_{\pii
_{k_2}(r),j_2}\Bigl|\pii_{k_1}\Biggr) =
\cases{ q^{-1}, &\quad  $|k_1|\neq|k_2|$,
\cr
\delta_{j_1j_2}, & \quad $|k_1|=|k_2|$,} %
\end{eqnarray}
it follows that\vspace*{-1pt}
\begin{eqnarray*}
T_n^* &=& \sum_{|k_1|\neq|k_2]}
w_{k_1,k_2}^*(\pii_{k_1},\pii_{k_2}) + \sum
_{|k_1|= |k_2]} w_{k_1,k_2}^*(\pii_{k_1},
\pii_{k_2})
\\
&=:& W_n^* + \sum_{|k_1|= |k_2]}
w_{k_1,k_2}^*(\pii_{k_1},\pii_{k_2}) = W_n^*
+ h^{-1/2}\mu_n^* + \mathrm{o}_{P^*}(1),
\end{eqnarray*}
see Theorem~\ref{theorem2} above. Setting $\widetilde
{W}_{k_1,k_2}^*:= w_{k_1,k_2}^*(\pii_{k_1},\pii_{k_2}) +
w_{k_2,k_1}^*(\pii_{k_2},\pii_{k_1})$, we obtain up to a negligible term
\begin{eqnarray*}
W_n^* &=& \sum_{0\leq k_1<k_2 \leq\lfloor n/2 \rfloor} \widetilde
W_{k_1,k_2}^*+\widetilde W_{-k_1,k_2}^* +\widetilde
W_{k_1,-k_2}^*+\widetilde W_{-k_1,-k_2}^*
\\
&=:& \sum_{0\leq k_1<k_2 \leq\lfloor n/2 \rfloor} W_{k_1,k_2}^*.
\end{eqnarray*}
Since the random variables $W_{k_1,k_2}^*$ are clean by \eqref{clean
condition} in the sense of Definition~2.1\vspace*{2pt} of De Jong \cite{Dejong:1987}, that is,
$E^*(W_{k_1,k_2}^*|\ppi_{k_1})= 0$ holds a.s., we are in the situation
to apply\vspace*{2pt} Proposition~3.2. of his paper. Hence, for obtaining
convergence in distribution (conditioned on the data)
%
\begin{eqnarray}
\label{conditionalCLTstud} \frac{1}{\tau_n^{*}} W_n^* \schwach N(0,1)
\end{eqnarray}
in probability, it remains to prove that
%
\begin{eqnarray}
\label{G1} G_1&:=&\sum_{0\leq k_1<k_2\leq\lfloor n/2\rfloor} E^*
\bigl(W_{k_1,k_2}^{*4}\bigr),
\\
\label{G2} G_2 &:=& \sum_{0\leq k_1<k_2<k_3\leq\lfloor n/2\rfloor} E^*
\bigl(W_{k_1,k_2}^{*2}W_{k_1,k_3}^{*2}+W_{k_2,k_1}^{*2}W_{k_2,k_3}^{*2}+W_{k_3,k_1}^{*2}W_{k_3,k_2}^{*2}
\bigr)
\end{eqnarray}
and
%
\begin{eqnarray}
\label{G4} G_4& :=&  \sum_{0\leq k_1<k_2<k_3<k_4\leq\lfloor n/2\rfloor} E^*
\bigl(W_{k_1,k_2}^*W_{k_1,k_3}^*W_{k_4,k_2}^*W_{k_4,k_3}^*+ W_{k_1,k_2}^*W_{k_1,k_4}^*W_{k_3,k_2}^*W_{k_3,k_4}^*\nonumber
\\[-8pt]\\[-8pt]
&&\hphantom{\sum_{0\leq k_1<k_2<k_3<k_4\leq\lfloor n/2\rfloor} E^*
\bigl(}{}
+W_{k_1,k_3}^*W_{k_4,k_1}^*W_{k_2,k_3}^*W_{k_2,k_4}^*
\bigr)
\nonumber
\end{eqnarray}
are all of lower order (in probability) than $\tau_n^{* 2}$, that is,
are all negligible. We prove the requested result only for the most
contributing case (\ref{G4}). Because of the inherited symmetries and
the constant number of periodogram factors with frequencies $\omega
_{k_1},\omega_{k_2},\omega_{k_3}$ and $\omega_{k_4}$, it suffices to
consider the representative
\begin{eqnarray*}
\widetilde G_4 &:=& \sum_{0\leq k_1<k_2<k_3<k_4\leq\lfloor n/2\rfloor} E^*
\bigl( w_{k_1,k_2}^*(\pii_{k_1},\pii_{k_2})w_{k_1,k_3}^*(
\pii_{k_1},\pii _{k_3})w_{k_4,k_2}^* (\pii_{k_4},
\pii_{k_2})w_{k_4,k_3}^*(\pii_{k_4},\pii_{k_3})
\bigr)
\\
&=& \sum_{0\leq k_1<k_2<k_3<k_4\leq\lfloor n/2\rfloor} \frac
{h^2}{n^4} \sum
_{m_1,\dots,m_8=0}^{p-1}\sum_{j_1,\dots,j_8=1}^q
\int_{-\ppi}^\ppi K_h(
\lambda_1-\omega_{k_1})K_h(\lambda
_1-\omega_{k_2})\,\mathrm{d}\lambda_1 \\
&&\hspace*{-3pt}\hphantom{\sum_{0\leq k_1<k_2<k_3<k_4\leq\lfloor n/2\rfloor} \frac
{h^2}{n^4} \sum
_{m_1,\dots,m_8=0}^{p-1}\sum_{j_1,\dots,j_8=1}^q}{}\times \int
_{-\ppi}^\ppi K_h(\lambda_2-
\omega_{k_1})K_h(\lambda_2-\omega
_{k_3})\,\mathrm{d}\lambda_2
\\
&& \hspace*{-3pt}\hphantom{\sum_{0\leq k_1<k_2<k_3<k_4\leq\lfloor n/2\rfloor} \frac
{h^2}{n^4} \sum
_{m_1,\dots,m_8=0}^{p-1}\sum_{j_1,\dots,j_8=1}^q}{}\times\int_{-\ppi}^\ppi K_h(
\lambda_3-\omega_{k_4})K_h(\lambda
_3-\omega_{k_2})\,\mathrm{d}\lambda_3\\
&&\hspace*{-3pt}\hphantom{\sum_{0\leq k_1<k_2<k_3<k_4\leq\lfloor n/2\rfloor} \frac
{h^2}{n^4} \sum
_{m_1,\dots,m_8=0}^{p-1}\sum_{j_1,\dots,j_8=1}^q}{}\times \int
_{-\ppi}^\ppi K_h(\lambda_4-
\omega_{k_4})K_h(\lambda_4-\omega
_{k_3})\,\mathrm{d}\lambda_4
\\
&&\hspace*{-3pt}\hphantom{\sum_{0\leq k_1<k_2<k_3<k_4\leq\lfloor n/2\rfloor} \frac
{h^2}{n^4} \sum
_{m_1,\dots,m_8=0}^{p-1}\sum_{j_1,\dots,j_8=1}^q}{} \times I_{j_1,m_1,m_2;k_1}\overline{I_{j_2,m_1,m_2;k_2}}I_{j_3,m_3,m_4;k_1}\\
&&\hspace*{-3pt}\hphantom{\sum_{0\leq k_1<k_2<k_3<k_4\leq\lfloor n/2\rfloor} \frac
{h^2}{n^4} \sum
_{m_1,\dots,m_8=0}^{p-1}\sum_{j_1,\dots,j_8=1}^q}{} \times
\overline{I_{j_4,m_3,m_4;k_3}}  I_{j_5,m_5,m_6;k_4}\overline{I_{j_6,m_5,m_6;k_2}}\\
&&\hspace*{-3pt}\hphantom{\sum_{0\leq k_1<k_2<k_3<k_4\leq\lfloor n/2\rfloor} \frac
{h^2}{n^4} \sum
_{m_1,\dots,m_8=0}^{p-1}\sum_{j_1,\dots,j_8=1}^q}{} \times I_{j_7,m_7,m_8;k_4}
\overline{I_{j_8,m_7,m_8;k_3}}
\\
&&\hspace*{-3pt} \hphantom{\sum_{0\leq k_1<k_2<k_3<k_4\leq\lfloor n/2\rfloor} \frac
{h^2}{n^4} \sum
_{m_1,\dots,m_8=0}^{p-1}\sum_{j_1,\dots,j_8=1}^q}{}\times E^* \Biggl(\sum_{r_1,\dots,r_4=1}^q \bigl(
\delta_{\pii
_{k_1}(r_1),j_1}\delta_{\pii_{k_2}(r_1),j_2} -q^{-1}\bigr)\cdots\\
&&\hspace*{-8pt} \hphantom{\sum_{0\leq k_1<k_2<k_3<k_4\leq\lfloor n/2\rfloor} \frac
{h^2}{n^4} \sum
_{m_1,\dots,m_8=0}^{p-1}\sum_{j_1,\dots,j_8=1}^q\times E^* \Biggl(\sum_{r_1,\dots,r_4=1}^q}{}\times\bigl(\delta_{\pii_{k_4}(r_4),j_7}\delta_{\pii
_{k_3}(r_4),j_8}
-q^{-1}\bigr)
\Biggr).
\end{eqnarray*}
Since the last conditional expectation above is bounded by $0$ from
below and by $q^4$ from above, taking the expectation above and
applying (\ref{brillinger}) for $k=16$ results in
\begin{eqnarray*}
E(\widetilde G_4) &=& \Biggl(\sum_{0\leq k_1<k_2<k_3<k_4\leq\lfloor
n/2\rfloor}
\frac{h^2}{n^4} \sum_{m_1,\dots,m_8=0}^{p-1}\sum
_{j_1,\dots,j_8=1}^q
\int_{-\ppi}^\ppi K_h(
\lambda_1-\omega_{k_1})K_h(\lambda
_1-\omega_{k_2})\,\mathrm{d}\lambda_1\\
&& \hphantom{\Biggl(\sum_{0\leq k_1<k_2<k_3<k_4\leq\lfloor
n/2\rfloor}
\frac{h^2}{n^4} \sum_{m_1,\dots,m_8=0}^{p-1}\sum
_{j_1,\dots,j_8=1}^q}{}\times \int
_{-\ppi}^\ppi K_h(\lambda_2-
\omega_{k_1})K_h(\lambda_2-\omega
_{k_3})\,\mathrm{d}\lambda_2
\\
&&\hphantom{\Biggl(\sum_{0\leq k_1<k_2<k_3<k_4\leq\lfloor
n/2\rfloor}
\frac{h^2}{n^4} \sum_{m_1,\dots,m_8=0}^{p-1}\sum
_{j_1,\dots,j_8=1}^q}{} \times\int_{-\ppi}^\ppi K_h(
\lambda_3-\omega _{k_4})K_h(
\lambda_3-\omega_{k_2})\,\mathrm{d}\lambda_3 \\
&&\hphantom{\Biggl(\sum_{0\leq k_1<k_2<k_3<k_4\leq\lfloor
n/2\rfloor}
\frac{h^2}{n^4} \sum_{m_1,\dots,m_8=0}^{p-1}\sum
_{j_1,\dots,j_8=1}^q}{}\times\int
_{-\ppi}^\ppi K_h(\lambda_4-
\omega_{k_4})K_h(\lambda_4-\omega
_{k_3})\,\mathrm{d}\lambda_4 \Biggr)\\
&&{}\times\mathrm{O}(1).
\end{eqnarray*}
Approximating all Riemann sums by their corresponding integrals and by
using standard substitutions, the expression above becomes an $\mathrm{O}(h)$
term. Similar arguments and using (\ref{brillinger}) for $k=32$ yield
$\Var(\widetilde G_4)=\mathrm{O}(h^2)$. This completes the proof.
\end{pf*}

\begin{pf*}{Proof of Corollary \protect\ref{corollary1}}
By comparing the results in Theorem~\ref{theorem1} and Theorem~\ref
{theorem4}, we have asymptotic exactness of $\varphi_{n,\mathrm{cent}}^*$ if
and only if $\tau^2_0=\tau_0^{* 2}$ holds under $H_0$. Rearrangements
of the summations over $j_1,j_2,j_3,j_4$ in the integrands of \eqref
{eq7} and \eqref{eq6} yield to the condition that
%
\begin{eqnarray}\label{eq8}
& & (q-1)\bigl|\tr\bigl(\mathbf{F}_{11}(\omega)\overline{\mathbf
{F}_{11}(\omega)}^T\bigr)\bigr|^2 -
\frac{2}{q}\mathop{\mathop{\sum}_{j_2,j_4=1}}_{j_2\neq j_4}^q\bigl|
\tr \bigl(\mathbf {F}_{11}(\omega)\overline{\mathbf{F}_{j_2j_4}(
\omega)}^T\bigr)\bigr|^2
\nonumber
\\[-8pt]\\[-8pt]
& &\quad {} +\mathop{\mathop{\sum}_{j_1,j_2,j_3,j_4=1}}_{j_1\neq j_3,j_2\neq
j_4}^q
\bigl(-q+q^2\delta_{j_1j_2}\bigr) \bigl(-q+q^2
\delta_{j_3j_4}\bigr) \bigl|\tr\bigl(\mathbf{F}_{j_1j_3}(\omega)
\overline{\mathbf {F}_{j_2j_4}(\omega)}^T\bigr)\bigr|^2
\nonumber
\end{eqnarray}
and
%
\begin{eqnarray}\label{eq9}
& & (q-1)\bigl|\tr\bigl(\mathbf{F}_{11}(\omega)\overline{\mathbf
{F}_{11}(\omega)}^T\bigr)\bigr|^2-\frac{2}{q}
\mathop{\mathop{\sum }_{j_2,j_4=1}}_{j_2\neq j_4}^q\bigl|\tr
\bigl(\mathbf{F}_{11}(\omega)\overline {\mathbf {F}_{j_2j_4}(
\omega)}^T\bigr)\bigr|^2
\nonumber
\\[-8pt]\\[-8pt]
& &\quad{}  +\mathop{\mathop{\sum}_{j_1,j_2,j_3,j_4=1}}_{j_1\neq j_3,j_2\neq
j_4}^q
\frac{1}{q^2(q-1)}\bigl|\tr\bigl(\mathbf {F}_{j_1j_3}(\omega)\overline{
\mathbf{F}_{j_2j_4}(\omega)}^T\bigr)\bigr|^2,
\nonumber
\end{eqnarray}
have to be equal. Equalizing both quantities and another rearrangement
of both last sums in \eqref{eq8} and in \eqref{eq9} gives the desired
result. By Lemma~1 in Janssen and Pauls \cite
{janssen/pauls:2003a}, this even shows the
asymptotic equivalence of the tests in \eqref{asequiv}.
Since we have $\widehat{\mu}-\mu_n^*=\mathrm{o}_P(h^{1/2})$ under $H_0$, the
same holds true for the other version $\varphi_n^*$.

Consistency for $\varphi_{n,\mathrm{cent}}^*$ follows as in the proof of
Theorem~\ref{theorem1} since $T_n-h^{-1/2}\widehat{\mu}$ converges
in probability to $+\infty$ under $H_1$ but the
critical value still converges in probability to a quantile of a
normal distribution. For the other test, we can rewrite $\varphi_n^*$
asymptotically equivalent as $\mathbf{1}_{(c_n^*(\alpha),\infty
)}(T_n-h^{-1/2}\mu_n^*)$ and consistency follows as above since
$T_n-h^{-1/2}\mu_n^*$ also converges in probability to $+\infty$
under $H_1$.\end{pf*}

\begin{pf*}{Proof of Corollary \protect\ref{corollary2}}
For $q=2$ the ratios occurring in the first and in the
third sum in \eqref{eq11} are zero and the second sum does not occur
at all. In the cases (ii) and (iii), we can treat $|\tr(\mathbf
{F}_{j_1j_3}(\omega)\overline{\mathbf{F}_{j_2j_4}(\omega)}^T)|^2$
as a constant, which causes the three sums in \eqref{eq11} to cumulate
to zero. For (iv) the right-hand side of Corollary~\ref
{corollary1}(ii) becomes
\begin{eqnarray*}
& & \frac{8}{9} \bigl(\bigl|f_{12}(\omega)\bigr|^4+\bigl|f_{13}(
\omega )\bigr|^4+\bigl|f_{23}(\omega)\bigr|^4 \bigr)
\\
& &\qquad {} -\frac{8}{9} \bigl(\bigl|f_{12}(\omega)\bigr|^2\bigl|f_{13}(
\omega )\bigr|^2+\bigl|f_{12}(\omega)\bigr|^2\bigl|f_{23}(
\omega)\bigr|^2+\bigl|f_{13}(\omega )\bigr|^2\bigl|f_{23}(
\omega)\bigr|^2 \bigr)
\\
&&\quad = \frac{4}{9} \bigl\{\bigl(\bigl|f_{12}(\omega)\bigr|^2-\bigl|f_{13}(
\omega )\bigr|^2\bigr)^2+\bigl(\bigl|f_{12}(
\omega)\bigr|^2-\bigl|f_{23}(\omega)\bigr|^2
\bigr)^2+\bigl(\bigl|f_{13}(\omega )\bigr|^2-\bigl|f_{23}(
\omega)\bigr|^2\bigr)^2 \bigr\},
\end{eqnarray*}
which gives the desired result.\end{pf*}

\begin{pf*}{Proof of Corollary \protect\ref{corollary3}}
Since $\widehat{\tau}^{* 2}$ converges in $\mathbb{P}\otimes
\widetilde{\mathbb{P}}$ probability to $\tau^{* 2}$, the critical
value $c_{n,\mathrm{stud}}^*(\alpha)$ always converges in probability to
$u_{1-\alpha}$. Hence, we can close the proof as in Theorem~\ref{theorem1}.
\end{pf*}

\begin{pf*}{Proof of Remark \protect\ref{linear processes}}
Note that the assumptions guarantee that the conditional CLT, see
Theorem~\ref{theorem4}, holds true as above. Moreover, the test
statistic also satisfies the CLT stated in Theorem~\ref{theorem1}
under the null, see Equations (2.4)--(2.6) in Dette and Paparoditis \cite{Dette:2009}, which shows the asymptotic exactness
under the null. Finally,
consistency of the tests follows from Theorem~2 of their paper.
\end{pf*}

\begin{pf*}{Proof of Corollary \protect\ref{corollary4}}
Checking through the proof of Theorem~\ref{theorem4} (where the
central limit theorem (\ref{CLT Tn*}) holds under the null as well as
for fixed alternatives)
we see that all results even remain valid for the processes $\underline
{X}_1^n,\dots,\underline{X}_n^n$ with spectral density given by
$\mathbf{f}^n$, that is, we also have convergence in distribution
given the data
\[
\bigl(\widehat{\tau}^{*}\bigr)^{-1} \biggl(T_n^*-
\frac{ \widehat{\mu}^*}{\sqrt {h}} \biggr)\stackrel{\mathcal{D}} {\longrightarrow}\mathcal{N}(0,1)
\]
in probability in this case.
Since the test statistic is asymptotically normally distributed with
mean $\nu$ and variance $\tau_0^2$, see the proof of Theorem~\ref
{theorem: loclaAlternatives}, the result (i) follows as in the proof of
Corollary~\ref{corollary3}.

For proving (ii), recall that our mixing conditions imply the $\sqrt {nh}$-consistency of the kernel density estimators for the processes.
This shows that $\widehat{\mu}-\widehat{\mu}^*=\mathrm{o}_P(h^{1/2})$ even
holds under the given local alternatives (\ref{localAlternatives})
with $\alpha_n = h^{-1/4}n^{-1/2}$.
Since we also have $\widehat{\tau}-\widehat{\tau}^*=\mathrm{o}_P(h^{1/2})$,
the result follows from Slutsky's lemma as in the proof of Corollary~\ref{corollary1}.
\end{pf*}

\section*{Acknowledgements}
The authors are very grateful to the Associate Editor and two referees
for their careful and detailed comments that led to a considerable
improvement of the paper.
Parts of this paper have been developed while the second author was
supported by a fellowship within the Postdoc-Programme of the
German Academic Exchange Service (DAAD) carried out at the University
of Bern, Switzerland.

\begin{supplement}
\stitle{Supplement to ``Testing equality of spectral densities using randomization techniques''}
\slink[doi]{10.3150/13-BEJ584SUPP} 
\sdatatype{.pdf}
\sfilename{BEJ584\_supp.pdf}
\sdescription{In the supplement to the current paper (cf. Jentsch and Pauly \cite{Jentsch/Pauly:2013}),
we provide additional supporting simulations for the asymptotic test
and all three randomization tests under consideration in a variety of examples.}
\end{supplement}

%

\printhistory

\end{document}